\documentclass[9pt]{amsart}
\usepackage{amsmath}
\usepackage{amsxtra}
\usepackage{amscd}
\usepackage{amsthm}
\usepackage{amsfonts}
\usepackage{amssymb}
\usepackage[foot]{amsaddr}
\usepackage{cite}
\usepackage{url}
\usepackage{rotating}
\usepackage{eucal}
\usepackage[all,2cell]{xy}
\UseAllTwocells
\usepackage{graphicx}
\usepackage{verbatim}
\usepackage{hyperref}
\usepackage{xparse}
\usepackage{upgreek}
\usepackage{MnSymbol}
\usepackage[usenames,dvipsnames,svgnames,table]{xcolor}

\newtheorem{cor}[subsubsection]{Corollary}
\newtheorem{lem}[subsubsection]{Lemma}
\newtheorem{goal}[subsubsection]{Goal}
\newtheorem{lem-defn}[subsubsection]{Lemma-Definition}
\newtheorem{prop}[subsubsection]{Proposition}

\newtheorem{warn}[subsubsection]{Warning}

\newtheorem{facts}[subsubsection]{Fact}

\newtheorem{propconstr}[subsubsection]{Proposition-Construction}

\newtheorem{thm}[subsubsection]{Theorem}
\newtheorem{defn}[subsubsection]{Definition}
\newtheorem{notn}[subsubsection]{Notation}
\newtheorem{convn}[subsubsection]{Convension}
\newtheorem{constr}[subsubsection]{Construction}

\numberwithin{equation}{section}
\theoremstyle{remark}
\newtheorem{rem}[subsubsection]{Remark}
\newtheorem{exam}[subsubsection]{Example}

\newcommand\nc{\newcommand}
\nc\on{\operatorname}
\nc\renc{\renewcommand}

\nc\ssec{\subsection}
\nc\sssec{\subsubsection}

\nc\blongeqn{\[ \begin{aligned}}
\nc\elongeqn{\end{aligned} \]}
\nc\mBA{{\mathbb A}}
\nc\mBB{{\mathbb B}}
\nc\mBC{{\mathbb C}}
\nc\mBD{{\mathbb D}}
\nc\mBE{{\mathbb E}}
\nc\mBF{{\mathbb F}}
\nc\mBG{{\mathbb G}}
\nc\mBH{{\mathbb H}}
\nc\mBI{{\mathbb I}}
\nc\mBJ{{\mathbb J}}
\nc\mBK{{\mathbb K}}
\nc\mBL{{\mathbb L}}
\nc\mBM{{\mathbb M}}
\nc\mBN{{\mathbb N}}
\nc\mBO{{\mathbb O}}
\nc\mBP{{\mathbb P}}
\nc\mBQ{{\mathbb Q}}
\nc\mBR{{\mathbb R}}
\nc\mBS{{\mathbb S}}
\nc\mBT{{\mathbb T}}
\nc\mBU{{\mathbb U}}
\nc\mBV{{\mathbb V}}
\nc\mBW{{\mathbb W}}
\nc\mBX{{\mathbb X}}
\nc\mBY{{\mathbb Y}}
\nc\mBZ{{\mathbb Z}}

\nc\mCA{{\mathcal A}}
\nc\mCB{{\mathcal B}}
\nc\mCC{{\mathcal C}}
\nc\mCD{{\mathcal D}}
\nc\mCE{{\mathcal E}}
\nc\mCF{{\mathcal F}}
\nc\mCG{{\mathcal G}}
\nc\mCH{{\mathcal H}}
\nc\mCI{{\mathcal I}}
\nc\mCJ{{\mathcal J}}
\nc\mCK{{\mathcal K}}
\nc\mCL{{\mathcal L}}
\nc\mCM{{\mathcal M}}
\nc\mCN{{\mathcal N}}
\nc\mCO{{\mathcal O}}
\nc\mCP{{\mathcal P}}
\nc\mCQ{{\mathcal Q}}
\nc\mCR{{\mathcal R}}
\nc\mCS{{\mathcal S}}
\nc\mCT{{\mathcal T}}
\nc\mCU{{\mathcal U}}
\nc\mCV{{\mathcal V}}
\nc\mCW{{\mathcal W}}
\nc\mCX{{\mathcal X}}
\nc\mCY{{\mathcal Y}}
\nc\mCZ{{\mathcal Z}}

\nc\mbA{{\mathbf A}}
\nc\mbB{{\mathbf B}}
\nc\mbC{{\mathbf C}}
\nc\mbD{{\mathbf D}}
\nc\mbE{{\mathbf E}}
\nc\mbF{{\mathbf F}}
\nc\mbG{{\mathbf G}}
\nc\mbH{{\mathbf H}}
\nc\mbI{{\mathbf I}}
\nc\mbJ{{\mathbf J}}
\nc\mbK{{\mathbf K}}
\nc\mbL{{\mathbf L}}
\nc\mbM{{\mathbf M}}
\nc\mbN{{\mathbf N}}
\nc\mbO{{\mathbf O}}
\nc\mbP{{\mathbf P}}
\nc\mbQ{{\mathbf Q}}
\nc\mbR{{\mathbf R}}
\nc\mbS{{\mathbf S}}
\nc\mbT{{\mathbf T}}
\nc\mbU{{\mathbf U}}
\nc\mbV{{\mathbf V}}
\nc\mbW{{\mathbf W}}
\nc\mbX{{\mathbf X}}
\nc\mbY{{\mathbf Y}}
\nc\mbZ{{\mathbf Z}}

\nc\mba{{\mathbf a}}
\nc\mbb{{\mathbf b}}
\nc\mbc{{\mathbf c}}
\nc\mbd{{\mathbf d}}
\nc\mbe{{\mathbf e}}
\nc\mbf{{\mathbf f}}
\nc\mbg{{\mathbf g}}
\nc\mbh{{\mathbf h}}
\nc\mbi{{\mathbf i}}
\nc\mbj{{\mathbf j}}
\nc\mbk{{\mathbf k}}
\nc\mbl{{\mathbf l}}
\nc\mbm{{\mathbf m}}
\nc\mbn{{\mathbf n}}
\nc\mbo{{\mathbf o}}
\nc\mbp{{\mathbf p}}
\nc\mbq{{\mathbf q}}
\nc\mbr{{\mathbf r}}
\nc\mbs{{\mathbf s}}
\nc\mbt{{\mathbf t}}
\nc\mbu{{\mathbf u}}
\nc\mbv{{\mathbf v}}
\nc\mbw{{\mathbf w}}
\nc\mbx{{\mathbf x}}
\nc\mby{{\mathbf y}}
\nc\mbz{{\mathbf z}}

\nc\mfa{{\mathfrak a}}
\nc\mfb{{\mathfrak b}}
\nc\mfc{{\mathfrak c}}
\nc\mfd{{\mathfrak d}}
\nc\mfe{{\mathfrak e}}
\nc\mff{{\mathfrak f}}
\nc\mfg{{\mathfrak g}}
\nc\mfh{{\mathfrak h}}
\nc\mfi{{\mathfrak i}}
\nc\mfj{{\mathfrak j}}
\nc\mfk{{\mathfrak k}}
\nc\mfl{{\mathfrak l}}
\nc\mfm{{\mathfrak m}}
\nc\mfn{{\mathfrak n}}
\nc\mfo{{\mathfrak o}}
\nc\mfp{{\mathfrak p}}
\nc\mfq{{\mathfrak q}}
\nc\mfr{{\mathfrak r}}
\nc\mfs{{\mathfrak s}}
\nc\mft{{\mathfrak t}}
\nc\mfu{{\mathfrak u}}
\nc\mfv{{\mathfrak v}}
\nc\mfw{{\mathfrak w}}
\nc\mfx{{\mathfrak x}}
\nc\mfy{{\mathfrak y}}
\nc\mfz{{\mathfrak z}}

\nc\pt{\mathrm{pt}}
\nc{\one}{{\mathbf{1}}}

\nc\clambda{ {\check{\lambda} }}
\nc\cmu{ {\check{\mu} }}

\nc\loccit{\emph{loc.cit.}}

\nc{\ot}{ \mathop{\otimes}\displaylimits  }
\nc{\mt}{ \mathop{\times}\displaylimits  }
\nc{\colim}{ \mathop{\on{colim}\,}\displaylimits  }

\makeatletter
\newcommand{\laxto}{\dashedrightarrow}
\newcommand{\xrightleftarrows}[1]{\mathrel{\substack{\xrightarrow{#1} \\[-.9ex] \xleftarrow{#1}}}}
\newcommand{\adj}{\xrightleftarrows{\rule{0.5cm}{0cm}}}

\newcommand*{\da@rightarrow}{\mathchar"0\hexnumber@\symAMSa 4B }
\newcommand*{\da@leftarrow}{\mathchar"0\hexnumber@\symAMSa 4C }
\newcommand*{\xlaxto}[2][]{%
  \mathrel{%
    \mathpalette{\da@xarrow{#1}{#2}{}\da@rightarrow{\,}{}}{}%
  }%
}
\newcommand{\xlaxgets}[2][]{%
  \mathrel{%
    \mathpalette{\da@xarrow{#1}{#2}\da@leftarrow{}{}{\,}}{}%
  }%
}
\newcommand*{\da@xarrow}[7]{%
  \sbox0{$\ifx#7\scriptstyle\scriptscriptstyle\else\scriptstyle\fi#5#1#6\m@th$}%
  \sbox2{$\ifx#7\scriptstyle\scriptscriptstyle\else\scriptstyle\fi#5#2#6\m@th$}%
  \sbox4{$#7\dabar@\m@th$}%
  \dimen@=\wd0 %
  \ifdim\wd2 >\dimen@
    \dimen@=\wd2 %   
  \fi
  \count@=2 %
  \def\da@bars{\dabar@\dabar@}%
  \@whiledim\count@\wd4<\dimen@\do{%
    \advance\count@\@ne
    \expandafter\def\expandafter\da@bars\expandafter{%
      \da@bars
      \dabar@ 
    }%
  }%  
  \mathrel{#3}%
  \mathrel{%   
    \mathop{\da@bars}\limits
    \ifx\\#1\\%
    \else
      _{\copy0}%
    \fi
    \ifx\\#2\\%
    \else
      ^{\copy2}%
    \fi
  }%   
  \mathrel{#4}%
}
\makeatother

\nc{\Hom}{\on{Hom}}
\nc{\bHom}{\mathbf{Hom}}
\nc{\End}{\on{End}}
\nc{\Sym}{\on{Sym}}
\nc{\Tot}{\on{Tot}}
\nc{\DGCat}{\on{DGCat}}
\nc{\bDGCat}{\mathbf{DGCat}}
\renc{\Pr}{\on{Pr}}
\nc{\onpr}{\on{pr}}
\nc{\Spec}{\on{Spec}}
\nc{\Reg}{\on{Reg}}
\nc{\Ad}{\on{Ad}}
\nc{\Rep}{\on{Rep}}
\nc{\Specm}{\on{Specm}}
\renc{\mod}{\on{-mod}}
\nc{\comod}{\on{-comod}}
\nc{\bimod}{\on{BiMod}}
\renc{\bmod}{\on{-}\mathbf{mod}}
\nc{\alg}{\on{-alg}}
\nc{\id}{\mathrm{id}}
\nc{\Vect}{\on{Vect}}
\nc{\Res}{\on{Res}}
\nc{\Ind}{\on{Ind}}
\nc{\ind}{\mathbf{ind}}
\nc{\coind}{\mathbf{coind}}
\nc{\res}{\mathbf{res}}
\nc{\inv}{\mathbf{inv}}
\nc{\coinv}{\mathbf{coinv}}
\nc{\oninv}{{\on{inv}}}

\nc{\unit}{\mathbf{unit}}
\nc{\counit}{\mathbf{counit}}

\nc{\Sch}{\on{Sch}}
\nc{\IndSch}{\on{IndSch}}
\nc{\PreStk}{\on{PreStk}}
\nc{\AlgStk}{\on{AlgStk}}
\nc{\QCoh}{\on{QCoh}}
\nc{\Coh}{\on{Coh}}
\nc{\Shv}{\on{Shv}}
\nc{\Dmod}{\on{D}}
\nc{\DMod}{\on{Dmod}}
\nc{\bCorr}{\mathbf{Corr}}
\nc{\Corr}{\on{Corr}}

\nc{\Funct}{\on{Funct}}
\nc{\LFun}{\on{LFun}}
\nc{\bFun}{\mathbf{Fun}}
\nc{\affSch}{\on{Sch}^{\on{aff}}}
\nc{\oblv}{\mathbf{oblv}}
\nc{\Av}{\mathbf{Av}}
\nc{\pr}{\mathbf{pr}}
\renc{\ker}{\mathbf{ker}}

\nc{\triv}{\mathbf{triv}}
\nc{\mult}{\mathbf{mult}}
\nc{\comult}{\mathbf{comult}}
\nc{\Id}{ \mathbf{Id} }

\nc{\Cat}{\on{-Cat}}
\nc{\Map}{\on{Maps}}
\nc{\bMap}{\mathbf{Maps}}
\nc{\bCat}{\on{-}\mathbf{Cat}}
\nc{\oneCat}{1\Cat}
\nc{\inftyone}{(\infty,1)}
\nc{\inftytwo}{(\infty,2)}
\nc{\act}{\curvearrowright}
\nc{\ract}{\curvearrowleft}
\nc{\bact}{\mathbf{act}}
\nc{\bcoact}{\mathbf{coact}}

\nc{\Gr}{\on{Gr}}
\nc{\dualize}{\mathbf{dualize}}

\nc{\Pro}{\on{Pro}}
\nc{\inj}{\hookrightarrow}
\nc{\surj}{\twoheadrightarrow}
\nc{\Ran}{\on{Ran}}

\nc{\givesto}{\rightsquigarrow}
\nc{\toto}{\longrightarrow}
\nc{\xto}{\xrightarrow}
\nc{\xgets}{\xleftarrow}
\nc{\os}{\overset}
\nc{\us}{\underset}
\nc{\supp}{\on{supp}}
\nc{\reg}{{\on{reg}}}
\nc{\opreg}{{\on{opreg}}}

\nc{\Grp}{{\on{Grp}}}

\nc{\Zas}{  \on{Zas}_\rel }
\nc{\Schi}{ \on{Sch}_\rel}  

\nc{\Tw}{{\on{Tw}}}
\nc{\Del}{ \mathbf{\Delta}}
\nc\abso{{\on{abs}}}
\nc\rel{{\on{rel}}}
\nc\ol{\overline}
\nc\wt{\widetilde}
\nc{\ul}{\underline}
\nc{\wh}{\widehat}
\nc{\hs}{\heartsuit}
\nc{\lnilp}{{_{\on{loc.nilp.}}}}
\nc{\cohb}{{\le\infty,\ge-\infty}}
\nc{\st}{{\on{st}}}
\nc{\cocomplt}{{\on{cocomplt}}}
\nc{\cont}{{\on{cont}}}
\nc{\pres}{{\on{pres}}}
\nc{\lax}{{\on{lax}}}
\nc{\marked}{{\on{marked}}}
\nc{\enh}{{\on{enh}}}
\nc{\un}{{\on{un}}}
\nc{\ad}{{\on{ad}}}
\nc{\pos}{{\on{pos}}}
\nc{\disj}{{\on{disj}}}
\nc{\ora}{\overrightarrow}
\nc{\ola}{\overleftarrow}
\nc{\strict}{{\on{strict}}}
\nc{\red}{{\on{red}}}
\nc{\fact}{{\on{fact}}}
\nc{\co}{{\on{co}}}
\nc{\op}{{\on{op}}}
\nc{\conj}{ {\on{conj}} }
\nc{\ld}{{ {ld}}}
\nc{\rd}{{ {rd}}}
\nc{\rev}{{\on{rev}}}
\nc{\dR}{{\on{dR}}}
\nc{\fp}{{\on{fp}}}
\nc{\ft}{{\on{ft}}}
\nc{\ift}{{\on{ift}}}
\nc{\lft}{{\on{lft}}}
\nc{\qcqs}{{\on{qcqs}}}
\nc{\gen}{{\on{gen}}}
\nc{\hgen}{{\on{-gen}}}
\nc{\rotshriek}{{\rotatebox[origin=c]{180}{!}}}
\nc{\glob}{{\on{glob}}}
\nc{\att}{{\on{att}}}
\nc{\rep}{{\on{rep}}}
\nc{\fix}{{\on{fix}}}
\nc{\Bru}{{\on{Bruhat}}}
\nc{\open}{{\on{open}}}

\nc{\proper}{\on{proper}}
\nc{\all}{\on{all}}
\nc{\closed}{\on{closed}}
\nc{\oso}{\os{\circ}}
\nc{\twocat}{2\on{-Cat}}
\nc{\placid}{{\on{placid}}}
\nc{\hol}{{\on{indhol}}}
\nc{\univ}{{\on{univ}}}
\nc{\diag}{{\on{diag}}}

\nc{\level}{{\on{level}}}
\nc{\xyshort}{\xymatrixrowsep{0.5cm}}
\nc{\xysshort}{\xymatrixcolsep{0.5cm}}

\nc{\semiinf}{{\frac{\infty}{2}}}
\nc{\str}{{\on{str}}}
\nc{\pair}{{\on{pair}}}
\nc{\conv}{{\on{conv}}}
\nc{\bConv}{{\mathbf{Conv}}}
\nc{\good}{{\on{good}}}
\nc{\ext}{{\on{ext}}}
\nc{\bt}{{\blacktriangle}}
\nc{\dtriv}{{\on{triv}}}
\nc{\diff}{{\on{diff}}}

\nc{\etale}{$\acute{\on{e}}$tale }
\nc{\Kunneth}{K$\on{\ddot{u}}$nneth }
\nc{\cech}{$\on{\breve{C}}$ech }
\nc{\Plucker}{Pl$\on{\ddot{u}}$cker }

\nc{\Levi}{{\on{Levi}}}
\nc{\semi}{{\on{semi}}}
\nc{\QCAD}{{\on{QCAD}}}

\nc{\ltri}{{\vartriangleleft}}
\nc{\rtri}{{\vartriangleright}}

\nc{\aug}{{\on{aug}}}
\nc{\Vin}{\on{Vin}}
\nc{\VinH}{\on{VinH}}
\nc{\VinGr}{\on{VinGr}}
\nc{\VinG}{\on{Vin}_G}
\nc{\Bun}{\on{Bun}}
\nc{\BunG}{\Bun_G}
\nc{\BunGG}{\Bun_{G\times G}}
\nc{\BunM}{\Bun_M}
\nc{\BunMM}{\Bun_{M\times M}}
\nc{\BunP}{\Bun_P}
\nc{\BunPm}{\Bun_{P^-}}
\nc{\BunPPm}{\Bun_{P\times P^-}}
\nc{\BunPmP}{\Bun_{P^-\times P}}
\nc{\VinBun}{\on{VinBun}}
\nc{\Br}{\on{Br}}

\nc{\MGPos}{{M,G\on{-}\pos}}
\nc{\TGPos}{{T,G\on{-}\pos}}

\nc{\GrG}{\Gr_G}
\nc{\GrP}{\Gr_P}
\nc{\GrPm}{\Gr_{P^-}}
\nc{\GrM}{\Gr_M}
\nc{\GrGI}{\Gr_{G,I}}
\nc{\GrGGI}{\Gr_{G\times G,I}}
\nc{\GrGJ}{\Gr_{G,J}}
\nc{\GrPI}{\Gr_{P,I}}
\nc{\GrPmI}{\Gr_{P^-,I}}
\nc{\GrPPmI}{\Gr_{P\times P^-,I}}
\nc{\GrPmPI}{\Gr_{P^-\times P,I}}
\nc{\GrMI}{\Gr_{M,I}}
\nc{\GrMMI}{\Gr_{M\times M,I}}

\nc{\LUI}{{\mCL U_I}}
\nc{\LUmI}{{\mCL U^-_I}}
\nc{\LMI}{{\mCL M_I}}
\nc{\LGI}{{\mCL G_I}}
\nc{\LPI}{{\mCL P_I}}
\nc{\LPmI}{{\mCL P^-_I}}
\nc{\LpMI}{{\mCL^+ M_I}}
\nc{\UKMO}{{\mCL U\mCL^+M_I}}

\nc{\mon}{\on{-um}}
\nc{\sect}{{\on{sect}}}
\nc{\ontriv}{{\on{triv}}}
\nc{\df}{{\on{dfstr}}}
\nc{\defect}{\df }
\nc{\inftyx}{{\infty\cdot x}}
\renc{\setminus}{{-}}

\nc{\IGP}{\mbI(G,P)}
\nc{\IGQ}{\mbI(G,Q)}
\nc{\Eis}{ \on{Eis} }
\nc{\CT}{ \on{CT} }
\nc{\Par}{ \on{Par} }
\nc{\PsId}{ \on{Ps-Id} }
\nc{\Dri}{ \mathbf{Dri} }

\nc{\Exo}{ {\on{Exo}} }

\nc{\wtG}{ {\wt{G}} }

\nc{\Tadp}{ {T_\ad^+} }
\nc{\Act}{ { \on{ActAlgStk}_\lft} }
\nc{\ActPair}{  {\on{ActArrAlgStk}_\lft} }
\nc{\IsoGrp}{ {\on{IsoGrp} }}

\nc{\gCP}{ {\ge C_P} }
\nc{\Arr}{ {\on{Arr}} }

\nc{\onoblv}{\on{oblv}}

\begin{document}

\newpage
\title[Deligne--Lusztig duality on $\BunG$]{Deligne--Lusztig duality on the moduli stack of bundles}
\author{Lin Chen}
\address{Harvard Mathematics Department, 1 Oxford Street, Cambridge 02138, MA, USA}
\email{linchen@math.harvard.edu}
\begin{abstract}
Let $\BunG(X)$ be the moduli stack of $G$-torsors on a smooth projective curve $X$ for a reductive group $G$. We prove a conjecture made by Drinfeld-Wang and Gaitsgory on the Deligne--Lusztig duality for D-modules on $\BunG(X)$. This conjecture relates the pseudo-identity functors in \cite{gaitsgory2017strange}, \cite{drinfeld2015compact} to the enhanced Eisenstein series and geometric constant term functors in \cite{gaitsgory2015outline}. We also prove a ``second adjointness'' result for these enhanced functors.
\end{abstract}

\maketitle
\tableofcontents

\setcounter{section}{-1}
\section{Introduction}

\ssec{Motivation: Deligne--Lusztig duality}

The following pattern, which are called Deligne--Lusztig duality, has been observed in several representation-theoretic contexts: The composition of two different duality functors on the category $\mCC_G$ attached to a reductive group $G$ is isomorphic to a Deligne--Lusztig functor, given by a complex indexed by standard parabolic subgroups\footnote{We fix a Borel subgroup $B$ of $G$. A parabolic subgroup is \emph{standard} if it contains $B$.} $P$ of $G$, whose terms are compositions
$$ \mCC_G \xto{\CT_P}  \mCC_M \xto{\Eis_P} \mCC_G,$$
where

\begin{itemize}
	\item $M$ is the Levi quotient group of $P$;
	
	\item $\mCC_M$ is the category attached to $M$;
	
	\item $\CT_P$ and $\Eis_P$ are adjoint functors connecting $\mCC_G$ and $\mCC_M$.
\end{itemize}
The study of such pattern can be dated back to the Alvis--Curtis duality operation on characters of finite Chevalley group (\!\!\cite{alvis1979duality}, \cite{curtis1980truncation}). Here are some examples of this patter:
\begin{itemize}
  \item The work of Bernstein-Bezrukavnikov-Kazhdan (\!\!\cite{bernstein2018deligne}), where $\mCC_G$ is the category of representations of the group $G(K)$, where $K$ is a non-Archimedian local field.

  \item The work of Yom Din (\!\!\cite{din2019deligne}), where $\mCC_G$ is the category of character D-modules on $G$, i.e., D-modules on the quotient stack $G/\Ad(G)$ of $G$ by its adjoint action.

  \item The work of Drinfeld-Wang (\!\!\cite{drinfeld2016strange} and \cite{wang2018invariant}), where $\mCC_G$ is the space of automorphic functions for the group $G$. Note that this example is actually one categorical level down from the above partern (i.e., one needs to replace ``categories'' by ``spaces'' and ``functors'' by ``operators'').
\end{itemize}

In the present paper we establish yet another incarnation of this pattern. Namely, we take $\mCC_G$ to be the category of automorphic D-modules\footnote{Our method can also be applied to the category of automorphic sheaves with suitable modifications.}, i.e., $\mCC_G = \DMod(\Bun_G(X))$, where $\Bun_G(X)$ is the moduli stack of $G$-torsors on a smooth complete curve $X$.

Our context can be viewed as directly categorifying that of Drinfeld-Wang. It is also closely connected to that of Yom Din because $G/\Ad(G)$ is an open substack of $\Bun_G(Y)$ where $Y$ is the \emph{nodal curve} obtained by gluing two points of $\mBP^1$.

Below we will review the contexts mentioned above.

\sssec{Work of \cite{bernstein2018deligne}}

Let $G$ be defined over a number field and $K$ be a non-archimedian local field. In \cite{bernstein2018deligne}, the authors proved the following result about the derived category $G(K)\mod$ of admissible representations of the $p$-adic group $G(K)$.

For any object $\mCM\in G(K)\mod$, consider the corresponding \emph{Deligne--Lusztig complex}\footnote{Analogous complexes for finite Chevalley group were firstly studied by Deligne and Lusztig in \cite{deligne1982duality}.}
$$  \on{DL}(\mCM):=  \; [  \mCM \to \bigoplus_P i_P^G\circ r_P^G(\mCM) \to \cdots \to i_B^G\circ r_B^G(\mCM)  ]    $$
where $(r_P^G,i_P^G)$ is the adjoint pair for the parabolic induction and Jacquet functors, and where the direct sum in the $k$-th term of the complex is taken over standard parabolic subgroups of corank $k$. The main theorem of \cite{bernstein2018deligne} says that
$$  \on{DL} \simeq \mBD^{\on{coh}} \circ  \mBD^{\on{contr}}[\on{rank}(G)],$$
where $\mBD^{\on{contr}}$ and $\mBD^{\on{coh}}$ are the contravariant endofunctors on $G(K)\mod$ for the \emph{contragredient} and \emph{cohomological} dualities. In other words,
$$\mBD^{\on{contr}}(\mCM):=\mCM^\vee,\;  \mBD^{\on{coh}}(\mCM) := \on{RHom}_{G(K)}(\mCM,\mCH),$$
where $\mCM^\vee$ is the admissible dual, and where $\mCH=C_c^\infty(G(K))$ is the regular bimodule for the Hecke algebra. The proof in \loccit\,used an explicit resolution for $\mCH$ coming from the geometry of the wonderful compactification of $G$.

\sssec{Work of \cite{din2019deligne}}
\label{sssec-work-yomdim}

Let $G$ be defined over an algebraically closed field $k$ of characteristic $0$. In \cite{din2019deligne}, the author proved the following result about the DG-category\footnote{See Notation \ref{notn-DG-categories} for our conventions for DG-categories.} $\Dmod(G/\!\Ad(G))$ of character D-modules on $G$.

Let $P$ be a standard parabolic subgroup and $M$ be its Levi quotient group. Consider the diagram
$$ G/\!\Ad(G) \xgets{p} P/\!\Ad(P) \xto{q} M/\!\Ad(M). $$
The map $p$ is projective and $q$ is smooth. Hence we have the \emph{parabolic restriction} functor
$$ \on{res}_P= q_*\circ p^! : \Dmod(G/\!\Ad(G)) \to \Dmod(M/\!\Ad(M)) $$
and its left adjoint, a.k.a. the \emph{parabolic induction} functor
$$ \on{ind}_P= p_*\circ q^! : \Dmod(M/\!\Ad(M)) \to \Dmod(G/\!\Ad(G)). $$
These functors are t-exact by \cite{bezrukavnikov2018parabolic}. Let $P^-$ be an opposite parabolic subgroup\footnote{In order to define $P^-$, we fix a Carton subgroup of $G$.}. It is known (see \cite[$\S$ 0.2.1]{drinfeld2014theorem}) that $ \on{res}_P$ is left adjoint to $\on{ind}_{P^-}$. This is analogous to the well-known Bernstein's second adjointness.

Consider the diagonal map $\Delta:G/\!\Ad(G)\to G/\!\Ad(G)\mt G/\!\Ad(G)$ and the endofunctor on $\Dmod(G/\!\Ad(G))$ given by the kernel $\Delta_!(k_{G/\!\Ad(G)})$, where $k_{G/\!\Ad(G)}$ is the constant D-module. Explicitly, this endofunctor is
$$ \onpr_{1,\bt}(  \Delta_!(k_{G/\!\Ad(G)})\ot^! \onpr_2^! (-) ) ,  $$
where $\onpr_{1,\bt}$ is the \emph{renormalized pushforward functor} in \cite{drinfeld2013some}. This endofunctor is the so-called \emph{Drinfeld-Gaitsgory functor} for $\Dmod( G/\!\Ad(G))$ in \cite{din2019deligne}.

One of the main results of \cite{din2019deligne} says that the above Drinfeld-Gaitsgory functor can be ``glued''\footnote{Roughly speaking, this means that up to a cohomological shift, the Drinfeld-Gaitsgory functor sends an object $\mCF\in  \Dmod(G/\!\Ad(G))^{\heartsuit}$ in the heart of the t-structure to a \emph{certain} complex
$$ \mCF \to \bigoplus_{P} \on{ind}_{P^-}\circ \on{res}_{P}(\mCF) \to \cdots \to \on{ind}_{B^-}\circ \on{res}_{B}(\mCF).$$
However, \cite{din2019deligne} did \emph{not} describe the connecting morphisms in the above complex. Nevertheless, we have confidence that one can use the method in the current paper to show that these connecting maps are given by the adjunction natural transformations of the second adjointness.} from the functors
$$ \bigoplus_{\on{rank}(P)=l} \on{ind}_{P^-}\circ \on{res}_{P}[l-\dim(T)]. $$
The proof in \loccit\,used a filtration of $ \Delta_!(k_{G/\!\Ad(G)})$ coming from the geometry of the wonderful compactification of $G$.

\sssec{The pseudo-identity functor(s)}

As explained in \cite{din2019deligne}, the above Drinfeld-Gaitsgory functor can be identified with
$$ \PsId_{ G/\!\Ad(G),! } \circ  (\PsId_{ G/\!\Ad(G),\on{naive} })^{-1} ,$$
where $\PsId$ are the pseudo-identity functors constructed in \cite{drinfeld2013some}. Here are more details.

Let $Y$ be a QCA algebraic stack\footnote{This means $Y$ is a quasi-compact algebraic stack whose automorphism groups of geometric points are affine.} in the sense of \loccit. Consider the cocomplete DG-category $\Dmod(Y)$ and its full subcategory $\Dmod(Y)^c$ of compact objects. Verdier duality provides an equivalence
$$ \mBD^{\on{Ver}}: \Dmod(Y)^{c} \to \Dmod(Y)^{c,\op}.$$
Using ind-completion, we obtain an equivalence
$$\PsId_{Y,\on{naive}}:  \Dmod(Y)^\vee \simeq \Dmod(Y), $$
where $ \Dmod(Y)^\vee$ is the \emph{Lurie dual}\footnote{See Notation \ref{notn-DG-categories} for what this means.} of $\Dmod(Y)$.

On the other hand, we have \emph{the product formula}:
$$ \Dmod(Y\mt Y) \simeq \Dmod(Y)\ot_k \Dmod(Y),$$
where $\ot_k$ is the \emph{Lurie tensor product} of cocomplete DG-categories. The RHS can be identified with $\LFun_k(\Dmod(Y)^\vee,\Dmod(Y))$, i.e., the category of $k$-linear colimit-preserving functors $\Dmod(Y)^\vee\to \Dmod(Y)$. Hence the object $\Delta_!(k_{Y})\in \Dmod(Y\mt Y)$ provides a functor 
$$ \PsId_{Y,!}:\Dmod(Y)^\vee \to \Dmod(Y).$$
In general, the functor $\PsId_{Y,!}$ is not an equivalence. We say $Y$ is \emph{miraculous} if it is an equivalence.

It is known that the equivalence $\PsId_{Y,\on{naive}}$ can be obtained in the same way by replacing $\Delta_!(k_{Y})$ by $\Delta_*(\omega_{Y})$, where $\omega_Y$ is the dualizing D-module on $Y$. For this reason, the functor $\PsId_{Y,\on{naive}}$ is also denoted by $\PsId_{Y,*}$ in the literature.

It follows from definitions that the composition $\PsId_{Y,!}\circ (\PsId_{Y,\on{naive}})^{-1}$ is the functor given by the kernel $\Delta_!(k_{Y})$. In other words, it is the Drinfeld-Gaitsgory functor for $\Dmod(Y)$.

\begin{rem} By \cite[Proposition 5.5]{din2019deligne}, the Drinfeld-Gaitsgory functor for $\Dmod(G/\!\Ad(G))$ is invertible, hence so is $\PsId_{!,G/\!\Ad(G)}$. In other words, $G/\!\Ad(G)$ is miraculous.
\end{rem}

\sssec{The current work}
Let $G$ be as in $\S$ \ref{sssec-work-yomdim} and $X$ be a connected smooth projective curve over $k$. Let $\BunG$ be the moduli stack of $G$-torsors on $X$. The purpose of this paper is to describe the Deligne--Lusztig duality on the DG-category $\Dmod(\BunG)$ of D-modules on $\BunG$.

Unlike $G/\!\Ad(G)$, the stack $\BunG$ is \emph{not} quasi-compact. Nevertheless, the main theorem of \cite{drinfeld2015compact} says that $\Dmod(\BunG)$ is compactly generated and hence dualizable. Also, the product formula
$$ \Dmod(\BunG\mt \BunG) \simeq \Dmod(\BunG)\ot_k \Dmod(\BunG) $$
still holds (see \cite[Remark 2.2.9]{drinfeld2015compact}). Hence as before, we have 
equivalences
\begin{equation} \label{eqn-duality-product}
\begin{aligned}
\LFun_k( \Dmod(\BunG)^\vee, \Dmod(\BunG) ) \simeq  \Dmod(\BunG) \ot_k  (\Dmod(\BunG)^\vee)^\vee  \simeq \\
\simeq \Dmod(\BunG)\ot_k \Dmod(\BunG) \simeq \Dmod(\BunG\mt \BunG),
\end{aligned}
\end{equation}
and we use the objects $\Delta_*(\omega_{\BunG}),\, \Delta_!(k_{\BunG})$ in the RHS to define functors
$$ \PsId_{\BunG,\on{naive}},\, \PsId_{\BunG,!}: \Dmod(\BunG)^\vee \to \Dmod(\BunG). $$
From now on, we write them just as $\PsId_{\on{naive}},\, \PsId_{!}$.

Unlike the quasi-compact case, the functor $\PsId_{\on{naive}}$ is \emph{not} an equivalence. On the contrary, the main theorem of \cite{gaitsgory2017strange} says:

\begin{itemize} 
	\item  The functor $\PsId_!: \Dmod(\BunG)^\vee \to \Dmod(\BunG) $ is an equivalence, i.e., $\BunG$ is miraculous.
\end{itemize}

Accordingly, the \emph{Deligne--Lusztig duality for $\BunG$} in this paper is actually analogous to the ``left adjoint version'' of \cite{bernstein2018deligne} and \cite{din2019deligne} because our $\Eis$-functor is \emph{left} adjoint to the $\CT$-functor but the contrary is true in their work. Namely, we will show

\begin{itemize}
	\item Up to cohomological shifts, the endofunctor $\PsId_{\on{naive}} \circ \PsId_{!}^{-1}$ on $\Dmod(\BunG)$ can be ``glued'' from the functors
	$$ \bigoplus_{\on{corank}(P)=l} \Eis^\enh_{P\to G}\circ \CT^\enh_{G\gets P},$$
	where $\Eis^\enh_{P\to G}$ (resp. $\CT^\enh_{G\gets P}$) is the enhanced Eisenstein series (resp. enhanced constant term) functor\footnote{We review the definitions of them in $\S$ \ref{intro-IGP}.} on $\Dmod(\BunG)$.
\end{itemize}
The precise statement of our main theorem will be given in $\S$ \ref{ssec-main-theorem}. See Theorem \ref{thm-main}.

\begin{rem} For $G=\on{SL}_2$, our main theorem was conjectured by V. Drinfeld and J. Wang in \cite[Conjecture C.2.1]{drinfeld2016strange}. For general $G$, Wang made the following remark in \cite[Remark 6.6.5]{wang2018invariant}: 
\begin{itemize}
	\item[] \emph{...describing the functor inverse to $\PsId_{\BunG,!}$ (we expect that one can mimic the construction of the Deligne--Lusztig complex using the functors $\Eis_P^\enh,\CT_P^\enh$)}. 
\end{itemize}
However, as far as we know, the precise formulation\footnote{Since the functors $\Eis_P^\enh,\CT_P^\enh$ are not t-exact, the naïve formation of the Deligne--Lusztig complex is not a well-defined object in the DG-category $\Dmod(\BunG)$.} of the conjecture for general $G$ was first made by D. Gaitsgory and recorded by D. Beraldo in \cite[$\S$ 1.5.5]{beraldo2019deligne}.
\end{rem}

\sssec{Relation with the Drinfeld-Wang operator}
\label{sssec-Drinfeld-Wang}
Drinfeld and Wang made their conjecture according to an analogous result on the space of automorphic forms proved by them in \cite{drinfeld2016strange} and \cite{wang2018invariant}. Let us briefly explain their work.

Let $F$ be a global function field over $\mBF_q$ and $\mBA$ be the adele ring of $F$. Let $G$ be a split reductive group over $\mBF_q$ and $G(\mBO)$ be the maximal compact subgroup of $G(\mBA)$. Let $\mCA_c$ be the space of compactly supported smooth $G(\mBO)$-finite functions on $G(\mBA)/G(F)$. As explained in \cite[Appendix A]{drinfeld2016strange}, the DG-category $\Dmod(\BunG)^\vee$ (when $G$ is in characteristic $0$) can be viewed as a geometric analogue of $\mCA_c^{G(\mBO)}$ (the subspace of $\mCA_c$ fixed by $G(\mBO)$).

Drinfeld and Wang also defined a subspace $\mCA_{pc\on{-}c}\subset \mCA_c$ such that $\mCA_{pc\on{-}c}^{G(\mBO)}$ is analogous to $\Dmod(\BunG)$. They also constructed a $G(\mBA)$-linear operator 
$$L:\mCA_c\to \mCA_{pc\on{-}c}$$
such that $L^{G(\mBO)}: \mCA_c^{G(\mBO)}\to \mCA_{pc\on{-}c}^{G(\mBO)}$ is analogous to the functor $\PsId_!$. Moreover, they proved $L$ is invertible and gave the following explicit formula for its inverse:
$$ L^{-1}f = \sum_P (-1)^{\dim Z_M} \Eis_P \circ \CT_P(f),$$
where $\Eis_P,\CT_P$ are the Eisenstein and constant term operators, and where $Z_M$ is the center of $M$.

Our main theorem can be viewed as a categorification of the above formula (when restricted to $G(\mBO)$-invariant subspaces). We refer the reader to \cite[Appendix C]{drinfeld2016strange} for more details on this analogy.

\ssec{Recollections: The parabolic category \texorpdfstring{$\IGP$}{I(G,P)}}
\label{intro-IGP}
From now on, we fix a connected reductive group\footnote{See Notation \ref{notn-convension-reductive-group} for our notations for concepts related to $G$.} $G$ defined over an algebraically closed field $k$ of characteristic $0$. For simplicity, we assume $[G,G]$ to be simply connected. We also fix a connected smooth projective curve $X$ over $k$.

Let $P$ be a parabolic subgroup of $G$ and $M$ be its Levi quotient group. Consider the diagram
$$  \BunG \xgets{\mfp_P} \BunP \xto{\mfq_P} \BunM.$$
In \cite{braverman2002geometric} and \cite{drinfeld2016geometric}, the authors constructed the \emph{geometric Eisenstein series functor} and the \emph{geometric constant term functor}
\begin{eqnarray*}
 \Eis_{P,!}= \mfp_{P,!}\circ \mfq_P^* :& \Dmod(\BunM) \to \Dmod(\BunG)  \\
 \CT_{P,*}= \mfq_{P,*}\circ \mfp_P^! :& \Dmod(\BunG) \to \Dmod(\BunM).
\end{eqnarray*}
However, they are not the functors that we will use (otherwise the main theorem would be false). Instead, we need to replace $\Dmod(\BunM)$ by the category $\IGP$ defined in \cite[$\S$ 6]{gaitsgory2015outline}, and accordingly use the ``enhanced'' functors $(\Eis_{P}^\enh,\CT^\enh_{P})$ defined there. We review these functors in this subsection.

\begin{rem}
As explained in \loccit, one can think of $\IGP$ as the DG-category of D-modules on a \emph{non-existent} stack obtained by gluing all the connected components of $\BunM$ together. Since this imaginary stack has the same field-valued points as $\BunM$, the difference between $\IGP$ and $\Dmod(\BunM)$ can not be seen in their ``de-categorifications''. In other words, both $\Eis_{P,!}$ and $\Eis_{P}^\enh$ are analogous to the Eisenstein series operator for automorphic functions.
\end{rem}

\sssec{Prestack of generic reductions}
Let $K$ be any algebraic group and $H$ be a subgroup of $K$. In \cite[Example 2.2.5]{barlev2014d}, J. Barlev constructed a prestack\footnote{See Definition \ref{defn-prestack} for what this means.} $\Bun_K^{H\on{-gen}}$ classifying $K$-torsors on $X$ equipped with a generic reduction to $H$. In the notation of \loccit, it is defined by
$$ \Bun_K^{H\on{-gen}} := \bMap(X, \mBB K ) \mt_{\mathbf{GMap}(X,\mBB K)} \mathbf{GMap}(X,\mBB H),$$
where $\mathbf{GMap}(X,Y)$ is the prestack classifying generic maps from $X$ to $Y$. To simplify the notation, we write the RHS as $\bMap_\gen(X,\mBB K\gets \mBB H)$. More generally, for any map between lft prestacks $Y_1\to Y_2$, we define
$$\bMap_\gen(X,Y_2\gets Y_1):=  \bMap(X, Y_2 ) \mt_{\mathbf{GMap}(X,Y_2)} \mathbf{GMap}(X,Y_1)$$
For future use, let us mention that the functor $\bMap_\gen(X,-)$ commutes with finite limits (see \cite[Remark 2.2.6]{barlev2014d}).

Applying the above construction to the diagram
\[ (G,P) \gets (P,P) \to (M,M), \]
we obtain a diagram
\[ \Bun_G^{P\on{-gen}} \xgets{\iota_P} \BunP \xto{\mfq_P}\ \BunM.\]

\begin{rem} The prestack $\Bun_G^{P\on{-gen}}$ has the same field-valued points as $\BunP$.
\end{rem}

\begin{defn} The DG-category $\IGP$ is defined as the fiber product of the following diagram:
\[
\xyshort
\xymatrix{
	\IGP \ar@{.>}[r] \ar@{.>}[d] &
	\Dmod(\BunM) \ar[d]^-{\mfq_P^*} \\
	\Dmod(\Bun_G^{P\on{-gen}} ) 
	\ar[r]^-{\iota_P^!}
	&\Dmod(\BunP).
}
\]
\end{defn}

\begin{rem} The above definition is equivalent to that in \cite[$\S$ 6.1]{gaitsgory2015outline} by \cite[Lemma 6.3.3]{gaitsgory2015outline}.
\end{rem}

\begin{rem} \label{rem-conservative-iota^!}
By \cite[Lemma 6.1.2]{gaitsgory2015outline}, the functor $\iota_P^!$ is conservative. By \cite[$\S$ 6.2.1]{gaitsgory2015outline}, the functor $\mfq_P^*$ is fully faithful. Therefore the functor
$$ \IGP \to \Dmod(\Bun_G^{P\on{-gen}} ) $$
is fully faithful and the functor $\IGP  \to \Dmod(\BunM)$ is conservative. Following \loccit, we denote the last functor by 
$$\iota_M^!: \IGP  \to \Dmod(\BunM) . $$
\end{rem}

The following result was claimed in \cite[$\S$ 6.2.5]{gaitsgory2015outline}. We provide a proof in Appendix \ref{appendix-proof-prop-well-defined-iota_M_!}.

\begin{prop} \label{prop-well-defined-iota_M_!}
(Gaitsgory)

 The partially defined left adjoint $\iota_{P,!}$ to $\iota_P^!$ is well-defined on the essential image of $\mfq_P^*$, and its image is contained in $\IGP$.
\end{prop}

\begin{cor} \label{cor-well-defined-iota_M_!}
The functor $\iota_M^!: \IGP  \to \Dmod(\BunM)$ has a left adjoint
$$\iota_{M,!}: \Dmod(\BunM)\to  \IGP.$$
\end{cor}

\proof By Proposition \ref{prop-well-defined-iota_M_!}, the functor $\iota_{P,!} \circ \mfq_P^*$ uniquely factors through a functor $\Dmod(\BunM)\to \IGP$, which is the desired left adjoint.

\qed[Corollary \ref{cor-well-defined-iota_M_!}]

\sssec{Enhanced Eisenstein series functor} Let $Q$ be another parabolic subgroup of $G$ that contains $P$. Consider the map
\[ \mfp_{P\to Q}^\enh:  \Bun_G^{P\on{-gen}} \to \Bun_G^{Q\on{-gen}}\]
and the functor
\[ \mfp_{P\to Q}^{\enh,!}: \Dmod( \Bun_G^{Q\on{-gen}} ) \to \Dmod (\Bun_G^{P\on{-gen}})  \]
The special case (when $Q=G$) of the following result was claimed in \cite[Lemma 6.3.3]{gaitsgory2015outline}. We provide a proof in Appendix \ref{appendix-proof-prop-well-defined-Eis_enh}.

\begin{prop} \label{prop-well-defined-Eis_enh}
(Gaitsgory) 

\begin{itemize}
	\item[(1)] The partially defined left adjoint $\mfp_{P\to Q,!}^\enh$ to $\mfp_{P\to Q}^{\enh,!}$ is well-defined on $\IGP\subset \Dmod(  \Bun_G^{P\on{-gen}} )$, and sends $\IGP$ into $\IGQ$.

	\item[(2)] Let
	$$\Eis^\enh_{P\to Q}: \IGP \to \IGQ $$
	be the functor obtained from $\mfp_{P\to Q,!}^\enh$. Then $\Eis^\enh_{P\to Q}$ has a \emph{continuous} right adjoint 
	$$ \CT^\enh_{Q\gets P}: \IGQ \to \IGP.$$
\end{itemize}
\end{prop}

\begin{rem} When $Q=G$, we also denote the adjoint pair $(\Eis^\enh_{P\to G},\CT^\enh_{G\gets P} )$ by $(\Eis^\enh_P, \CT^\enh_P)$.
\end{rem}

\begin{warn} The functor $\mfp_{P\to Q}^{\enh,!}$ does not send $\IGQ$ into $\IGP$. Hence the functor $\CT^\enh_{Q\gets P}$ is not the restriction of $\mfp_{P\to Q}^{\enh,!}$. Instead, it is given by $\on{Av}^{\mbU_P}\circ \mfp_{P\to Q}^{\enh,!}$, where $\on{Av}^{\mbU_P}$ is the right adjoint to the inclusion $\IGP\subset \Dmod(\Bun_G^{P\on{-gen}} )$. We refer the reader to \cite[$\S$ 6.1.3]{gaitsgory2015outline} for the meaning of the symbol $\on{Av}^{\mbU_P}$.
\end{warn}

\ssec{The main theorem}
\label{ssec-main-theorem}

We fix a Borel subgroup $B$ of $G$. Let $\Par$ be the poset of standard parabolic subgroups (i.e., parabolic subgroups containing $B$) of $G$. Let $\Par'$ be $\Par\setminus\{G\}$. We view posets as categories in the standard way. It follows formally from Proposition \ref{prop-well-defined-Eis_enh} that we have a functor
$$ \mbD\mbL: \Par \to \LFun_k(\Dmod(\BunG), \Dmod(\BunG) ),\; P\mapsto \Eis^\enh_{P\to G} \circ \CT^\enh_{G\gets P}$$
such that a morphism $P\to Q$ in $\Par$ is sent to the composition 
$$ \Eis^\enh_{P\to G} \circ \CT^\enh_{G\gets P} \simeq \Eis^\enh_{Q\to G} \circ \Eis^\enh_{P\to Q} \circ \CT^\enh_{Q\gets P} \circ  \CT^\enh_{G\gets Q} \to  \Eis^\enh_{Q\to G} \circ  \CT^\enh_{G\gets Q}.$$
Note that $\mbD\mbL(G)$ is the identity functor $\Id$ on $\Dmod(\BunG)$.

Since $G$ is the final object in $\Par$. The above functor $\mbD\mbL$ induces a natural transformation $\colim_{P\in \Par'} \mbD\mbL(P)  \to \mbD\mbL(G)$. Our main theorem is
\begin{thm} \label{thm-main}
There is a canonical equivalence
\begin{equation} \label{eqn-thm-main}
 \on{coFib} (  \colim_{P\in \Par'} \mbD\mbL(P)  \to \mbD\mbL(G) ) \simeq \PsId_{\on{naive}} \circ \PsId_!^{-1}[-2\dim(\BunG)-\dim(Z_G)]
 \end{equation}
in $\LFun_k(\Dmod(\BunG), \Dmod(\BunG) )$.
\end{thm}

\begin{rem} \label{rem-thm-main-complex}
Let $\mCF\in \Dmod(\BunG)^\heartsuit$ be an object in the heart of the t-structure. If the functors $\mbD\mbL(P)$ were t-exact (which is not true), then the value of the LHS of (\ref{eqn-thm-main}) on $\mCF$ could be calculated by a complex
$$ \mbD\mbL(B)(\mCF) \to  \cdots \to \bigoplus_{\on{corank}(P)=1} \mbD\mbL(P)(\mCF)  \to  \mbD\mbL(G)(\mCF).$$
Hence the LHS of (\ref{eqn-thm-main}) can be viewed as an $\infty$-categorical analogue for the Deligne--Lusztig complex.
\end{rem}

\sssec{A stronger result}
As mentioned in \cite[Remark 1.5.6]{beraldo2019deligne}, Gaitsgory's strategy for the proof of the above theorem is to express both sides via the Drinfeld's compactification $\ol{\Bun}_G:=  \VinBun_G/T$, where
$$ \VinBun_G:= \bMap_\gen(X, G\backslash \Vin_G/G \supset G\backslash\, _0\!\Vin_G/G ).$$
We refer the reader to \cite[$\S$ 2.2.4]{schieder2016geometric}, \cite[$\S$ 2.3.3]{finkelberg2020drinfeld} and \cite[$\S$ 1]{chen2021thesis} for a detailed discussion about $\Vin_G$ and $\VinBun_G$. For now, it is enough to know
\footnote{\label{Footnote:SL2-olBunG}
It is helpful to have the case $G=\on{SL}_2$ in mind. In this case, $\VinBun_{G}(S)$ classifies $\phi:E_1\to E_2$, where $E_1$ and $E_2$ are rank $2$ vector bundles on $X\mt S$ with trivialized determinant line bundles, and $\phi$ is a map between \emph{coherent} $\mCO_{X\mt S}$-modules such that for any geometric point $s$ of $S$, the map $\phi|_{X\mt s}$ is nonzero. The Cartan subgroup $\mBG_m$ acts on $\VinBun_{\on{SL}_2}$ by scalar multiplication on $\phi$.

We warn the reader that the projection $\VinBun_G\to \Bun_{G} \mt \Bun_{G}$ below sends the above data to $(E_2,E_1)$ rather than $(E_1,E_2)$. See \cite[Warning 1.1.3]{chen2021thesis} for the reason of this convention.}:

\begin{itemize}
	\item $\ol{\Bun}_G$ is an algebraic stack;
	
	\item There is a canonical map $\ol{\Bun}_G \to T_\ad^+ /T$, where $T_\ad^+:=\mBA^\mCI \supset \mBG_m^\mCI \simeq T_\ad$ is the semi-group completion\footnote{See Notation \ref{notn-convension-semi-group-completion} for more information about $T_\ad^+$.} of the adjoint torus $T_\ad$;
	
	\item The diagonal map $\Delta:\BunG\to \BunG\mt \BunG$ canonically factors as
	$$ \BunG \xto{b}  \ol{\Bun}_G  \xto{\ol\Delta} \BunG\mt \BunG $$
	and the map $\ol\Delta$ is schematic and proper.
	
	\item Let $Z_G$ be the center of $G$. The map $b:\BunG\to \ol{\Bun}_G$ canonically factors as
	$$ \BunG \xto{r} \BunG \mt \mBB Z_G \xto{j_G} \ol{\Bun}_G $$
	and the map $j_G$ is a schematic open embedding onto
	$$ \ol{\Bun}_{G,G}:= \ol{\Bun}_G \mt_{ T_\ad^+ /T } T_\ad /T.$$
\end{itemize}

The coordinate stratification on $T_\ad^+:=\mBA^\mCI$ can be labelled by $\Par$ such that the $G$-stratum is the open stratum $\mBG_m^\mCI\subset \mBA^\mCI$ while the $B$-stratum is the closed stratum $0\inj \mBA^\mCI$ (see Notation \ref{notn-convension-semi-group-completion}). This stratification induces a stratification on $\ol{\Bun}_G$ labelled by $\Par$, known as \emph{the parabolic stratification}.

For a standard parabolic $P\in \Par$, the $P$-stratum is defined as
$$ \ol{\Bun}_{G,P}:= \ol{\Bun}_G \mt_{ T_\ad^+/T} T_{\ad,P}^+/T. $$
We deonte the corresponding locally closed embedding by
$$ i_P: \ol{\Bun}_{G,P} \to \ol{\Bun}_G. $$

Note that $\Par$ is isomorphic to the power poset $P(\mCI)$ of $\mCI$. By a general construction (see Corollary \ref{Cor-gluing-functors}) for stacks stratified by power posets, we have a canonically defined functor\footnote{The functor $\mbK$ is given by $\mbG^*_{ r_!(k_{\Bun_G}), \ol{\Bun}_G }$, which is defined in Corollary \ref{Cor-gluing-functors}. Here we use the full category $\Dmod_\hol(\ol{\Bun}_G)\subset \Dmod(\ol{\Bun}_G)$ of ind-holonomic D-modules to make sure the $*$-pullback and $!$-pushforward are well-defined.
}
\begin{eqnarray} \label{eqn-mbK}
\mbK&:& \Par \to \Dmod_\hol(\ol{\Bun}_G),
\\\nonumber
 & & P \mapsto i_{P,!}\circ i_{P}^* \circ j_{G,*}\circ r_!(k_{\Bun_G}) [\on{rank}(M)-\on{rank}(G)],
\end{eqnarray}
and a canonical isomorphism (see Lemma \ref{lem-gluing-functors-colimite})
\begin{equation} \label{eqn-resolution-*-extension-on-olBunG}
\on{coFib} ( \colim_{ P\in \Par' }  \mbK(P) \to \mbK(G) ) \simeq j_{G,*}\circ r_!(k_{\Bun_G}).
\end{equation}
Consider the composition
\begin{equation} \label{eqn-functor-mbE}
 \mbE: \LFun_k( \Dmod(\BunG),\Dmod(\BunG) ) \to  \LFun_k( \Dmod(\BunG)^\vee,\Dmod(\BunG) ) \simeq  \Dmod(\BunG\mt \BunG), 
\end{equation}
where the first functor is given by precomposition with $\PsId_!$ and the last equivalence is (\ref{eqn-duality-product}). Equivalently, $\mbE$ sends an endo-functor $F$ to 
$$\mbE(F)\simeq (F\ot\Id)\circ \Delta_!(k_{\Bun_G}),$$
where we view $F\ot\Id$ as an endo-functor of $\Dmod(\BunG\mt \BunG)\simeq \Dmod(\BunG) \ot_k \Dmod(\BunG)$ .

Let us first deduce Theorem \ref{thm-main} from the following stronger result, which will be proved in $\S$ \ref{ssec-proof-main-theorem} (but with some details postponed to the later sections).

\begin{thm} \label{thm-main-variant}
There is a canonical commutative diagram
$$
\xyshort
\xymatrix{
	\Par \ar[r]^-{\mbD\mbL} \ar[d]^-{\mbK} &
	\LFun_k( \Dmod(\BunG),\Dmod(\BunG) ) \ar[d]^-{\mbE} \\
	\Dmod_\hol(\ol{\Bun}_G) \ar[r]^-{\ol{\Delta}_!}
	& \Dmod(\BunG\mt \BunG).
}
$$
\end{thm}

\sssec{Deduction of Theorem \ref{thm-main} from Theorem \ref{thm-main-variant}} We first deduce our main theorem from Theorem \ref{thm-main-variant}. This step is due to Gaitsgory.

First note that $\mbE$ is an equivalence because $\PsId_!$ is. By definition,
$$\mbE^{-1}( \Delta_*(\omega_{\BunG}) ) \simeq \PsId_{\on{naive}}\circ \PsId_!^{-1}.$$

On the other hand, as in \cite[$\S$ 3.2.3]{gaitsgory2017strange}, we have an isomorphism
$$ \Delta_*(\omega_{\BunG})[ -2\dim(\BunG)-\dim(Z_G) ] \simeq \ol{\Delta}_! \circ j_{G,*}\circ  r_!( k_{\BunG} ),$$
where the cohomological shift by $-2\dim(\BunG)$ is due to the difference between $\omega_{\BunG}$ and $k_{\BunG}$, while that by $-\dim(Z_G)$ is due to the difference between $r_*$ and $r_!$. Hence the isomorphism (\ref{eqn-resolution-*-extension-on-olBunG}) implies
$$ \on{coFib} (  \colim_{ P\in \Par' }  \mbE^{-1} \circ \ol{\Delta}_!\circ \mbK(P) \to \mbE^{-1} \circ \ol{\Delta}_!\circ \mbK(G) )  \simeq \PsId_{\on{naive}}\circ \PsId_!^{-1}[  -2\dim(\BunG)-\dim(Z_G) ].$$
Then we are done because $ \mbE^{-1} \circ \ol{\Delta}_!\circ \mbK\simeq \mbD\mbL$.

\qed[Theorem \ref{thm-main}]

\begin{rem} As a first test for Theorem \ref{thm-main-variant}, let us evaluate the above diagram at $G\in \Par$. By definition, $\mbK(G) \simeq j_{G,!}\circ r_!(k_{\BunG})$. Hence $\ol{\Delta}_!\circ \mbK(G) \simeq \Delta_!(k_{\BunG})$. On the other hand $\mbD\mbL(G) \simeq \Id$, hence $\mbE\circ \mbD\mbL(G) \simeq \Delta_!(k_{\BunG})$ by the definition of $\PsId_!$.
\end{rem}

\begin{rem} The statement of Theorem \ref{thm-main} depends on the miraculous duality on $\BunG$ (i.e., $\PsId_!$ is invertible) but that of Theorem \ref{thm-main-variant} does not. Our proof of the latter will not depend on the miraculous duality either.
\end{rem}

\begin{rem} \label{rem-relation-with-unit-nearby}
The following claim is neither proved nor used in this paper. The object $\mbK(P)$ can be obtained by the following nearby cycles construction. Choose a homomorphism $\gamma:\mBA^1\to T_\ad^+$ connecting the unit point $C_G$ and the point $C_P$. In \cite{schieder2016geometric}, S. Schieder calculated the nearby cycles of the constant D-module for the $\mBA^1$-family
$$ \VinBun_G^\gamma: =\VinBun_G\mt_{T_\ad^+} \mBA^1.$$
Denote this nearby cycles sheaf by $\Psi^\gamma\in \Dmod_{\hol}(\VinBun_G|_{C_P})$. Then up to a cohomological shift, $\mbK(P)$ is isomorphic to the $!$-pushforward of $\Psi^\gamma$ along $\VinBun_G|_{C_P} \to \ol{\Bun}_G$. Moreover, one can use $\Psi^\gamma$ to construct a duality bewteen $\IGP$ and $\mathbf{I}(G,P^-)$. See \cite[Theorem E, Theorem H]{chen2021thesis} for more details.
\end{rem}

\ssec{Organization of this paper}
	The outline of the proof for Theorem \ref{thm-main-variant} will be provided in $\S$ \ref{ssec-proof-main-theorem}. Each other section corresponds to a step in that proof.

	In Appendix \ref{Appendix-corr}, we review the theory of D-modules.

	In Appendix \ref{appendix-well-defined}, we provide proofs for the results mentioned in $\S$ \ref{intro-IGP} (which are due to Gaitsgory).

	In Appendix \ref{appendix-stratification}, we review the gluing functors for D-modules on stratified stacks.

	In Appendix \ref{appendix-wtG}, we prove some results about the group scheme $\wt{G}$, which is the stabilizer of the canonical section $\mfs:T_\ad^+ \to \Vin_G$ for the $(G\mt G)$-action on $\Vin_G$.
	
% !TEX root =  ../main.tex 
\ssec{Notations and conventions}

\begin{notn}[$\infty$-categories]

	We use the theory of $(\infty,1)$-categories developed in \cite{HTT}. We use same symbols to denote a $(1,1)$-category and its simplicial nerve. The reader can distinguish them according to the context.

For two objects $c_1,c_2\in C$ in an $(\infty,1)$-category $C$, we write $\Map_C(c_1,c_2)$ for the \emph{mapping space} between them, which is an object in the homotopy category of spaces. We omit the subscript $C$ if there is no ambiguity.

A \emph{correspondence} from $c_1$ to $c_2$ in $C$ is a diagram $c_2\gets d\to c_1$. The \emph{composition} of two correspondences $c_2\gets d\to c_1$ and $c_3\gets e\to c_2$ is $c_3\gets e\mt_{c_2} d\to c_1$ (when the fiber product $e\mt_{c_2} d$ exists). A $2$-morphism from $c_2\gets d\to c_1$ to $c_2\gets d'\to c_1$ is a morphism $d\to d'$ defined over $c_2$ and $c_1$.

We also need the theory of $(\infty,2)$-categories developed in \cite{GR-DAG1}. For two objects $c_1,c_2\in \mBS$ in an $(\infty,2)$-category $\mBS$, we write $\bMap_\mBS(c_1,c_2)$ for the \emph{mapping $(\infty,1)$-category} between them.
\end{notn}

\begin{defn}[Algebraic geometry] \label{defn-prestack}

Unless otherwise stated, all algebro-geometric objects are defined over a fixed algebraically closed ground field $k$ of characteristic $0$, and are classical (i.e. non-derived).

A \emph{locally finite type prestack} or \emph{lft prestack} is a contravariant (accessible) functor 
$$(\affSch_\ft)^{\op} \to \on{Groupoids}$$
from the category of finite type affine $k$-schemes to the category of groupoids. The collection of them form a $(2,1)$-category $\PreStk_\lft$.

An \emph{algebraic stack} is a lft $1$-Artin stack in the sense of \cite[Chapter 2, $\S$ 4.1]{GR-DAG1}. 

Following \cite{drinfeld2013some}, an algebraic stack $Y$ is \emph{QCA} if it is quasi-compact and the automorphism groups of geometric points are affine. $Y$ is \emph{safe} if it is QCA and the automorphism groups of geometric points are unipotent.

An algebraic stack $Y$ is \emph{QCAD} if it is QCA and the diagonal map $Y\to Y\mt Y$ is affine.

A morphism $Y_1\to Y_2$ between lft prestacks is \emph{stacky} (resp. \emph{QCA}, \emph{QCAD}, \emph{safe}) if it is relatively represented by (resp. QCA, QCAD, safe) algebraic stacks, i.e., for any finite type affine scheme $T$ over $Y_2$, the fiber product $Y_1\mt_{Y_2} T$ is a (resp. QCA, QCAD, safe) algebraic stack.
\end{defn}

\begin{notn} We fix a connected smooth projective curve $X$ over $k$.
\end{notn}

\begin{notn}[DG-categories] \label{notn-DG-categories}

We study \emph{DG-categories} over $k$. Unless otherwise stated, DG-categories are assumed to be \emph{cocomplete} (i.e., containing small colimits), and functors between them are assumed to be $k$-linear and \emph{continuous} (i.e. preserving small colimits). The $(\infty,1)$-category formed by them is denoted by $\DGCat_\cont$. The corresponding $(\infty,2)$-category is denoted by $\bDGCat_\cont$.

$\DGCat_\cont$ carries a closed symmetric monoidal structure, known as the \emph{Lurie tensor product} $\otimes_k$. The unit object for it is the DG-category $\Vect_k$ of $k$-vector spaces. For $\mCC,\mCD\in \DGCat_\cont$, we write $\LFun_k(\mCC,\mCD)$ for the object in $\DGCat_\cont$ characterized by the universal property 
$$\Map(\mCE,\LFun_k(\mCC,\mCD)) \simeq \Map(\mCE\ot_k\mCC,\mCD).$$

A DG-category $\mCM$ is \emph{dualizable} if it is a dualizable object in $\DGCat_\cont$. We write $\mCM^\vee$ for its dual DG-category, which is canonically equivalent to $\LFun_k( \mCM,\Vect_k)$. It is well-known that $\mCM$ is dualizable if it is compactly generated, and there is a canonical identification $\mCM^\vee \simeq \Ind(\mCM^{c,\op})$.
\end{notn}

\begin{notn}[D-modules]

We will use the theory of D-modules on lft prestacks, which is reviewed in Appendix \ref{Appendix-corr}. Here are some notations for this theory.

Let $Y$ be a lft prestack. We write $\Dmod(Y)$ for the DG-category of D-modules on $Y$ We write $\omega_Y$ for the dualizing D-module on $Y$. When $Y$ is an algebraic stack, we write $k_Y$ for the constant D-module. 

Let $f:Y_1\to Y_2$ be a morphism between lft prestacks. We have the \emph{!-pullback functor} $f^!:\Dmod(Y_2)\to \Dmod(Y_1)$.

If $f$ is QCA, we have the \emph{renormalized pushforward functor}, or \emph{$\bt$-pushforward functor}, $f_\bt:\Dmod(Y_1)\to \Dmod(Y_2)$ defined in \cite{drinfeld2013some}. When $f$ is safe, $f_\bt$ coincides with the usual \emph{$*$-pushforward functor}\footnote{In fact, in this paper we only use $\bt$-pushforward functors along safe maps, hence the reader can safely ignore the difference between them and usual $*$-pushforward functor. However, the theory of $\bt$-pushforward functors allow us to prove results for \emph{a priori} non-safe maps} $f_*$.

We have base-change isomorphisms between $\bt$-pushforward functors along QCAD morphisms and all $!$-pullback functors\footnote{The \emph{higher compatibilities} for base-change isomorphisms are claimed in \cite[Remark 9.3.13]{drinfeld2013some} for $\bt$-pushforward functors along any QCA maps. However, we cannot verify their claim unless restricting to QCAD maps. See Appendix \ref{Appendix-corr} for more details.}. For any correspondence of lft prestacks $\alpha: (Y_2\overset{f}\gets Z\overset{g}\to Y_1)$ such that $f$ is QCAD, we write $ \DMod^{\bt,!}(\alpha):= f_\bt\circ g^!:\Dmod(Y_1)\to \Dmod(Y_2)$. The assignment $\alpha\mapsto \DMod^{\bt,!}(\alpha)$ is compatible with compositions of correspondences. For any $2$-morphism from $\alpha: (Y_2\overset{f}\gets Z\overset{g}\to Y_1)$ to $\alpha': (Y_2\overset{f}\gets Z'\overset{g}\to Y_1)$ such that $Z\to Z'$ is an open embedding, we have a natural transformation $\DMod^{\bt,!}(\alpha')\to \DMod^{\bt,!}(\alpha')$.

The left adjoints of $!$-pullback and $*$-pushforward functors are \emph{not} well-defined in general. Hence sometimes we restrict to the full subcategory $\Dmod_\hol(Y)\subset \Dmod(Y)$ of ind-holonomic D-modules on $Y$. For a morphism $f:Y_1\to Y_2$ between lft prestacks, we have the $!$-pushforward functor $f_!:\Dmod_\hol(Y_1)\to \Dmod_\hol(Y_2)$ which is left adjoint to the restriction of $f^!$. If $f$ is stacky, we have the $*$-pullback functor $f^*: \Dmod_\hol(Y_2)\to \Dmod_\hol(Y_1)$. If $f$ is safe, $f^*$ is left adjoint to the restriction of $f_*\simeq f_\bt$.

We have base-change isomorphisms between all $!$-pushforward functors and $*$-pullback functors along stacky maps. For any correspondence of lft prestacks $Y_2\overset{f}\gets Z\overset{g}\to Y_1 $ such that $g$ is stacky, we write $ \DMod^{!,*}(\alpha):= f_!\circ g^*:\Dmod_\hol(Y_1)\to \Dmod_\hol(Y_2)$. The assignment $\alpha\mapsto \DMod_\hol^{!,*}(\alpha)$ is compatible with compositions of correspondences. For any $2$-morphism from $\alpha: (Y_2\overset{f}\gets Z\overset{g}\to Y_1)$ to $\alpha': (Y_2\overset{f}\gets Z'\overset{g}\to Y_1)$ such that $Z\to Z'$ is an open embedding, we have a natural transformation $\DMod^{!,*}(\alpha)\to \DMod^{!,*}(\alpha')$.

\end{notn}

\begin{defn}[Stratifications] \label{defn-stratified-space}
Let $Y$ be an algebraic stack and $I$ be a finite set. A \emph{stratification of $Y$ labelled by the power poset $P(I)$} is an assignment of open substacks $U_i\subset Y$ for any $i\in I$. For such an assignment, we define $i_J:Y_J\to Y$ to be the reduced locally closed substack of $Y$ given by
$$ (\bigcap_{ j\in J } U_j)\setminus (\bigcup_{ i\notin J } U_i).$$
We call $Y_J$ \emph{the stratum labelled by $J$}. Every geometric point of $Y$ is contained in exactly one stratum.

For any object $J\in P(I)$, there is a unique open substack $Y_{\ge J}\subset Y$ whose geometric points are exactly those contained in $\bigcup_{K\supset J} Y_K  = \bigcap_{ j\in J } U_j$. Similarly, we define the reduced closed substack $Y_{\le J}$.

The stratum $Y_I$ is an open substack of $Y$. Hence we also write $j_I:=i_I$ for this open embedding.

\end{defn}

\begin{convn} Sometimes we also use the notations $Y_J$ and $Y_{\le J}$ for certain stacks whose underlying reduced stacks are defined above. For example, if we have a map $Y\to Z$ and a stratification of $Z$ labelled by $P(I)$, then we obtain a stratification of $Y$ labelled by $P(I)$ by pulling back the open substacks. We often write $Y_J:= Y\mt_Z Z_J$ although it is not necessarily reduced. This ambiguity is not a problem because the category of D-modules on a stack only depends on the underlying reduced stack.

\end{convn}

\begin{defn}[Coordinate stratification] Let $Y$ be a finite type scheme and $\{f_i\}_{i\in I}$ be regular functions on $Y$. Then we obtain a stratification of $Y$ labelled by $P(I)$ with $U_i$ given by the non-vanishing locus of $f_i$. In particular, the coordinate functions induce a stratification of the affine space $\mBA^I$ labelled by $P(I)$. This stratification is known as the \emph{coordinate stratification}.
\end{defn}

\begin{notn}[Reductive groups] \label{notn-convension-reductive-group}

We fix a connected reductive group $G$. For simplicity, we assume $[G,G]$ to be simply connected\footnote{Such assumption was made in many references that we cite, but our proof depends on such assumption only because they use the geometric constructions in \cite{braverman2002geometric}. However, \cite[Section 7]{schieder2015harder} explained how to modify such constructions such that they work for any reductive group. Applying such modifications, our proof works for any reductive group.}. 

We fix a pair of \emph{opposite Borel subgroups} $(B,B^-)$ of it, therefore a \emph{Cartan subgroup} $T$. We write $Z_G$ for the center of $G$ and $T_\ad:=T/Z_G$ for the \emph{adjoint torus}.

We write $\mCI$ for the \emph{set of vertices in the Dynkin diagram} of $G$, $\Lambda_G$ (resp. $\check{\Lambda}_G$) for the \emph{coweight (resp. weight) lattice}, and $\Lambda_G^\pos\subset \Lambda_G$ for the sub-monoid spanned by all positive simple co-roots $(\alpha_i)_{i\in \mCI}$.

We often use $P$ to denote a \emph{standard parabolic subgroup} of $G$ (i.e. a parabolic subgroup containing $B$). We write $P^-$ for the corresponding \emph{standard opposite parabolic subgroup} and $M:=P\cap P^-$ for the \emph{Levi subgroup}. We write $U$ (resp. $U^-$) for the \emph{unipotent radical} of $P$ (resp. $P^-$). When we need to use a different parabolic subgroup, we often denote it by $Q$ and its Levi subgroup by $L$.

We write $\Par$ for the partially ordered set of standard parabolic subgroups of $G$. We write $\Par'=\Par\setminus\{G\}$. We view them as categories in the standard way.
\end{notn}

\begin{notn}[Semi-group completion] \label{notn-convension-semi-group-completion}

The collection of simple positive roots of $G$ provides an identification $T_\ad\simeq \mBG_m^\mCI := \prod_{i\in \mCI} \mBG_m$. We define $T_\ad^+:=\mBA^\mCI\supset \mBG_m^\mCI\simeq T_\ad$, which is a semi-group completion of the adjoint torus $T_\ad$.

Consider the coordinate stratification of the affine space $T_\ad^+:=\mBA^\mCI$. The poset of this stratification can be identified with $\Par$ such that for $P\subset \Par$ corresponding to a subset $\mCJ\subset \mCI$ of the Dynkin diagram, the $P$-stratum is given by $\mBG_m^\mCJ \mt \{0\}^{\mCI\setminus \mCJ}\to \mBA^\mCI$. In particular, we have $T_{\ad,P}^+\simeq \mBG_m^\mCJ \mt \{0\}^{\mCI\setminus \mCJ}$, $T_{\ad,\ge P}^+\simeq \mBG_m^\mCJ \mt \mBA^{\mCI\setminus \mCJ}$ and $T_{\ad,\le P}^+\simeq \mBA^\mCJ \mt \{0\}^{\mCI\setminus \mCJ}$

Write $C_P$ for the unique point in $T_{\ad,P}^+$ whose coordinates are either $0$ or $1$. In particular $C_B$ is the zero element in $T_\ad^+$ and $C_G$ is the unit element. We use the same symbols to denote the images of these points in the quotient stack $T_{\ad}^+/T$.

Consider the homomorphism $Z_M/Z_G \to T_\ad$. Let\footnote{It was denoted by $T^+_{\ad,\ge P,\on{strict}}$ in \cite{schieder2016geometric}.} $T^+_{\ad,\ge C_P}$ be its closure in $T_\ad^+$. In other words, we have $T^+_{\ad,\ge C_P}\simeq \{1\}^\mCJ \mt \mBA^{\mCI\setminus \mCJ}$

Note that it is a sub-semi-group of $T^+_{\ad,\ge P}$ that contains $C_P$ as an idempotent element.
\end{notn}

\textbf{Acknowledgements:} I want to thank my advisor Dennis Gaitsgory for teaching me all the important concepts in this paper, such as the pseudo-identity functor, the Vinberg-degeneration, Braden's theorem, etc.. I am also grateful for his sharing of notes on the category $\IGP$ and his comments on the first draft of this paper.

I am very grateful to the anonymous referees for their numerous comments and suggestions.

\section{Outline of the proof} 
\label{ssec-proof-main-theorem}

In this subsection, we reduce Theorem \ref{thm-main-variant} to a series of partial results, which will be proved in the later sections.

\ssec{Step 1: constructing the natural transformation}

The first step is to construct a natural transformation from $\ol{\Delta}_!\circ \mbK$ to $\mbE\circ \mbD\mbL$. Let us first explain how to construct the morphism
\begin{equation} \label{nt-level-objects}
  \ol{\Delta}_!\circ \mbK(P) \to \mbE\circ \mbD\mbL(P).
 \end{equation}
For $P\in \Par$, let $\ol{\Bun}_{G,P}$ be the $P$-stratum of $\ol{\Bun}_G$. We will construct (see Proposition-Construction \ref{propconstr-factor-through-gen} and the remark below it) a canonical commutative diagram\footnote{In the case $G=\on{SL}_2$, using the notations in Footnote \ref{Footnote:SL2-olBunG}, $\VinBun_{G,B}$ classifies (up to nil-isomorphisms) objects ($\phi:E_1\to E_2$) with $\det(\phi)=0$. It follows from the definition that the subsheaf $\on{im}(\phi)$ is a generic line bundle. Then the map $\on{im}(\phi)\to E_2$ provides a generic $B$-reduction to the $\on{SL}_2$-torsor for $E_2$. This provides a map $\VinBun_{G,B}\to \Bun_G^{B\hgen}$ that factors through the quotient $\VinBun_{G,B}/T=\ol{\Bun}_{G,B}$.}
\begin{equation} \label{eqn-begining-of-everything}
\xyshort
\xymatrix{
	\ol{\Bun}_{G,P} \ar[r]^-{ \ol{\Delta}_P^{\enh,l} } \ar[d]^-{i_P} &
	\Bun_G^{P\hgen}\mt \BunG \ar[d]^-{  \mbp^\enh_{P\mt G\to G\mt G} } \\
	\ol{\Bun}_G \ar[r]^-{\ol{\Delta}} & \BunG\mt \BunG.
}
\end{equation}
Consider the object
$$\mCF_P:=\ol{\Delta}^{\enh,l}_{P,!} \circ i_{ P}^!\circ \mbK(P) \in \Dmod_\hol(  \Bun_G^{P\hgen}\mt \BunG  ). $$
Note that 
\begin{equation} \label{eqn-olDel-K=Eis-mCF-pre}
\ol{\Delta}_! \circ \mbK(P) \simeq   \mfp^\enh_{P\mt G\to G\mt G,!}(\mCF_P)
\end{equation}
because $\mbK(P)\simeq i_{P,!}\circ i_P^!(\mbK(P))$.

The following result will be proved in $\S$ \ref{ssec-step-1-2}:

\begin{prop} \label{prop-FP-equivariant}
The object $\mCF_P:= \ol{\Delta}^{\enh,l}_{P,!} \circ i_{ P}^!\circ \mbK(P) \in \Dmod_\hol(  \Bun_G^{P\hgen}\mt \BunG  )$ is contained in the full subcategory
$$  \mbI(G\mt G,P\mt G) \subset \Dmod(  \Bun_G^{P\hgen}\mt \BunG  ).$$
\end{prop}

Let $\mCF_P'$ be the corresponding object in $\mbI(G\mt G,P\mt G)$. By (\ref{eqn-olDel-K=Eis-mCF-pre}), we have
\begin{equation} \label{eqn-olDel-K=Eis-mCF}
\ol{\Delta}_!\circ \mbK(P) \simeq \Eis^\enh_{ P\mt G\to G\mt G }(\mCF_P').
\end{equation}
Hence by functoriality of the LHS, we obtain a morphism
\begin{equation} \label{eqn-Eis-mCF-to-Del}
 \Eis^\enh_{ P\mt G\to G\mt G }(\mCF_P') \to \mCF_G'.
\end{equation}
By adjunction, we have a morphism
\begin{equation} \label{eqn-mCF-to-CT-Del}
  \theta_P': \mCF_P'\to \CT^\enh_{G\mt G\gets P\mt G}(\mCF_G' ).
 \end{equation}
Note that $\mCF_G'=\mCF_G\simeq \Delta_!(k_{\BunG})$.

On the other hand, by definition we have
\begin{equation}\label{eqn-E-DL=to-Eis-CT-mCFG}
 \mbE\circ \mbD\mbL(P) \simeq ((\Eis_{P\to G}^\enh\circ \CT_{G\gets P}^\enh)\ot \Id)( \Delta_!(k_{\Bun_G}) ) \simeq  \Eis^\enh_{ P\mt G\to G\mt G } \circ \CT^\enh_{ G\mt G\gets P\mt G }(\Delta_!(k_{\BunG})),
\end{equation}
where the last isomorphism is because $\mbI(G,P)\ot_k \Dmod(\Bun_G) \simeq \mbI(G\mt G,P\mt G)$ (see Lemma \ref{lem-double-the-group} below). Now we declare the morphism (\ref{nt-level-objects}) to be given by
\begin{equation} \label{eqn-value-natural-transformation}
 \ol{\Delta}_!\circ \mbK(P) \simeq \Eis^\enh_{ P\mt G\to G\mt G }(\mCF_P') \xto{ \Eis^\enh(\theta_P') }  \Eis^\enh_{ P\mt G\to G\mt G } \circ \CT^\enh_{ G\mt G\gets P\mt G }(\mCF_G') \simeq  \mbE\circ \mbD\mbL(P). 
 \end{equation}

In order to obtain the desired natural transformation, we need the following construction:

\begin{propconstr} \label{propconstr-factor-through-gen}
	Let $\ol{\Bun}_{G,\le P}$ be the reduced closed substack of $\ol{\Bun}_G$ containing all the $P'$-strata with $P'\subset P$. Then there exist canonical maps
	$$\ol{\Delta}^\enh_{\le P}:\ol{\Bun}_{G,\le P} \to \Bun_G^{P\hgen}\mt \BunG^{P^-\gen}$$
	that are functorial\footnote{Note that for any $P\subset Q$, we have maps 
	$$\ol{\Bun}_{G,\le P}\to \ol{\Bun}_{G,\le Q}\;\on{ and }\; \Bun_G^{P\hgen}\mt \BunG^{P^-\hgen} \to \Bun_G^{Q\hgen}\mt \BunG^{Q^-\hgen}.$$} in $P$ such that when $P=G$ we have $\ol{\Delta}^\enh_{\le G}=\ol{\Delta}: \ol{\Bun}_G \to \BunG\mt \BunG $.
\end{propconstr}

\begin{rem}
In particular, we have functorial maps
$$\ol{\Delta}^{\enh,l}_{\le P}:\ol{\Bun}_{G,\le P} \to \Bun_G^{P\hgen}\mt \BunG.$$
The map $\ol{\Delta}^{\enh,l}_{P}$ in (\ref{eqn-begining-of-everything}) is defined to be its restriction to the $P$-stratum. 
\end{rem}

Note that we also have
$$\mCF_P \simeq \ol{\Delta}^{\enh,l}_{\le P,!} \circ i_{\le P}^!\circ \mbK(P).$$
Using this, it is purely formal to show that the morphisms (\ref{eqn-value-natural-transformation}) constructed above are functorial in $P$. Namely, in $\S$ \ref{ssec-step-1-3}, we will use the theory of (co)Cartesian fibrations to prove:

\begin{propconstr} \label{propconstr-nt-cartesian}
There exists a canonical natural transformation $\ol{\Delta}_!\circ \mbK\to \mbE\circ \mbD\mbL$ whose value at $P\in \Par$ is equivalent to the morphism (\ref{eqn-value-natural-transformation}).
\end{propconstr}

\ssec{Step 1.5: reducing to compare two objects}

After obtaining the natural transformation, we only need to show its value at each $P\in \Par$ is invertible. From this step on, we fix such a standard parabolic $P$.

Unwinding the definitions, we need to show the morphism (\ref{eqn-mCF-to-CT-Del})
$$\theta_P': \mCF_P'\to \CT^\enh_{G\mt G\gets P\mt G}(\mCF_G' )$$
 is invertible. Recall (see Remark \ref{rem-conservative-iota^!}) that the functor
$$ \iota_{M\mt G}^!: \mbI(G\mt G,P\mt G) \to \Dmod( \BunM \mt \BunG )$$
is conservative. Hence we only need to show the map $\iota_{M\mt G}^!(\theta'_P)$ is invertible. By definition, $\iota_{M\mt G}^!$ is isomorphic to
$$ \mbI(G\mt G,P\mt G) \to \Dmod(\Bun_G^{P\hgen}\mt \BunG) \xto{\iota^!_{P\mt G}}  \Dmod(\BunP\mt \BunG) \xto{\mfq_{P\mt G,*}} \Dmod(\BunM\mt \BunG).$$
We denote the composition of the latter two functors by
$$ \CT^{\gen}_{P\mt G,*}: \Dmod(\Bun_G^{P\hgen}\mt \BunG) \xto{\iota^!_{P\mt G}}  \Dmod(\BunP\mt \BunG) \xto{\mfq_{P\mt G,*}} \Dmod(\BunM\mt \BunG). $$
Then the source of $\iota_{M\mt G}^!(\theta'_P)$ is isomorphic to $\CT^{\gen}_{P\mt G,*}(\mCF_P)$.

On the other hand, the functor $\iota_{M\mt G}^!\circ \CT^\enh_{G\mt G\gets P\mt G}$ is isomorphic to the usual geometric constant term functor
$$ \CT_{P\mt G,*}: \Dmod(\BunG\mt \BunG) \to \Dmod(\BunM\mt \BunG)$$
(as can be seen by passing to left adjoints). Hence the target of $\iota_{M\mt G}^!(\theta'_P)$ is isomorphic to $\CT_{P\mt G,*}(\mCF_G)$.
Let\footnote{We will give a more direct description of $\gamma_P$ in $\S$ \ref{sssec-description-gammaP}.}
\begin{equation} \label{eqn-gamma-P}
 \gamma_P:  \CT^{\gen}_{P\mt G,*}(\mCF_P) \to \CT_{P\mt G,*}(\mCF_G).
\end{equation}
be the morphism obtained from $\iota_{M\mt G}^!(\theta'_P)$ via the above isomorphisms. Then we have reduced the main theorem to:
\begin{itemize}
  \item We only need to show $\gamma_P$ is invertible.
\end{itemize}

\ssec{Step 2: translating by the second adjointness}

Recall that the main theorem of \cite{drinfeld2016geometric} says that when restricted to each connected component $\Bun_{M,\lambda}$ of $\BunM$, the functor 
$$ \CT_{P,*,\lambda}:\Dmod(\BunG) \xto{!\on{-pull}} \Dmod(\Bun_{P,\lambda}) \xto{*\on{-push}} \Dmod(\Bun_{M,\lambda}) $$
is canonically left adjoint to 
$$\Eis_{P^-,*,\lambda}: \Dmod(\Bun_{M,\lambda}) \xto{!\on{-pull}} \Dmod(\Bun_{P^-,\lambda}) \xto{*\on{-push}}  \Dmod(\BunG).$$

In particular, the functor $\CT_{P,*}\simeq \bigoplus \CT_{P,*,\lambda}$ preserves ind-holonomic objects and its restriction to $\Dmod_\hol(\BunG)$ is isomorphic to
$$ \CT_{P^-,!}: \Dmod_\hol( \BunG ) \xto{\mfp^{-,*}} \Dmod_\hol(\BunPm)  \xto{\mfq^-_{!}} \Dmod_\hol(\BunM).$$
Hence we obtain an isomorphism
$$ \CT_{P\mt G,*}(\mCF_G) \simeq \, \CT_{P^-\mt G,!}(\mCF_G).$$
As explained in \loccit, this result is similar to Braden's theorem on hyperbolic localizations (\!\!\cite{braden2003hyperbolic}) because $P$, $P^-$ and $M$ are respectively the attractor, repeller and fixed loci for a certain $\mBG_m$-action on $G$.

Now there is a similar story when we replace $\BunG$ by $\BunG^{P\hgen}$. To explain it, we view $\BunG^{P\hgen}$ as applying $\Bun_\bullet^{\bullet\hgen}$ to the pair $(G,P)$. Then the attractor, repeller and fixed loci for this pair are respectively $(P,P)$, $(P^-,M)$ and $(M,M)$. Therefore one would expect the main theorem of \cite{drinfeld2016geometric} holds if $\BunG$, $\BunP$, $\BunPm$ and $\BunM$ are replaced by $\BunG^{P\hgen},\BunP,\BunPm^{M\hgen},\BunM$. We will prove the following result in $\S$ \ref{section-2nd-adjointness}, which says this expectation is correct:

\begin{thm} \label{thm-generic-second-adjointness} The functor 
\[
\CT_{P,*,\lambda}^\gen:\Dmod(\BunG^{P\hgen}) \xto{!\on{-pull}} \Dmod(\Bun_{P,\lambda}) \xto{*\on{-push}} \Dmod(\Bun_{M,\lambda})
\]
is canonically left adjoint to
\[
\Dmod(\Bun_{M,\lambda}) \xto{!\on{-pull}} \Dmod(\Bun_{P^-,\lambda}^{M\hgen}) \xto{*\on{-push}}  \Dmod(\BunG^{P\hgen}).
\]
Therefore the restriction of $\CT_{P,*}^\gen$ on ind-holonomic objects is isomorphic to
\[\CT_{ P^-,! }^\gen: \Dmod_\hol(   \BunG^{P\hgen}  ) \xto{*\on{-pull}}  \Dmod_\hol(\BunPm^{M\hgen} ) \xto{!\on{-push}}  \Dmod_\hol(  \BunM ).\]
\end{thm}

\begin{rem} Part of the statement says all the functors in Theorem \ref{thm-generic-second-adjointness} are well-defined. For example, we need to show the map $\Bun_{P^-,\lambda}^{M\hgen}\to \BunG^{P\hgen}$ is safe.
\end{rem}

\begin{rem} Unlike the main theorem of \cite{drinfeld2016geometric}, Theorem \ref{thm-generic-second-adjointness} does \emph{not} imply the similar result when $\BunP$ and  $\BunPm^{M\hgen}$ are exchanged.

\end{rem}

Applying the above theorem to the reductive group $G\mt G$ and its parabolic subgroup $P\mt G$, we obtain an isomorphism
$$ \CT_{P\mt G,*}^\gen(\mCF_G) \simeq \, \CT_{P^-\mt G,!}^\gen(\mCF_G).$$
Hence the morphism (\ref{eqn-gamma-P}) is equivalent to a certain morphism
$$ \gamma'_P: \,\CT_{ P^-\mt G,! }^\gen(\mCF_P) \to \,\CT_{P^-\mt G,!}(\mCF_G).$$
Hence we have reduced the main theorem to the following problem:
\begin{itemize}
	\item We only need to show $\gamma'_P$ is invertible.
\end{itemize}

\begin{rem} It is easier to study $\gamma'_P$ than $\gamma_P$ because we can use the base-change isomorphims.
\end{rem}

\begin{rem} We believe Theorem \ref{thm-generic-second-adjointness} (and its proof) is of independent interest. For example, we can use them to give a description of the monad structure of $\wt{\Omega}_P:= \iota_{M}^!\circ \iota_{M,!}$ via Verdier (co)specialization along Schieder's local models. This monad was the central concept in the paper \cite{gaitsgory2011acts}. The details of it will be provided elsewhere.
\end{rem}

\ssec{Step 3: diagram chasing}

Using the base-change isomorphisms, and using the facts that $\mbK(P)$ is a $!$-extension along $\ol{\Bun}_{G,P}\to \ol{\Bun}_G$, one can simplify the source and target of $\gamma'_P$. Let us state the result directly\footnote{The result below only serves as motivation and will be incorporated into Lemma \ref{lem-diagram-chasing}.}.

Consider the correspondences
\begin{eqnarray*}
 \beta_P: (\BunM\mt \BunG  &\gets \BunPm^{M\hgen} \mt_{  \BunG^{P\hgen} } \ol{\Bun}_{G,P}  \to& \ol{\Bun}_{G})\\
\beta_{G}:( \BunM\mt \BunG  &\gets \BunPm \mt_{  \BunG } \ol{\Bun}_{G,G}  \to& \ol{\Bun}_{G}),
\end{eqnarray*}
where the map $\ol{\Bun}_{G,P}\to \BunG^{P\hgen} $ was contructed in Proposition-Construction \ref{propconstr-factor-through-gen}, and the left arm of $\beta_P$ is given by
$$ \BunPm^{M\hgen} \mt_{  \BunG^{P\hgen} } \ol{\Bun}_{G,P}   \to \BunM \mt \BunG\mt \BunG \xto{\on{pr}_{13}}   \BunM\mt \BunG. $$
Then the base-change isomorphisms provide
\begin{eqnarray*}
\CT_{ P^-\mt G,! }^\gen(\mCF_P) &\simeq& \DMod_\hol^{!,*}( \beta_P ) \circ \mbK(P) ,\\
\CT_{P^-\mt G,!}(\mCF_G) &\simeq& \DMod_\hol^{!,*}( \beta_{G} ) \circ \mbK(G).
\end{eqnarray*}

This motivates the following construction (see $\S$ \ref{ssec-step-3-1}):

\begin{propconstr} \label{prop-constr-BunPm-olBunG-gen}
	There exists an open substack\footnote{\label{fn-Pm-geP-gen-SL2}
In the case $G=\on{SL}_2$, recall that $\Bun_{B^-}\mt_{\BunG} \ol{\Bun}_{G}$ classifies certain chains $E_1\to E_2\to L_2$. Then the desired open substack classifies those chains such that the restriction of $E_1\to L_2$ at any geometric point of $S$ is nonzero.} 
	$$ (\BunPm\mt_{\BunG} \ol{\Bun}_{G,\ge P})^\gen \subset \BunPm\mt_{\BunG} \ol{\Bun}_{G,\ge P}$$
	such that the parameterized correspondence
	$$ 
	\xyshort
	\xymatrix{
		\beta: & \BunM\mt \BunG  &
		(\BunPm\mt_{\BunG} \ol{\Bun}_{G,\ge P})^\gen  \ar[r] \ar[l]\ar[d] &
		\ol{\Bun}_{G} \\
		& & T_{\ad,\ge P}^+/T
	}
 	$$
 	captures $\beta_P$ (resp. $\beta_G$) as its restriction to the $P$-stratum (resp. $G$-stratum) of $T_{\ad,\ge P}^+/T$.
\end{propconstr}

Using the fact that $\mbK(P)$ is a $!$-extension along $\ol{\Bun}_{G,P}\to \ol{\Bun}_G$ again, we obtain isomorphisms\footnote{See (\ref{eqn-Dhol-ult}) for the notation $\DMod_\hol^{!,*}$.}
\begin{eqnarray*}
\DMod_\hol^{!,*}( \beta_P ) \circ \mbK(P) &\simeq& \DMod_\hol^{!,*}( \beta ) \circ \mbK(P) ,\\
\DMod_\hol^{!,*}( \beta_{G} ) \circ \mbK(G) &\simeq& \DMod_\hol^{!,*}( \beta ) \circ \mbK(G).
\end{eqnarray*}
We will prove the following result in $\S$ \ref{ssec-step-3-2} by a routine diagram chasing:

\begin{lem} \label{lem-diagram-chasing}
	The morphisms $\gamma_P$ and $\gamma'_P$ are both equivalent to the morphism
	$$ \DMod_\hol^{!,*}( \beta ) \circ \mbK(P)\to \DMod_\hol^{!,*}( \beta ) \circ \mbK(G)  $$
	given by the functor $\DMod_\hol^{!,*}( \beta ) \circ \mbK: \Par \to \Dmod(\BunM\mt \BunG)$.
\end{lem}

Hence we have reduced the main theorem to the following problem:
\begin{itemize}
	\item We only need to show the functor $\DMod_\hol^{!,*}( \beta ) \circ \mbK$ sends the arrow $P\to G$ to an isomorphism.
\end{itemize}

\ssec{Step 4: restoring the symmetry}

In the previous step, we have reduced the main theorem to proving a morphism in $\Dmod(\BunM\mt \BunG)$ is invertible. Note that we have broken the symmetry in the sense that the two factors $\BunM$ and $\BunG$ are different. In this step, we restore the symmetry and reducing the main theorem to proving a morphism in $\Dmod(\BunM\mt \BunM)$ is invertible. In fact, we think the previous steps are more or less natural to people familiar with this field, but it is the desire to restore this symmetry that led us to discover the proof of the main theorem.

In $\S$ \ref{ssec-step-4-1}, we will show
\begin{propconstr} \label{propconstr-step-4-1}
	There exists a canonical factorization of the map\footnote{In the case $G=\on{SL}_2$, the map $(\Bun_{B^-}\mt_{\BunG} \ol{\Bun}_{G})^\gen \to  \Bun_G^{B^-\hgen}$ sends a chain $E_1\to E_2\to L_2$ in Footnote \ref{fn-Pm-geP-gen-SL2} to the generic $B^-$-reduction provided by the map $E_1\to L_2$.} 
	$$(\BunPm\mt_{\BunG} \ol{\Bun}_{G,\ge P})^\gen \to \BunM\mt \BunG$$
	via $\BunM\mt \Bun_G^{P^-\hgen}$.
\end{propconstr}

In particular we obtain a correspondence
$$ \beta' : ( \BunM\mt \Bun_G^{P^-\hgen}  \gets (\BunPm\mt_{\BunG} \ol{\Bun}_{G,\ge P})^\gen \to \ol{\Bun}_G  )$$
and we only need to show $\DMod_\hol^{!,*}( \beta' ) \circ \mbK$ sends the arrow $P\to G$ to an isomorphism.

The following result will be proved in $\S$ \ref{ssec-step-4-2}:
\begin{prop} \label{prop-step-4-2}
The objects $\DMod_\hol^{!,*}( \beta' ) \circ \mbK(P)$ and $\DMod_\hol^{!,*}( \beta' ) \circ \mbK(G)$ are both contained in the full subcategory
$$ \mbI(M\mt G, M\mt P^- ) \subset \Dmod( \BunM\mt \Bun_G^{P^-\hgen} ). $$
\end{prop}
Consider the functor
\[  \CT_{M\mt P^-,*}^\gen:\Dmod( \BunM\mt \Bun_G^{P^-\hgen} ) \xto{!\on{-pull}} \Dmod(\BunM\mt \BunPm) \xto{*\on{-push}} \Dmod(\BunM\mt \BunM).\]
By definition, its restriction to $\mbI(M\mt G, M\mt P^- ) $ is the conservative functor $\iota^!: \mbI(M\mt G, M\mt P^- )\to \Dmod(\BunM\mt \BunM)$. Hence we only need to show $ \CT_{M\mt P^-,*}^\gen\circ \DMod_\hol^{!,*}( \beta' ) \circ \mbK$ sends the arrow $P\to G$ to an isomorphism.

Now consider the correspondence
$$\delta:( \BunM\mt \BunM \gets  \BunM\mt \Bun_P^{M\hgen}  \to   \BunM\mt \Bun_G^{P^-\hgen})$$
and the functor
$$ \CT_{M\mt P,!}^\gen:= \Dmod_\hol^{!,*}( \delta ).$$
Applying Theorem \ref{thm-generic-second-adjointness} to the reductive group $M\mt G$ and its parabolic subgroup $M\mt P^-$, we obtain $\CT_{M\mt P^-,*}^\gen\simeq \CT_{M\mt P,!}^\gen$. Hence we have reduced the main theorem to the following problem:
\begin{itemize}
	\item We only need to show the functor $ \DMod_\hol^{!,*}(\delta\circ \beta' ) \circ \mbK$ sends the arrow $P\to G$ to an isomorphism.
\end{itemize}

\ssec{Step 5: calculating via the local models}

Now comes the fruit of the symmetry restored in Step 4. Namely, the correspondence $\delta\circ \beta'$ was in fact studied in the literature. To describe it, we need to use \emph{Vinberg semi-group} $\Vin_G$. We refer the reader to \cite[Appendix D]{drinfeld2016geometric}, \cite{wang2017reductive}, \cite[\S 1.2]{nearby-unit} for details about it. To continue our proof, it is enough to know the following:
\begin{itemize}
  \item $\Vin_G$ is an affine normal semi-group equipped with a flat homomorphism to $T_\ad^+$.
  \item $\Vin_G$ is equipped with a $T$-action compatible with the $T$-action on $T_\ad^+$.
  \item $\Vin_G$ is equipped with a $(G,G)$-action which preserves the fibers of $\Vin_G\to T_\ad^+$.
  \item $\Vin_G\to T_\ad^+$ is equipped with a semi-group splitting $\mfs:T_\ad^+\to \Vin_G$.
  \item The open subgroup of invertible elements of $\Vin_G$ can be identified with $G_\enh \simeq (G\mt T)/Z_G$ in a way compatible with the $T$-action, $(G,G)$-action and the projection to $T/Z_G\subset T_\ad^+$. Via this identification, the section $\mfs$ sends $t\in T/Z_G$ to $(t^{-1},t)\in (G\mt T)/Z_G$.
  \item The $(G,G)$-orbit of the section $\mfs$ is an open subscheme $_0\,\Vin_G\subset \Vin_G$.
  \item In the case $G = \on{SL}_2$, $\Vin_G$ is the semi-group of $2\mt 2$ matrices.
\end{itemize}
We also need Schieder's (relative) local model\footnote{\label{fn-local-model-sl2}
In the case $G=\on{SL}_2$, $Y^B_\rel$ classifies chains $L_1\to E_1\to E_2\to L_2$ where $L_1\to E_1$, $E_1\to E_2$ and $E_2\to L_2$ are respectively $S$-points of $\Bun_B$, $\VinBun_G$ and $\Bun_{B^-}$ such that the restriction of $L_1\to L_2$ at any geometric point of $S$ is nonzero.} in \cite{schieder2016geometric}:
$$ Y^P_\rel:= \bMap_\gen(X, P^-\backslash \Vin_{G,\ge C_P}/P \supset P^-\backslash \Vin^\Bru_{G,\ge C_P}/P ),$$
where (see Notation \ref{notn-convension-semi-group-completion})
$$ \Vin_{G,\ge C_P}:=  \Vin_G\mt_{T_\ad^+} T_{\ad,\gCP}^+$$
and $\Vin^\Bru_{G,\ge C_P}\subset \Vin_{G,\ge C_P}$ is the open subscheme $P^- \mfs(   T_{\ad,\gCP}^+)P$.

The $T$-action on $G\backslash \Vin_{G,\ge P}/G$ induces a $Z_M$-action on $P^-\backslash \Vin_{G,\ge C_P}/P$ which preserves $P^-\backslash \Vin^\Bru_{G,\ge C_P}/P$, hence we obtain a $Z_M$-action on $Y^P_\rel$. Note that we have the following diagram of stacks equipped with group actions:
$$  (\pt \act \BunM\mt \BunM) \gets ( Z_M \act Y^P_\rel  ) \to (T\act \VinBun_{G}).$$
By taking quotients, we obtain
\[ \BunM\mt \BunM \gets Y^P_\rel/Z_M \to \ol{\Bun}_G. \]

The following result is proved in $\S$ \ref{ssec-step-4-1}.

\begin{lem} \label{lem-appear-local-model}
The composition $\delta\circ \beta'$ is isomorphic to
$$ \BunM\mt \BunM \gets Y^P_\rel/Z_M \to \ol{\Bun}_G.$$
\end{lem}

It is known (see Construction \ref{constr-project-local-model-hecke}) that $Y^P_\rel/Z_M\to \BunM\mt \BunM$ factors via $H_\MGPos/Z_M$, where $H_\MGPos$ is the $G$-positive part of Hecke stack for $M$-torsors\footnote{In the case $G=\on{SL}_2$, $H_{T,G\on{-pos}}$ classifies morphisms between line bundles $L_1\to L_2$ whose restriction at any geometric point of $S$ is nonzero. The map $Y^P_\rel \to H_\MGPos$ sends the chain $L_1\to E_1\to E_2\to L_2$ in Footnote \ref{fn-local-model-sl2} to $L_1\to L_2$.}:
$$H_\MGPos := \bMap_\gen(X, M\backslash \ol{M}/M \supset M\backslash M/M).$$
Here $\ol{M}\simeq \mfs(C_P)\cdot \Vin_{G,C_P}\cdot \mfs(C_P)$ is a closed subscheme of $\Vin_{G,C_P}$ (see \cite{wang2017reductive}).

Consider the correspondence
$$ \psi_P: H_\MGPos/Z_M \gets  Y^P_\rel/Z_M \to \ol{\Bun}_{G}.$$
We have reduced the main theorem to 
\begin{goal}  \label{prop-finish}
The functor $\DMod_\hol^{!,*}(  \psi_P)\circ\mbK $ sends the arrow $P\to G$ to an isomorphism.
\end{goal}

We will prove a stronger result:

\begin{goal} For any $Q\in \Par_{\ge P}$, the functor $\DMod_\hol^{!,*}(  \psi_P)\circ\mbK $ sends the arrow $Q\to G$ to an isomorphism.
\end{goal}

We prove this by induction on the relative rank between $Q$ and $G$. When $Q=G$, there is nothing to prove. Hence we assume $Q\ne G$ and assume the above claim is correct for any $Q'$ strictly greater than $Q$. Let $L$ be the Levi subgroup of $Q$.

Consider the object
$$\mCD_{Q}:= \on{coFib}( \colim_{Q'\in \Par'\cap \Par_{\ge Q} } \mbK(Q') \to \mbK(G) ) $$
We claim
\begin{equation} \label{eqn-claim-DQ}
\DMod_\hol^{!,*}(  \psi_P)(\mCD_Q) \simeq 0.
\end{equation}
Let us execute the induction step using this claim. Note that the underlying simplicial set of the category $\Par'\cap \Par_{\ge Q}$ is weakly contractible, hence by \cite[Corollary 4.4.4.10]{HTT}
\[ \colim_{Q'\in \Par'\cap \Par_{\ge Q} } \mbK(G) \simeq \mbK(G) . \]
Hence
$$ \mCD_{Q} \simeq \colim_{Q'\in \Par'\cap \Par_{\ge Q} } \on{coFib}( \mbK(Q') \to \mbK(G)).$$
By induction hypothesis, the functor $\DMod_\hol^{!,*}(  \psi_P)$ sends $\on{coFib}( \mbK(Q') \to \mbK(G))$ to $0$ unless $Q'= Q$. Hence $\DMod_\hol^{!,*}(  \psi_P)(\mCD_Q)$ is isomorphic to\footnote{We use the following formal fact. Let $[1]$ be the poset $\{0,1\}$ and $I$ be an index category obtained by removing the final object from $[1]^r$ ($r\ge 1$). Let $C$ be any stable category. Suppose $F:I\to C$ is a functor such that $F(x)\simeq 0$ unless $x$ is the initial object $i_0$. Then $\colim F \simeq F(i_0)[r-1]$. This fact can be proven by induction on $r$.}
$$ \DMod_\hol^{!,*}(  \psi_P)(\on{coFib}( \mbK(Q) \to \mbK(G)) [ \on{rank}(G)-\on{rank}(L) +1].$$
Then the claim (\ref{eqn-claim-DQ}) implies 
\[ \DMod_\hol^{!,*}(  \psi_P)\circ \mbK(Q) \to \DMod_\hol^{!,*}(  \psi_P)\circ \mbK(G) \]
is an isomorphism as desired.

It remains to prove (\ref{eqn-claim-DQ}). Consider the maps
$$ \ol{\Bun}_{G,G} \xto{j_{G,\ge Q} } \ol{\Bun}_{G,\ge Q} \xto{j_{\ge Q}} \ol{\Bun}_{G}.$$
By Lemma \ref{lem-gluing-functors-colimite}, we have $\mCD_Q  \simeq j_{\ge Q,!} (\mCF)$, where 
$$\mCF := (j_{G,\ge Q})_* \circ r_!(k_{\Bun_G}).$$ 
Hence by the base-change isomorphism, $\DMod_\hol^{!,*}(  \psi_P)(\mCD_Q)$ is isomorphic to $\DMod_\hol^{!,*}(\psi_{P,\ge Q}) (\mCF) $, where
$$ \psi_{P,\ge Q}: (   H_{\MGPos}/Z_M \gets   (Y^P_{\rel}/Z_M)_{\ge Q}  \to   \ol{\Bun}_{G,\ge Q}  ) $$
and $ (Y^P_{\rel}/Z_M)_{\ge Q}$ is the open substack of $Y^P_{\rel}/Z_M$ containing those $Q'$-strata with $Q'\supset Q$. The following construction will be provided in $\S$ \ref{ssec-contraction}:

\begin{propconstr} \label{propconstr-contration}
The correspondence $\psi_{P,\ge Q}$ is isomorphic to the composition of 
$$\psi_{Q,\ge Q}: (  H_{ L,G\on{-pos}}/Z_L    \gets   Y^Q_{\rel}/Z_L   \to   \ol{\Bun}_{G,\ge Q}  )$$
with a certain correspondence from $H_{ L,G\on{-pos}}/Z_L $ to $H_{ M,G\on{-pos}}/Z_M $.
\end{propconstr}

Therefore we only need to show $\DMod_\hol^{!,*}(\psi_{Q,\ge Q}) (\mCF) \simeq 0$. We will prove the following stronger claim: for any $Q\in \Par$, we have
$$ \DMod_\hol^{!,*}(\psi_{Q,\ge Q})\circ (j_{G,\ge Q})_* \simeq 0,$$
where recall
\[j_{G,\ge Q}: \ol{\Bun}_{G,G}\to   \ol{\Bun}_{G,\ge Q},\, \psi_{Q,\ge Q}: (  H_{ L,G\on{-pos}}/Z_L    \gets   Y^Q_{\rel}/Z_L   \to   \ol{\Bun}_{G,\ge Q}  ). \]

To finish the proof, we need one more geometric input. In \cite{schieder2016geometric}, the author constructed a \emph{defect stratification} on the $P$-stratum $\ol{\Bun}_{G,P}$. Let $_\df\ol{\Bun}_{G, P}$ be the disjiont union of all the defect strata. It is known (see $\S$ \ref{sssec-equivariant-category}) that
$$ _\df\ol{\Bun}_{G, P} \simeq  \BunP \mt_{\BunM} ( H_\MGPos/Z_M ) \mt_{\BunM} \BunPm .$$ 
Consider the diagram
\begin{eqnarray*}
     H_{\MGPos}/Z_M & \xgets{\mfq^+_{P,\Vin}}  \, _\df\ol{\Bun}_{G,P} \xto{\mfp^+_{P,\Vin}}  &   \ol{\Bun}_{G,\ge P}, \\
     H_{\MGPos}/Z_M  &  \xgets{\mfq^-_{P,\Vin}}   \, Y_\rel^P/Z_M\xto{\mfp^-_{P,\Vin}}  &  \ol{\Bun}_{G,\ge P}  . \\
\end{eqnarray*}
In $\S$ \ref{ssec-interplay}, we will prove the following ``second-adjointness-style'' result: 
\begin{thm}  \label{thm-interplay} 
	The functor
	$$ \mfq^\mp_{P,\Vin,!}\circ \mfp^{\mp,*}_{P,\Vin}:\Dmod_\hol( \ol{\Bun}_{G,\ge P})  \to \Dmod_\hol(   H_{\MGPos}/Z_M)   $$
	is isomorphic to the restriction of the functor
	$$ \mfq_{P,\Vin,*}^\pm\circ \mfp_{P,\Vin}^{\pm,!}: \Dmod( \ol{\Bun}_{G,\ge P})  \to \Dmod(  H_{\MGPos}/Z_M).$$

\end{thm}

Now the $Q$-version of Theorem \ref{thm-interplay} says 
$$\DMod_\hol^{!,*}(\psi_{Q,\ge Q}) \simeq \mfq_{Q,\Vin,*}^+\circ \mfp_{Q,\Vin}^{+,!} .$$ Hence we have
$$ \DMod_\hol^{!,*}(\psi_{Q,\ge Q})\circ (j_{G,\ge Q})_*  \simeq \mfq_{Q,\Vin,*}^+\circ \mfp_{Q,\Vin}^{+,!}  \circ (j_{G,\ge Q})_*   .$$
Note that $_\df\ol{\Bun}_{G,Q}$ and $\ol{\Bun}_{G,G}$ have empty intersection (because $Q\ne G$). Hence $\mfp_{Q,\Vin}^{+,!}  \circ (j_{G,\ge Q})_*\simeq 0$. This finishes the proof. 

\qed[Theorem \ref{thm-main-variant}]

\begin{rem} 
In the case $G=\on{SL}_2$, one can use Theorem \ref{thm-interplay} to give a quicker proof of Goal \ref{prop-finish}. Namely, using the theorem, we only need to show $\mfq_{B,\Vin,*}^+\circ \mfp_{B,\Vin}^{+,!}\circ \mbK$ sends the arrow $B\to G$ to an isomorphism. Recall that $\mfp_{B,\Vin}^{+}$ factors through
$$  i_B:  \ol{\Bun}_{G,B} \to \ol{\Bun}_G.$$
Hence we only need to show $i_B^!\circ \mbK$ sends $B\to G$ to an isomorphism. However, this is obvious because the image of this arrow is the map (see Remark \ref{rem-compare-gluing-functor}):
$$
  i_B^* \circ j_{G,*}( r_!(k_{\BunG}) )[-1] \to i_B^!\circ j_{G,!}( r_!(k_{\BunG}) ),
$$
which is an isomorphism because $i_B$ and $j_G$ are complementary to each other.
\end{rem}

\section{Step 1}
We have three results to prove in this step: Proposition-Construction \ref{propconstr-factor-through-gen}, Proposition \ref{prop-FP-equivariant} and Proposition \ref{propconstr-nt-cartesian}. Each subsection corresponds to a result. Note that we have to do them in this order because the statement of the second result depends on the construction of the first.

The main idea of the proof of Proposition-Construction \ref{propconstr-factor-through-gen} is to construct a canonical map $G\backslash \,_0\!\Vin_{G,\le P}/G \to \mBB P\mt \mBB P^-$, which is achieved by studying the stablizer group $\wt{G}$ (reviewed in Appendix \ref{appendix-wtG}).

The main ingredient of the proof of Proposition \ref{prop-FP-equivariant} is the \emph{defect stratification} on $\ol{\Bun}_G$ defined in \cite{schieder2016geometric} which allows us to define certain full subcategory $\Dmod(\ol{\Bun}_{G,P})^{\mbU\mt \mbU^-}\subset \Dmod(\ol{\Bun}_{G,P})$ closely connected to the category $\mbI(G\mt G,P\mt G)$. Such connection is due to a comparison of the defect stratification of $\ol{\Bun}_G$ and that of the Drinfeld compactification $\wt{\Bun}_P$ defined in \cite{braverman2002geometric}.

The proof of Proposition \ref{propconstr-nt-cartesian} is purely formal.

\ssec{Proof of Proposition-Construction \ref{propconstr-factor-through-gen}}
\label{ssec-step-1-1}
Our goal is:
\begin{goal}
	Construct canonical maps
	$$\ol{\Del}^\enh_{\le P}:\ol{\Bun}_{G,\le P} \to \Bun_G^{P\hgen}\mt \BunG^{P^-\gen}$$
	that are functorial in $P$ such that when $P=G$ we have $\ol{\Del}^\enh_{\le G}=\ol{\Delta}$.
\end{goal}

	By definition, we have $\ol{\Bun}_{G,\le P}\simeq \VinBun_{G,\le P}/T$ and
	$$ \VinBun_{G,\le P} \simeq \bMap_\gen(X, G\backslash \Vin_{G,\le P}/G \supset G\backslash \,_0\!\Vin_{G,\le P}/G ),$$
	where the $T$-action on $\VinBun_{G,\le P}$ is induced from the $T$-action on $\Vin_{G,\le P}$. By Fact \ref{fact-wtg}(2), we have
	\begin{equation}\label{ean-proof-propconstr-factor-through-gen-1}
	  G\backslash \,_0\!\Vin_{G,\le P}/G \simeq \mBB \wtG_{\le P},
	 \end{equation}
	where $\wtG_{\le P}$ is a subgroup scheme of $G\mt G\mt T_{\ad,\le P}^+\to T_{\ad,\le P}^+$.

	By Lemma \ref{lem-wtG-leP-containing}, the map $\mBB \wtG_{\le P}  \to \mBB G\mt \mBB G$ factors as $ \mBB \wtG_{\le P} \to \mBB P\mt \mBB P^- \to \mBB G\mt \mBB G. $
	Also, the maps 
	\begin{equation}\label{ean-proof-propconstr-factor-through-gen-2} \mBB \wtG_{\le P} \to \mBB P\mt \mBB P^-
	\end{equation}
	 are functorial in $P$. Now we have the following commutative diagram of algebraic stacks:
	\begin{equation}\label{ean-proof-propconstr-factor-through-gen-3}
	\xyshort
	\xymatrix{
		G\backslash \Vin_{G,\le P}/G \ar[d] &
		G\backslash \,_0\!\Vin_{G,\le P}/G \ar[d]^-{ (\ref{ean-proof-propconstr-factor-through-gen-2})\circ (\ref{ean-proof-propconstr-factor-through-gen-1}) } \ar[l] \\
		\mBB G\mt \mBB G &
		\mBB P\mt \mBB P^-. \ar[l]
	}
	\end{equation}
	Taking $\bMap_\gen(X,-)$, we obtain maps
	$$\VinBun_{G,\le P} \to \Bun_G^{P\hgen}\mt \BunG^{P^-\gen}$$
	functorial in $P$. To finish the construction, we only need to show:

	\begin{lem} \label{lem-unique-lifting-group-action-VinleP-to-BPBPm}
	The map $G\backslash \,_0\!\Vin_{G,\le P}/G\to \mBB P\mt \mBB P^-$ constructed above can be uniquely lifted to a morphism
	$$ (T\act G\backslash \,_0\!\Vin_{G,\le P}/G)\to ( \pt \act  \mBB P\mt \mBB P^- ) $$
	fitting into the following commutative diagram
	$$
	\xyshort
	\xymatrix{
		(T\act G\backslash \Vin_{G,\le P}/G) \ar[d] &
		(T\act G\backslash \,_0\!\Vin_{G,\le P}/G) \ar[d] \ar[l] \\
		(\pt \act \mBB G\mt \mBB G) &
		(\pt \act \mBB P\mt \mBB P^-). \ar[l]
	}
	$$
	\end{lem}

	\proof The uniqueness follows from the fact that $ \mBB P\mt \mBB P^-\to \mBB G\mt \mBB G$ is schematic. It remains to prove the existence.

	The map $G\backslash \,_0\!\Vin_{G,\le P}/G\to \mBB P\mt \mBB P^-$ induces a $(G\mt G)$-equivariant map
	\begin{equation}\label{ean-proof-propconstr-factor-through-gen-4}
	 _0\!\Vin_{G,\le P} \to  G/P\mt G/P^-.
	 \end{equation}
	We only need to show the $T$-action on $_0\!\Vin_{G,\le P}$ preserves the fibers of this map.

	Recall that any closed point in $_0\!\Vin_{G,\le P}$ is of the form $g_1\cdot \mfs(s)\cdot g_2^{-1}$ where $g_1$ and $g_2$ are closed points of $G$, $s$ is a closed point of $T_{\ad,\le P}^+$ and $\mfs$ is the canonical section. Unwinding the definitions, the map (\ref{ean-proof-propconstr-factor-through-gen-4}) sends this point to $(g_1,g_2)$. Now consider the $T$-action on $\Vin_G$. It follows from definition that a closed point $t$ of $T$ sends the point $\mfs(s)$ to the point $\iota(t)\cdot \mfs(ts)$, where $\iota:T\inj G$ is the embedding. Since the $T$-action commutes with the $(G\mt G)$-action, the element $t$ sends $g_1\cdot \mfs(s)\cdot g_2^{-1}$ to $g_1\iota(t)\cdot \mfs(s)\cdot g_2^{-1}$. This makes the desired claim manifest.
	
	\qed[Lemma \ref{lem-unique-lifting-group-action-VinleP-to-BPBPm}]

	\qed[Proposition-Construction \ref{propconstr-factor-through-gen}]

\ssec{Proof of Proposition \ref{prop-FP-equivariant}} \label{ssec-step-1-2}
\begin{goal} 
	The object $\mCF_P$ is contained in the full subcategory
	$$  \mbI(G\mt G,P\mt G) \subset \Dmod(  \Bun_G^{P\hgen}\mt \BunG  ).$$
\end{goal}

\sssec{The $(\mbU,\mbU^-)$-equivariant categories}
	\label{sssec-equivariant-category}
	In order to prove Proposition \ref{prop-FP-equivariant}, we will introduce a subcategory
	$$ \Dmod(\ol{\Bun}_{G,P})^{\mbU\mt \mbU^-} \subset  \Dmod(\ol{\Bun}_{G,P}),$$
	whose definition is similar to 
	$$\mbI(G\mt G,P\mt P^-) \subset \Dmod( \Bun_G^{P\hgen}\mt \Bun_G^{P^-\hgen} ).$$

	To define this subcategory, we use the \emph{defect stratification} on $\ol{\Bun}_{G,P}$ defined\footnote{More precisly, \cite{schieder2016geometric} constructed the defect stratification on $\VinBun_{G,C_P}$. It follows from the construction that the $Z_M$-action on $\VinBun_{G,C_P}$ preserves the defect strata. Hence we obtain a stratification on 
	$$\ol{\Bun}_{G,P}\simeq \VinBun_{G,P}/T \simeq \VinBun_{G,C_P}/Z_M.$$} in \cite{schieder2016geometric}. Recall the disjoint union of its strata is given by\footnote{The corresponding $Z_M$-action on $\BunP \mt_{\BunM} H_{\MGPos} \mt_{\BunM} \BunPm$ is the one induced by the $Z_M$-action on $H_{\MGPos}$. Note that the map $H_\MGPos \to \BunM \mt \BunM$ is $Z_M$-equivariant for this action and the \emph{trivial} action on $\BunM\mt \BunM$.}
	$$ _\defect \ol{\Bun}_{G,P} \simeq \BunP \mt_{\BunM} ( H_\MGPos/Z_M ) \mt_{\BunM} \BunPm,$$
	or more abstractly
	$$ _\defect \ol{\Bun}_{G,P} \simeq \bMap_\gen( X, P\backslash \ol{M}/P^- \supset P\backslash M/P^- )/Z_M,$$
	where $\ol{M}$ is the closure of the locally closed embedding
	$$  M \simeq (P\mt P^-)/(P\mt_M P^-) \inj (G\mt G)/(P\mt_M P^-) \simeq \,_0\!\Vin_{G,C_P}.$$

	It is well-known that the map $\BunP \to \BunM$ is \emph{universally homological contractible}, or $\emph{UHC}$. In other words, for any lft prestack $Y\to \BunM$, the $!$-pullback functor $\Dmod(Y) \to \Dmod( Y\mt_{\BunM} \BunP )$ is fully faithful. In particular, the following $!$-pullback functor is fully faithful
	$$_\defect \ol{\Bun}_{G,P} \to  (H_{\MGPos}/Z_M) \mt_{\BunM} \BunPm .$$
	We denote its essential image by $\Dmod(_\defect \ol{\Bun}_{G,P} )^{\mbU_P}$. Similarly we define $\Dmod(_\defect \ol{\Bun}_{G,P} )^{\mbU_P^-}$ and $\Dmod(_\defect \ol{\Bun}_{G,P} )^{\mbU_P\mt \mbU_P^-}$.

	Since $\BunP\to \BunM$ is smooth, in the previous definition, we can also use $*$-pullbacks instead of the $!$-pullbacks. The resulting subcategories are the same.

	We define $\Dmod(\ol{\Bun}_{G,P})^{\mbU_P}$ to fit into the following pullback diagram
	$$
	\xyshort
	\xymatrix{
	\Dmod(\ol{\Bun}_{G,P})^{\mbU_P} \ar[r]^-\subset \ar[d] &
	\Dmod(\ol{\Bun}_{G,P}) \ar[d]^-{!\on{-pull}} \\
	 \Dmod( _\defect \ol{\Bun}_{G,P} )^{\mbU_P} \ar[r]^-\subset &
	  \Dmod( _\defect \ol{\Bun}_{G,P}).
	}
	$$
	Similarly we define $\Dmod(\ol{\Bun}_{G,P} )^{\mbU_P^-}$ and $\Dmod( \ol{\Bun}_{G,P} )^{\mbU_P\mt \mbU_P^-}$. We also define the version of these sub-categories for ind-holonomic D-modules.

	We will deduce Proposition \ref{prop-FP-equivariant} from the following three lemmas. The proof of the first one is completely similar to that in \cite[Appendix $\S$ G.1]{chen2021thesis}. We provide the proofs for the other two.

\begin{lem} \label{lem-equivariant-for-K}
	For any morphism $P\in \Par$, the object 
	$$i_P^*\circ  j_{G,*} \circ r_{!}(k_{\BunG}) \in \Dmod_\hol(\ol{\Bun}_{G,P})$$
	is contained in $\Dmod_\hol(\ol{\Bun}_{G,P})^{\mbU_P\mt \mbU_P^-}$.
\end{lem}

\begin{lem}  \label{lem-equivariant-for-!-extension}
	The $!$-pushforward functor
	$$ \Dmod_\hol (_\defect\ol{\Bun}_{G,P}) \to \Dmod_\hol (\ol{\Bun}_{G,P})$$
	preserves $(\mbU_P\mt \mbU_P^-)$-equivariant objects.
\end{lem}

\proof 
	It suffices to prove the similar version after replacing $\ol{\Bun}_{G,P}$ by its smooth cover $\VinBun_{G,C_P}$. By \cite[$\S$ 3.3.2]{schieder2016geometric}, the map $f:\, _\defect\VinBun_{G,C_P} \to \VinBun_{G,C_P}$ factors as
	$$ _\defect\VinBun_{G,C_P} \xto{j} \wt{\Bun}_P \mt_{\BunM} H_{\MGPos} \mt_{\BunM} \wt{\Bun}_{P^-} \xto{\ol{f}} \VinBun_{G,C_P}$$
	such that $j$ is a schematic open embedding and $\ol{f}$ is proper on each connected component. Recall that $\wt{\Bun}_P$ also has a defect stratification with
	$$_\defect\wt{\Bun}_P \simeq \BunP \mt_{\BunM} H_{\MGPos}.$$

	We define $\Dmod(_\defect\wt{\Bun}_P )^{\mbU_P}$ to be the full subcategory of $\Dmod(_\defect\wt{\Bun}_P )$ consisting of objects that are $!$-pullbacks from $\Dmod(H_{\MGPos})$. We define $\Dmod(\wt{\Bun}_P )^{\mbU_P}$ similarly as before. We also define 
	$$\Dmod(  \wt{\Bun}_P \mt_{\BunM} H_{\MGPos} \mt_{\BunM} \wt{\Bun}_{P^-} )^{  \mbU_P\mt \mbU_P^-}.$$

	We claim the functor $\ol{f}_!$ preserves $(\mbU_P\mt \mbU_P^-)$-equivariant objects. To prove the claim, we use the fact that $\ol{f}$ is compatible with the defect stratifications. In other words, we have
	\begin{eqnarray*}
	_\defect\VinBun_{G,C_P} \mt_{\VinBun_{G,C_P} } (\wt{\Bun}_P \mt_{\BunM} H_{\MGPos} \mt_{\BunM} \wt{\Bun}_{P^-}) \simeq \\
	\simeq (\BunP  \mt_{\BunM} H_{\MGPos}) \mt_{\BunM}   H_{\MGPos}   \mt_{\BunM} (H_{\MGPos} \mt_{\BunM} \BunPm),
	\end{eqnarray*}
	such that the porjection from the RHS to 
	$$_\defect\VinBun_{G,C_P} \simeq  \BunP  \mt_{\BunM} H_{\MGPos} \mt_{\BunM} \BunPm $$ is induced by the ``composition'' map
	$$H_{\MGPos} \mt_{\BunM}   H_{\MGPos}   \mt_{\BunM} H_{\MGPos}\to H_{\MGPos}.$$
	Then the claim follows from the base-change isomorphisms (which exist because $\ol{f}$ is proper on each connected component).

	It remains to show $j_!$ preserves $(\mbU_P\mt \mbU_P^-)$-equivariant objects. Using the base-change isomorphism, it suffices to show that the $!$-pushforward functor
	$$ \Dmod_\hol(\BunP) \to \Dmod_\hol(\wt{\Bun}_P ) $$
	preserves $\mbU_P$-equivariant object. However, this is well-known and was proved in $\S$ \ref{appendix-proof-prop-well-defined-iota_M_!}.

\qed[Lemma \ref{lem-equivariant-for-!-extension}]

\begin{lem} \label{lem-!-push-equivariant-to-IGP}
	The functor 
	$$\ol{\Delta}_{P,!}^{\enh,l}:\Dmod_\hol(\ol{\Bun}_{G,P}) \to \Dmod_\hol( \Bun_G^{P\hgen} \mt \BunG )$$
	sends objects in $\Dmod_\hol(\ol{\Bun}_{G,P})^{\mbU_P}$ to objects in $\mbI(G\mt G,P\mt G)$.
\end{lem}

\proof 
	Lemma \ref{lem-equivariant-for-!-extension} formally implies $\Dmod_\hol(\ol{\Bun}_{G,P})^{\mbU_P}$ is generated under colimits and extensions by the image of the $!$-pushforward funnctor
	$$\Dmod_\hol(_\df\ol{\Bun}_{G,P})^{\mbU_P}\to \Dmod_\hol(\ol{\Bun}_{G,P})^{\mbU_P}. $$
	Hence it suffices to show the $!$-pushforward along
	$$ _\df\ol{\Bun}_{G,P}\to \ol{\Bun}_{G,P} \to  \Bun_G^{P\hgen}\mt \BunG $$
	sends $\mbU_P$-equivariant objects into $\mbI(G\mt G,P\mt G)$. Unwinding the definitions, this map is isomorphic to
	$$ \BunP \mt_{\BunM} (H_\MGPos/Z_M) \mt_{\BunM} \BunPm \xto{a} \BunP \mt \BunG\xto{b} \Bun_G^{P\hgen}\mt \BunG.$$
	By the base-change isomorphism, $a_!$ preserves $\mbU_P$-equivariant objects. Then we are done because $b_!$ sends $\mbU_P$-equivariant objects into $\mbI(G\mt G,P\mt G)$ by Proposition \ref{prop-well-defined-Eis_enh}(1).

\qed[Lemma \ref{lem-!-push-equivariant-to-IGP}]

\sssec{Finish the proof}
	Recall 
	$$\mCF_P \simeq \ol{\Del}^{\enh,l}_{P,!} \circ i_{P}^!\circ \mbK(P).$$
	By definition, we have 
	$$i_{P}^!\circ \mbK(P)\simeq i_{P}^* \circ j_{G,*}\circ r_!(k_{\Bun_G}) [\on{rank}(M)-\on{rank}(G)].$$
	By Lemma \ref{lem-equivariant-for-K}, this is an $\mbU_P$-equivariant object. Then we are done by Lemma \ref{lem-!-push-equivariant-to-IGP}.

\qed[Proposition \ref{prop-FP-equivariant}]

\ssec{Proof of Proposition-Construction \ref{propconstr-nt-cartesian}} \label{ssec-step-1-3} 
\begin{goal} 
Construct a canonical natural transformation $\ol{\Delta}_!\circ \mbK\to \mbE\circ \mbD\mbL$ whose value at $P\in \Par$ is equivalent to the morphism (\ref{eqn-value-natural-transformation})
\end{goal}

Proposition \ref{prop-well-defined-Eis_enh} provides a functor
$$ \mbI(G,-): \Par \to \DGCat_\cont$$
that sends an arrow $P\to Q$ to the functor $\Eis_{P\to Q}^\enh$. Hence we also have a functor
\begin{equation} \label{eqn-IGG-?G}
 \mbI(G\mt G,-\mt G):\Par \to \DGCat_\cont 
\end{equation}
that sends an arrow $P\to Q$ to the functor $\Eis_{P\mt G\to Q\mt G}^\enh$. 

\begin{lem} \label{lem-double-the-group} The functor (\ref{eqn-IGG-?G}) is canonically isomorphic to the functor
$$ \Par\to \DGCat_\cont,\; P\mapsto \mbI(G,P)\ot_k\Dmod(\BunG).$$
\end{lem}

\proof By the proof of \cite[Corollary 2.3.4]{drinfeld2013some}, the functor $\Dmod(Y)\ot_k\Dmod(\BunG) \to \Dmod(Y\mt \BunG)$ is an equivalence for any lft prestack $Y$. Then the lemma follows from definitions.

\qed[Lemma \ref{lem-double-the-group}]

Let $\wt{\mbI} \to \Par$ be the presentable fibration\footnote{A presentable fibration is both a Cartesian fibration and a coCartesian fibration whose fibers are presentable $(\infty,1)$-categories. See \cite[Definition 5.5.3.2]{HTT}.} classifying the functor (\ref{eqn-IGG-?G}). Note that $\Par$ has a final object $G$, and the fiber of $\wt\mbI$ at this object is $\wt\mbI_G := \Dmod(\BunG\mt \BunG)$. Consider the trivial fibration $ \wt{\mbI}_G \mt \Par \to \Par$. It follows formally that we have an adjoint pair
$$ \Eis^\enh :\, \wt{\mbI} \, \adj \, \wt\mbI_G\mt \Par \,: \CT^\enh,$$
where $\Eis^\enh$ (resp. $\CT^\enh$) preserves co-Cartesian (resp. Cartesian) arrows and its fiber at $P\in \Par$ is $\Eis^\enh_{P\mt G\to G\mt G}$ (resp. $\CT_{G\mt G\gets P\mt G}^\enh$). Using Lemma \ref{lem-double-the-group}, the functor $\mbE\circ\mbD\mbL$ is isomorphic to 
$$ \Par \xto{(\Delta_!(k_{\BunG}),-)}   \wt\mbI_G\mt \Par \xto{ \CT^\enh  }  \wt{\mbI} \xto{ \Eis^\enh  }   \wt\mbI_G\mt \Par \xto{\onpr} \wt\mbI_G.$$
Denote the composition of the first two functors by $\mbS_{\CT}:\Par \to \wt{\mbI}$. Note that it is the unique Cartesian section whose value at $G\in \Par$ is $\Delta_!(k_{\BunG})\in \wt{\mbI}_G$.

We also have a functor
$$ \Par^\op \to \DGCat_\cont,\; P\mapsto \Dmod(\Bun_G^{P\hgen}\mt \BunG)$$
that sends an arrow to the corresponding $!$-pullback functor. Let $\mCD_\gen\to \Par$ be the corresponding \emph{Cartesian} fibration. By Proposition \ref{prop-well-defined-Eis_enh}(1), we have a fully faithful functor $\wt{\mbI}\to \mCD_\gen$ that preserves \emph{co-Cartesian} arrows (although $\mCD_\gen$ is not a co-Cartesian fibration).

On the other hand, consider the functor
$$ \Par \to \DGCat_\cont,\; P\mapsto \Dmod_\hol(\ol{\Bun}_{G,\le P})$$
that sends an arrow to the corresponding $!$-extension functor. Let $ \ol\mCD \to \Par$ be the presentable fibration classifying this functor. We have a fully faithful functor
$$ \ol\mCD \to \Dmod_\hol(\ol{\Bun}_{G}) \mt \Par$$
whose fiber at $P\in \Par$ is the corresponding $!$-extension functor. The graph of the functor $\mbK$:
$$ \Par \to \Dmod_\hol(\ol{\Bun}_{G}) \mt \Par,\; P\mapsto ( \mbK(P),P)$$
is contained in the above full subcategory $\ol{\mCD}$. Hence we obtain a section $\mbS_\mbK:\Par \to \ol\mCD$ to the projection $\ol\mCD\to \Par$.

By Proposition \ref{propconstr-factor-through-gen}, we also have functorial maps
$$\ol{\Del}^{\enh,l}_{\le P}:\ol{\Bun}_{G,\le P} \to \Bun_G^{P\hgen}\mt \BunG.$$
Hence there is a functor 
$$\ol{\mCD} \to \mCD_\gen$$ that preserves co-Cartesian arrows such that its fiber at $P\in \Par$ is the composition
$$ \Dmod_\hol( \ol{\Bun}_{G,\le P}  ) \xto{ \ol{\Del}^{\enh,l}_{\le P,!} } \Dmod_\hol(\Bun_G^{P\hgen}\mt \BunG) \to \Dmod(\Bun_G^{P\hgen}\mt \BunG).$$
By construction, the composition
$$ \Par \xto{\mbS_\mbK}  \ol{\mCD} \to \mCD_\gen $$
sends $P$ to $\mCF_P$, viewed as an object in $\mCD_\gen$ over $P\in \Par$. Hence by Proposition \ref{prop-FP-equivariant}, this functor factors through the full subcategory $\wt{\mbI} \subset \mCD_\gen$. Let $\mbS_\mbK':\Par\to \wt{\mbI}$ be the corresponding functor. By constuction, $\ol{\Delta}_!\circ \mbK $ is isomorphic to the composition 
$$ \Par \xto{\mbS_\mbK' }   \wt{\mbI} \xto{ \Eis^\enh  }  \wt\mbI_G\mt \Par \xto{\onpr} \wt\mbI_G.$$

In summary, we have obtained two sections $\mbS_{\CT}$ and $\mbS_\mbK'$ to the Cartesian fibration $\wt{\mbI}\to \Par$ such that $\ol{\Delta}_!\circ \mbK$ and $\mbE\circ \mbD\mbL$ are obtained respectively by composing them with
$$  \wt{\mbI} \xto{ \Eis^\enh  }   \wt\mbI_G\mt \Par \xto{\onpr} \wt\mbI_G. $$

Now the identification $\mCF_G'=\mCF_G\simeq \Delta_!(k_{\BunG})$ provides an isomorphism $\mbS_\mbK'(G)\simeq \mbS_{\CT}(G)$. Since $G\in \Par$ is the final object and since $\mbS_{\CT}$ is a Cartesian section, we obtain a natural transformation $\mbS_\mbK'\to \mbS_{\CT}$ whose value at $P\in \Par$ is the unique arrow $\mbS_\mbK'(P)\to \mbS_{\CT}(P)$ fitting into the following commutative diagram
$$
\xyshort
\xymatrix{
	\mbS_\mbK'(P) \ar[r] \ar[d] &
	 \mbS_{\CT}(P) \ar[d] \\
	 \mbS_\mbK'(G) \ar[r]^-\simeq &  \mbS_{\CT}(G).
}
$$ 
By construction, when viewed as a morphism in $\wt{\mbI}_P \simeq \mbI(G\mt G,P\mt G)$, the arrow $\mbS_\mbK'(P)\to \mbS_{\CT}(P)$ is equivalent to (\ref{eqn-mCF-to-CT-Del}). Now the desired natural transformation $\ol{\Delta}_!\circ \mbK\to \mbE\circ \mbD\mbL$ is given by composing the above natural transformation $\mbS_\mbK'\to \mbS_{\CT}$ with the functor
$$  \wt{\mbI} \xto{ \Eis^\enh  }    \wt\mbI_G\mt \Par \xto{\onpr}\wt\mbI_G. $$

\qed[Proposition-Construction \ref{propconstr-nt-cartesian}]

\section{Step 2} \label{section-2nd-adjointness}

The goal of this section is to prove Theorem \ref{thm-generic-second-adjointness}. We first show the functors in the statement is well-defined in \S \ref{ssec-step-2-1}. Then in \S \ref{ssec-drinfeld} we review a general framework of Drinfeld in \cite{drinfeld2013algebraic} for proving adjointness. Finally, we verify Theorem \ref{thm-generic-second-adjointness} fits into the Drinfeld's framework in \S \ref{ssec-step-2-2}.

\begin{rem} Drinfeld's framework is highly abstract and uses $(\infty,2)$-categories. The readers might wonder why we use it rather than mimicking the \emph{second} proof in \cite{drinfeld2016geometric} or just applying the Braden's theorem. In this remark, we justify our choice.

First, Theorem \ref{thm-generic-second-adjointness} as well as the main theorem of \cite{drinfeld2016geometric} are not corollaries of the Braden's theorem because $\Bun_P$ is \emph{not} the attractor locus of any $\mBG_m$-action on $\Bun_G^{P\hgen}$ (resp. $\Bun_G$). The action on $\Bun_G$ induced by the $\mBG_m$-action on $G$ (whose attractor is $P$) is in fact trivial because the $\mBG_m$-action on $G$ is an inner-automorphism.

Then we explain why we do not mimic the second proof in \cite{drinfeld2016geometric}. In that proof, the authors use smooth descents to reduce to a problem about certain smooth atlases of $\BunG, \Bun_{P^\pm},\BunM$ where Braden's theorem can be applied. The problem for this method is $\Bun_G^{P\hgen}$ is not an algebraic stack and it does not have a smooth cover by schemes.

Let us also comment on the \emph{first} proof in \cite{drinfeld2016geometric}. As explained in \cite{drinfeld2013algebraic}, this proof, as well as Braden's proof of his theorem, in essence are examples of Drinfeld's framework. The general abstract framework has the advantage that one does not need to do the tedious and unmotivated step of verifying the adjunctions (like in \cite[Section 3]{drinfeld2016geometric}).
\end{rem}

\ssec{Well-definedness}
\label{ssec-step-2-1}

\begin{lem} \label{lem-quasi-compactness} We have:

\begin{itemize}
  \item[(1)]
The correspondences
\begin{eqnarray*}
  \alpha_{P,\lambda}^{+,\gen} &:& ( \Bun_{M,\lambda} \gets \Bun_{P,\lambda} \to \BunG^{P\hgen} ), \\
   \alpha_{P,\lambda}^{-,\gen} &:& ( \BunG^{P\hgen} \gets  \Bun_{P^-,\lambda}^{M\hgen} \to \Bun_{M,\lambda} )
\end{eqnarray*}
are morphisms in $\bCorr(\PreStk_\lft)_{\on{QCAD},\all}^{\open,2\on{-op}}$. In fact, the first leftward map is \emph{safe} and the second leftward map is quasi-compact and schematic. In particular, the functors in the statment of Theorem \ref{thm-generic-second-adjointness} are well-defined.

\item[(2)] There is a $2$-morphism $ \alpha_{P,\lambda}^{+,\gen} \circ \alpha_{P,\lambda}^{-,\gen} \to \on{Id}_{\Bun_{M,\lambda}}$ given by the map
$$ \Bun_{M,\lambda} \to \Bun_{P,\lambda}\mt_{\BunG^{P\hgen}} \Bun_{P^-,\lambda}^{M\hgen}.$$
In other words, this map is a schematic open embedding\footnote{In fact, $\Bun_{P}\mt_{\BunG^{P\hgen}} \Bun_{P^-}^{M\hgen}$ is the open Zastava space in the literature.}.

\end{itemize}
\end{lem}

The map $\Bun_{P,\lambda} \to \Bun_{M,\lambda}$ is safe by \cite[Footnote 2]{drinfeld2016geometric}. To prove the claim for the map $\BunG^{P\hgen} \gets  \Bun_{P^-,\lambda}^{M\hgen}$, consider the following commutative diagram
$$
\xyshort
\xymatrix{
	\Bun_{P^-}^{M\hgen} \ar[r] \ar[d] &
	 \BunG^{P\hgen} \ar[d] \\
	 \BunPm \ar[r] & \BunG.
}
$$
We claim it induces a schematic open embedding
$$ \Bun_{P^-}^{M\hgen}\to  \BunG^{P\hgen}\mt_{ \BunG} \BunPm.$$
Indeed, the RHS is isomorphic to $\bMap_\gen(X,  \mBB P^-\gets \mBB P\mt_{\mBB G} \mBB P^-)$ and the above map is isomorphic to the map
$$ \bMap_\gen(X,\mBB P^-\gets \mBB M ) \to  \bMap_\gen(X,  \mBB P^-\gets \mBB P\mt_{\mBB G} \mBB P^-)$$
induced by the map $\mBB M \to  \mBB P\mt_{\mBB G} \mBB P^-$. Then the claim follows from the fact that $\mBB M \to  \mBB P\mt_{\mBB G} \mBB P^-$ is a schematic open embedding.

Now (1) follows from the above claim and the well-known fact that $\Bun_{P^-,\lambda}\to \BunG$ is quasi-compact and schematic.

To prove (2), we only need to show
$$ \BunM \to \BunP \mt_{ \Bun_G^{P\hgen} } \Bun_{P^-}^{M\hgen} $$
is a schematic open embedding. As before, this follows from the fact that it is isomorphic to
$$ \bMap_\gen(X,\mBB M\gets \mBB M ) \to  \bMap_\gen(X,  \mBB P\mt_{\mBB G} \mBB P^- \gets \mBB M)$$
and the fact that $\mBB M \to  \mBB P\mt_{\mBB G} \mBB P^-$ is a schematic open embedding.

\qed[Lemma \ref{lem-quasi-compactness}]

\ssec{Recollections: Drinfeld's framework}
\label{ssec-drinfeld}
In \cite[Appendix C]{drinfeld2013algebraic}, Drinfeld set up a general framework to prove results like Theorem \ref{thm-generic-second-adjointness}. We review this framework in this subsection. In fact, we slightly generalize it to the case of lft prestacks.

\begin{defn}  \label{defn-PmbA1}
	We equip the category $\affSch_\ft$ with the Cartesian symmetric monoidal structure. Recall the notion of \emph{enriched categories}. Following \loccit, we define a category $\mbP_{\mBA^1}$ enriched in $\affSch_\ft$ as follows:

	\begin{itemize}
	\item It has two objects: the ``big'' one $\mbb$ and the ``small'' one $\mbs$.

	\item The mapping scheme $\bHom_{\mbP_{\mBA^1}}( \mbb,\mbb)$ is defined to be $\mBA^1$. The other three mapping shemes are defined to be $\pt$, viewed as the zero point in $\mBA^1$. The composition laws are all induced by the semi-group structure on $\mBA^1$.
	\end{itemize}

	The unique morphism $\mbs\to \mbb$ is denoted by $\alpha^+$ and the unique morphism $\mbb\to \mbs$ is denoted by $\alpha^-$.
\end{defn}

\begin{defn}  
	Let $\AlgStk_\QCAD$ be the $(2,1)$-category of QCAD algebraic stacks equipped with the Cartesian symmetric monoidal structure. We define a category\footnote{It was denoted by $\mbP_{\mBA^1}/\mBB\mBG_m$ in \cite{drinfeld2013algebraic}.} $\Dri$ enriched in $\AlgStk_\QCAD$ by replacing $\mBA^1$ in Definition \ref{defn-PmbA1} by the quotient stack $\mBA^1/\mBG_m$, and replacing the zero map $\pt\to \mBA^1$ by the map $\mBB \mBG_m \to \mBA^1/\mBG_m$ obtained by taking quotients.

	Note that there is an obvious functor $\mbP_{\mBA^1}\to \Dri$. We use the same symbols $\alpha^+$ and $\alpha^-$ to denote the corresponding morphisms in $\Dri$.
\end{defn}

\begin{defn} 
	Let $\mCO$ be a monoidal $(\infty,1)$-category, $\mCA$ be a category enriched in $\mCO$ and $\mCC$ be a module $(\infty,2)$-category of $\mCO$. As explained in \cite[$\S$ C.13.1]{drinfeld2013algebraic}, there is a notion of \emph{weakly $\mCO$-enriched (unital) right-lax functors}\footnote{It was called just by \emph{lax functors} in \loccit.} from $\mCA$ to $\mCC$. We will review its explicit meaning later in our particular examples. For now, let us give the formal definition.

	 We assume $\mCO$ is small. Consider the $(\infty,1)$-category $\Funct(\mCO^\op,(\infty,1)\on{-Cat})$ equipped with the Day convolution monoidal structure (see \cite[$\S$ 2.2.6]{HA}). Then $\mCC$ has a $\Funct(\mCO^\op,(\infty,1)\on{-Cat})$-enriched structure such that for any $x,y\in \mCC$, the object
	$$ \bHom_{\mCC}(x,y) \in  \Funct(\mCO^\op,(\infty,1)\on{-Cat}) $$
	is the functor $o\mapsto \bMap_\mCC(o\ot x,y)$.

	On the other hand, there is a \emph{right-lax} monoidal structure on the Yoneda functor
	$$  \mCO \to \Funct(\mCO^\op,(\infty,1)\on{-Cat}).$$
	Then a weakly $\mCO$-enriched functor (resp. right-lax functor) $F:\mCA\laxto \mCC$ is defined to be a functor (resp. right-lax functor) $F$ that intertwines the enrichment via the above right-lax monoidal functor.
\end{defn}

\begin{notn} 
	Consider the $(3,2)$-category $\bCorr(  \PreStk_\lft )_{\on{QCAD},\all}^{\open ,2\on{-op}}$. We equip it with the obvious $\AlgStk_\QCAD$-action.

	A \emph{Drinfeld pre-input} is a weakly $\AlgStk_\QCAD$-enriched	right-lax functor $F: \mbP_{\mBA^1} \dashrightarrow \bCorr$ such that it is strict at the composition $\alpha^+ \circ \alpha^-$, i.e., the $2$-morphism $F(\alpha^+) \circ F(\alpha^-) \to F( \alpha^+\circ \alpha^- )$ is invertible.

	A \emph{Drinfeld input} is a weakly $\AlgStk_\QCAD$-enriched	right-lax functor $F^\sharp:\Dri \dashrightarrow \bCorr$ such that the composition $\mbP_{\mBA^1}\to \Dri \dashrightarrow \bCorr$ is a Drinfeld pre-input.
\end{notn}

\begin{rem} \label{rem-translation-drinfeld-preinput}
	Unwinding the definitions, a Drinfeld pre-input provides
	
	\begin{itemize}
	\item Two lft prestacks $Z:=F(\mbb)$ and $Z^0:=F(\mbs)$;
	
	\item Two correspondences
	$$F(\alpha^+) :(Z \xgets{p^+} Z^+  \xto{q^+} Z^0)\;\on{ and }\; F(\alpha^-):(Z^0 \xgets{q^-} Z^-  \xto{p^-} Z) $$
	whose left arms are QCAD maps;
	
	\item An $\mBA^1$-family of correspondences:
	$$
	\xyshort
	\xymatrix{
	 Z & \wt{Z} \ar[r] \ar[l] \ar[d] & Z;\\
	 & \mBA^1  
	} $$
	given by $\Hom(\mbb,\mbb) \mt F(\mbb) \to F(\mbb)$;
	
	\item Isomorphisms 
	$$Z^+ \mt_{Z^0} Z^- \simeq \wt{Z}\mt_{\mBA^1} 0 \;\on{ and }\; Z\simeq \wt{Z}\mt_{\mBA^1} 1$$
	defined over $Z\mt Z$, given respectively by the invertible $2$-morphism $F(\alpha^+) \circ F(\alpha^-) \to F( \alpha^+\circ \alpha^- )$ and $\on{Id}_{F(\mbb)} \simeq F(\on{Id}_{\mbb})$
	
	\item An open embedding 
	$$j: Z^0 \to Z^-\mt_{Z} Z^+$$
	defined over $Z^0 \mt Z^0$, given by the lax composition law for $\mbs\gets\mbb\gets\mbs$;
	
	\item Open embeddings
	$$ Z^+ \mt \mBA^1 \to \wt{Z} \mt_Z Z^+ \;\on{ and }\; Z^- \mt \mBA^1 \to Z^- \mt_Z \wt{Z},  $$
	defined respectively over $Z\mt Z^0\mt \mBA^1$ and $Z^0\mt Z\mt \mBA^1$, given respectively by the lax composition laws for $\mbb \gets \mbb \gets \mbs$ and $\mbs \gets \mbb \gets \mbb$;
	
	\item An open embedding\footnote{The map $\mBA^2\to \mBA^1$ in the formula is the multiplication map.}
	$$  \wt{Z}\mt_{\mBA^1} \mBA^2 \to \wt{Z}\mt_Z \wt{Z} $$
	defined over $Z \mt Z \mt \mBA^2$, given by the lax composition law for $\mbb\gets \mbb \gets \mbb$.
	
	\item Some higher compatibilities.
	\end{itemize}
\end{rem}

\begin{exam} \label{exam-drinfeld-pre-input-Gm}
	For any finite type scheme $Z$ equipped with a $\mBG_m$-action, \cite{drinfeld2013algebraic} constructed a Drinfeld pre-input such that $Z^+$, $Z^-$ and $Z^0$ are respectively the attractor, repeller and fixed loci of $Z$. Also, $\wt{Z}$ is the so-called Drinfeld-Gaitsgory interpolation, which is an $\mBA^1$-degeneration from $Z$ to $Z^\att\mt_{Z^\fix} Z^\rep$. Moreover, this construction is functorial in $Z$ and compatible with Cartesian products.

	When $Z$ is affine, the corresponding right-lax functor $\mbP_{\mBA^1} \dashedrightarrow \bCorr $ is strict. In particular, we obtain a functor $\mbP_{\mBA^1} \to \Corr$.

	It was also shown in \loccit\,that there is a Drinfeld input with $F^\sharp(\mbb)=Z/\mBG_m$ and $F^\sharp(\mbs)=Z^\fix/\mBG_m$
\end{exam}

\sssec{Drinfeld's theorem on adjunctions}
	Let $F^\sharp:\Dri\dashrightarrow \bCorr(\PreStk_\lft)_{\on{QCAD},\all}^{\open,2\on{-op}}$ be a Drinfeld input and $F$ be the corresponding Drinfeld pre-input. We use the notations in Remark \ref{rem-translation-drinfeld-preinput}. Consider the composition
	$$\mbP_{\mBA^1} \xlaxto{F} \bCorr(\PreStk_\lft)_{\on{QCAD},\all}^{\open,2\on{-op}} \xto{\DMod^{\bt,!}} \bDGCat_\cont.$$
	By construction, it sends $\alpha^+$ and $\alpha^-$ respectively to the functors
	$$ \DMod^{\bt,!}\circ F(\alpha^+) \simeq p^{+}_\bt \circ q^{+,!}, \;  \DMod^{\bt,!}\circ F(\alpha^-) \simeq q^-_\bt \circ p^{-,!}.$$
	The $2$-morphism: $F(\alpha^-) \circ F(\alpha^+) \to F(\alpha^-\circ \alpha^+) = F(\on{Id}_\mbs)$ gives a natural transformation\footnote{Explicitly, the LHS is the $!$-pull-$\bt$-push along $Z^0 \gets Z^-\mt_Z Z^+\to Z^0$, while the RHS is that along $Z^0 \gets Z^0\to Z^0$. The desired natural transformation is induced by the adjoint pair $(j^!,j_\bt)$ for the open embedding $j:Z^0 \to Z^-\mt_Z Z^+$.}
	\begin{equation} \label{eqn-dri-counit}
 	q^-_\bt \circ p^{-,!} \circ  p^{+}_\bt \circ q^{+,!} \to \Id_{\Dmod(Z^0)}.
	\end{equation}

	The following result was proved in \cite[Appendix C]{drinfeld2013algebraic}.

	\begin{thm}\label{thm-drinfeld-framework}
 	(Drinfeld) In the above setting, there is an adjoint pair
	$$ q^-_\bt \circ p^{-,!}: \Dmod( Z) \adj  \Dmod(Z^0) : p^{+}_\bt \circ q^{+,!}$$
	with the counit adjunction natrual transformation given by (\ref{eqn-dri-counit}).
	\end{thm}

	\begin{rem} The unit adjunction is given by a specialization construction along $\wt{Z}\to \mBA^1$. We do not need it in this paper.
	\end{rem}

	\begin{rem} More precisly, \loccit\,focused on the problem of reproving the Braden's theorem (see \cite{braden2003hyperbolic}) using the Drinfeld input in Example \ref{exam-drinfeld-pre-input-Gm}. However, the proof there works for any Drinfeld input.
	\end{rem}

\ssec{Proof of Theorem \ref{thm-generic-second-adjointness}}
\label{ssec-step-2-2}

Throughout this subsection, we fix a co-character $\gamma:\mBG_m\to Z_M$ dominant and regular with respect to $P$. Note that the homomorphism $\mBG_m\to Z_M\to Z_M/Z_G$ can be uniquely extended to a homomorphism between semi-groups $\ol{\gamma}:\mBA^1\to T_{\ad,\ge C_P}^+$. Via $\gamma$, the adjoint action $T_\ad \act G$ induces an action $\mBG_m\act G$.

We first deduce the theorem from the following result:

\begin{propconstr} \label{propconstr-drinfeld-input}
There exists a canonical Drinfeld input 
$$F^\sharp:\Dri\dashrightarrow \bCorr(\PreStk_\lft)_{\on{QCAD},\all}^{\open,2\on{-op}}$$
such that\footnote{We also require that the $2$-morphism $F^\sharp(\alpha^+)\circ F^\sharp(\alpha^-)\to F^\sharp(\on{Id}_\mbs)$ is given by the obvious open embedding.
} it sends $\alpha^+$ and $\alpha^-$ respectively to 
\begin{eqnarray*}
  \Bun_G^{P^-\hgen}/\mBG_m \gets \Bun_{P,\lambda}^{M\hgen}/\mBG_m \to \Bun_{M,\lambda}/\mBG_m, \\
\Bun_{M,\lambda}/\mBG_m \gets \Bun_{P^-,\lambda}/\mBG_m \to  \Bun_G^{P^-\hgen}/\mBG_m.
\end{eqnarray*}
\end{propconstr}

\sssec{Deduction Theorem \ref{thm-generic-second-adjointness}}

We will use the mirror version of Proposition-Construction \ref{propconstr-drinfeld-input} by exchanging $P$ and $P^-$. Using Theorem \ref{thm-drinfeld-framework}, we obtain the version of Theorem \ref{thm-generic-second-adjointness} after replacing the relevant stacks by their $\mBG_m$-quotients. The same proof of \cite[Theorem 3.4.3]{drinfeld2014theorem} implies the following adjoint pair
$$  \DMod^{\bt,!}(\alpha_{P,\lambda}^{+,\gen}):  \Dmod(\BunG^{P\hgen})^{\mBG_m\on{-mon}} \adj \Dmod(\BunM) :\DMod^{\bt,!}(\alpha_{P,\lambda}^{-,\gen}), $$
where 
$$ \Dmod(\BunG^{P\hgen})^{\mBG_m\on{-mon}} \subset  \Dmod(\BunG^{P\hgen})$$ 
is the full subcategory generated by the essential image of the $!$-pullback functor
$$  \Dmod(\BunG^{P\hgen}/\mBG_m) \to  \Dmod(\BunG^{P\hgen}). $$
Then we are done because the $\mBG_m$-action on $\BunG^{P\hgen}$ can be trivialized.

\qed[Theorem \ref{thm-generic-second-adjointness}]

It remains to construct the Drinfeld input in Proposition-Construction \ref{propconstr-drinfeld-input}.

\begin{notn}
	Let $\Grp^{\on{aff}}_\ft$ be the category of group schemes $H\to S$ with $H$ and $S$ being finite type affine schemes. Consider its arrow category $\Arr(\Grp^{\on{aff}}_\ft)$. We equip the category
	$$\Corr(\Arr(\Grp^{\on{aff}}_\ft))_{\all,\all}$$
	with the obvious $\affSch_\ft$-action.
\end{notn}

\begin{constr} 
	Via the co-character $\gamma$, the adjoint actions $Z_M\act G$ and $Z_M\act P$ induces actions $\mBG_m\act G$ and $\mBG_m\act P$. The corresponding attractor, repeller and fixed loci are:
	$$ G^{\att,\gamma} = P,\;G^{\rep,\gamma} = P^-,\;G^{\fix,\gamma} = M,\;(P^-)^{\att,\gamma} = M,\;(P^-)^{\rep,\gamma} = P^-,\;(P^-)^{\fix,\gamma} = M.$$
	Using Example \ref{exam-drinfeld-pre-input-Gm}, we obtain a weakly $\affSch_\ft$-enriched functor 
	$$\Theta_{P^-\to G}:\mbP_{\mBA^1} \to \Corr(\Arr(\Grp^{\on{aff}}_\ft))_{\all,\all}$$
	sending $\alpha^+$ and $\alpha^-$ respectively to
	$$  (P^-\to G) \gets (M\to P) \to (M\to M),\; (M\to M)  \gets (P^-\to P^-) \to (P^-\to G).$$
\end{constr}

\begin{rem} 
	By construction, $\Hom(\mbb,\mbb) \mt \Theta_{P^-\to G}(\mbb) \to \Theta_{P^-\to G}(\mbb)$ corresponds to the following diagram
	$$
	\xyshort
	\xymatrix{
	 (P^-\to G) & (\wt{P^-}^{\gamma}\to \wtG^\gamma) \ar[r] \ar[l] \ar[d] & (P^-\to G);\\
	 & \mBA^1,  
	} $$
	where $\wtG^\gamma$ (resp. $\wt{P^-}^{\gamma}$) is the Drinfeld-Gaitsgory interpolation for the action $\mBG_m\act G$ (resp. $\mBG_m\act P$). Note that we have
	 $$\wt{G}^\gamma \simeq \wtG\mt_{\Tadp,\ol{\gamma}} \mBA^1,\;  \wt{P^-}^\gamma \simeq P^-\mt_G \wt{G}^\gamma. $$ 
\end{rem}

\begin{constr} 
	Consider the functor 
	$$\mBB:\Grp^{\on{aff}}_\ft \to \AlgStk_\lft,\;(H\to S)\mapsto \mBB H,$$
	where $\mBB H:=S/H$ is the quotient stack. Similarly we have a functor $\Arr(\Grp^{\on{aff}}_\ft)\to \Arr(\AlgStk_\lft)$. This functor does not commute with fiber products, hence we only have a \emph{right-lax} functor
	$$ \Corr(\Arr(\Grp^{\on{aff}}_\ft))_{\all,\all}\laxto \bCorr(\Arr(\AlgStk_\lft))_{\all,\all}^{\all,2\on{-op}}.$$
	This right-lax functor has an $\affSch_\ft$-linear structure. Hence by composing with $\Theta_{P^-\to G}$, we obtain a weakly $\affSch_\ft$-enriched right-lax functor 
	$$ \Theta_{\mBB P^- \to \mBB G}: \mbP_{\mBA^1} \laxto \bCorr(\Arr(\AlgStk_\lft))_{\all,\all}^{\all,2\on{-op}}. $$
\end{constr}

\begin{defn} A morphism $(Y_1\to Y_2)\to (Y_1'\to Y_2')$ in $\Arr(\AlgStk_\lft)$ is called an \emph{open embedding} if both $Y_1\to Y_1'$ and $Y_2\to Y_2'$ are schematic open embeddings.
\end{defn}

\begin{lem} \label{lem-opennes-drinfeld-input}
	The right-lax functor $\Theta_{\mBB P^- \to \mBB G}$ factors through $\bCorr(\Arr(\AlgStk_\lft))_{\all,\all}^{\open,2\on{-op}}$ and is strict at the composition $\alpha^+\circ \alpha^-$.
\end{lem}

\proof 
	Consider the two forgetful functors $\Arr(\AlgStk_\lft)\to \AlgStk_\lft,\;(Y_1\to Y_2)\mapsto Y_i$. We only need to prove the similar claims after applying these forgetful functors. Those claims for the first forgetful functor are obvious (because $\wt{P^-}^\gamma \simeq P^-\mt \mBA^1$). It remains to prove those for the second forgetful functor.

	To prove the claim on strictness, we only need to check $\mBB (P\mt_M P^-) \to \mBB P\mt_{\mBB M} \mBB P^- $ is an isomorphism. The target is just $P\backslash M/P^-$. Hence this isomorphism follows from the fact $M\simeq (P\mt P^-)/(P\mt_M P^-)$.
	
	To prove the claim on openness, we only need to check that the following four maps are schematic open embeddings:
	\begin{eqnarray*}
	  \mBB (P^- \mt_G  \wtG^\gamma) \to \mBB P^- \mt_{\mBB G}  \mBB \wtG^\gamma, &  \mBB (\wtG^\gamma\mt_G P) \to \mBB \wtG^\gamma\mt_{\mBB G}  \mBB P, \\
	   \mBB (\wtG^\gamma\mt_G  \wtG^\gamma) \to \mBB \wtG^\gamma\mt_{\mBB G}  \mBB \wtG^\gamma &  \mBB (P^-\mt_G P)\to  \mBB P^-\mt_{\mBB G} \mBB P 
	 \end{eqnarray*}
   We first prove the claim for the last one. The target is just $P^-\backslash G/P$, it contains $P^-\backslash P^-P/P$ as an open substack, which is isomorphic to $\mBB (P^-\mt_G P)$ because $P^-\mt P$ acts transitively on $P^-P$ and the stablizer group at the unit element is $P^-\mt_G P$. The claims for the first two maps follows from Corollary \ref{cor-open-BPmwtG}. The proof for the third one is similar. Namely, consider the action
		$$(G\mt G\mt G) \act \,_0\!\Vin_G^\gamma \mt \,_0\!\Vin_G^\gamma,\; (g_1,g_2,g_3)\cdot (x_1,x_2)\mapsto ( g_1 x_1 g_2^{-1},g_2 x_2 g_3^{-1} ).$$
	Its stablizer for the canonical section is the group scheme $\wtG^\gamma\mt_G \wtG^\gamma$. We only need to prove the similar version of Lemma \ref{lem-Pm-s-G-open-orbit}, i.e., to show
	$$  ( G\mt G\mt G \mt \mBA^2  )/(\wtG^\gamma\mt_G \wtG^\gamma) \to  \,_0\!\Vin_G^\gamma \mt \,_0\!\Vin_G^\gamma $$
	is an open embedding. As before, we only need to show the LHS is smooth. Now the functor $\Theta_{P^-\to G}$ provides an isomorphism $\wtG^\gamma\mt_G \wtG^\gamma \simeq \wtG^\gamma \mt_{ \mBA^1} \mBA^2$
	covering the map
	$$\onpr_{13}\mt \on{Id}_{\mBA^2}: (G\mt G\mt G)\mt \mBA^2 \to (G\mt G)\mt \mBA^2. $$
	Hence we have a map
	$$   ( G\mt G\mt G \mt \mBA^2  )/(\wtG^\gamma\mt_G \wtG^\gamma) \to (G\mt G\mt \mBA^2)/(\wtG^\gamma \mt_{ \mBA^1} \mBA^2) \simeq \,_0\!\Vin_G^\gamma \mt_{ \mBA^1} \mBA^2. $$
	Then we are done because it is a smooth map to a smooth scheme.

\qed[Lemma \ref{lem-opennes-drinfeld-input}]

\begin{constr}
	Consider the functor
	$$  \bMap_\gen(X,-): \Arr(\AlgStk_\lft) \to \PreStk_\lft,\; (Y_1\to Y_2)\mapsto \bMap_\gen(X,Y_1\gets Y_2). $$
	It is easy to see that it sends open embeddings to schematic open embeddings. Hence we obtain a functor 
	$$  \bCorr(\Arr(\AlgStk_\lft))_{\all,\all}^{\open,2\on{-op}} \to \bCorr(\PreStk_\lft)_{\all,\all}^{\open,2\on{-op}}.$$
	This functor has an $\affSch_\ft$-linear structure\footnote{This is because for affine schemes $Y$, we have $\bMap(X,Y) \simeq Y$.}. Hence by composing with $\Theta_{\mBB P^- \to \mBB G}$, we obtain a weakly $\affSch_\ft$-enriched right-lax functor 
	$$  \Theta: \mbP_{\mBA^1} \laxto \bCorr(\PreStk_\lft)_{\all,\all}^{\open,2\on{-op}}$$
	that is strict at the composition $\alpha^+\circ \alpha^-$.
\end{constr}

\begin{rem} Explicitly, we have:
\begin{itemize}
	\item The right-lax functor $\Theta$ sends $\alpha^+$ and $\alpha^-$ respectively to 
	$$\Bun_G^{P^-\hgen} \gets \Bun_{P}^{M\hgen} \to \Bun_{M},\; 
	\Bun_{M} \gets \Bun_{P^-} \to  \Bun_G^{P^-\hgen} .$$

	\item The map $\Hom(\mbb,\mbb) \mt \Theta(\mbb) \to \Theta(\mbb)$ is provided by the $\mBA^1$-family of correspondences:
	$$
	\xyshort
	\xymatrix{
	  \Bun_G^{P^-\hgen}  & \bMap_\gen(X, \mBB \wtG^\gamma \gets \mBB \wt{P^-}^\gamma ) \ar[r] \ar[l] \ar[d] & \Bun_G^{P^-\hgen}.\\
	 & \mBA^1  
	} $$
	\end{itemize}
\end{rem}

\begin{constr} \label{constr-vinbunGPmgen}
We write
\begin{eqnarray*}
  \VinBun_G^{P^-\hgen,\gamma} &:=&  \bMap_\gen(X, \mBB \wtG^\gamma \gets \mBB \wt{P^-}^\gamma ),\\
  \VinBun_G^{\gamma} &:=&  \bMap(X, \mBB \wtG^\gamma  ).
 \end{eqnarray*}
There is a map
$$\VinBun_G^{P^-\hgen,\gamma}\to \Bun_G^{P^-\hgen} \mt_{\BunG}  \VinBun_G^{\gamma}$$
induced by the map 
$$\mBB \wt{P^-}^\gamma \simeq \mBB ( P^-\mt_G \wtG^\gamma )\to  \mBB P^-\mt_{\mBB G} \mBB \wtG^\gamma.$$
By Corollary \ref{cor-open-BPmwtG}, these maps are schematic open embeddings.
\end{constr}

\begin{constr} Recall that $\VinBun_{G,C_P} \simeq \BunP \mt_{\BunM} \BunPm$. Hence there is a unique open substack $\VinBun_{G,\lambda}^{\gamma}$ of $\VinBun_G^{\gamma}$ obtained by removing all the connected components 
$$\Bun_{P,\mu} \mt_{\Bun_{M,\mu}} \Bun_{P^-,\mu} \subset \VinBun_{G,C_P} $$
with $\mu\ne \lambda$ from its $0$-fiber. Let $\VinBun_{G,\lambda}^{P^-\hgen,\gamma}$ be the corresponding open sub-prestack. It is easy to see we can modify $\Theta$ to obtain
$$  \Theta_\lambda: \mbP_{\mBA^1} \laxto \bCorr(\PreStk_\lft)_{\all,\all}^{\open,2\on{-op}}$$
such that
	\begin{itemize}
	\item It sends $\alpha^+$ and $\alpha^-$ respectively to 
	$$\Bun_G^{P^-\hgen} \gets \Bun_{P,\lambda}^{M\hgen} \to \Bun_{M,\lambda},\; 
	\Bun_{M,\lambda} \gets \Bun_{P^-,\lambda} \to  \Bun_{G}^{P^-\hgen} .$$

	\item The map $\Hom(\mbb,\mbb) \mt \Theta_\lambda(\mbb) \to \Theta_\lambda(\mbb)$ is provided by the $\mBA^1$-family of correspondences:
	$$
	\xyshort
	\xymatrix{
	  \Bun_G^{P^-\hgen}  & \VinBun_{G,\lambda}^{P^-\hgen,\gamma}  \ar[r] \ar[l] \ar[d] & \Bun_G^{P^-\hgen}.\\
	 & \mBA^1  
	} $$

	\item The other data are induced from $\Theta$.
	\end{itemize}
\end{constr}

\begin{lem} \label{lem-quasi-compactness-dri}
	The right-lax functor $ \Theta_\lambda$ factors through $\bCorr(\PreStk_\lft)_{\on{QCAD},\all}^{\open,2\on{-op}}$.
\end{lem}

\proof We only need to check all the three left arms in the above three correspondences are QCAD. The claims for the first two arms are just (the mirror version of) Lemma \ref{lem-quasi-compactness}(1). To prove the claim for the third arm, using the open embedding in Construction \ref{constr-vinbunGPmgen}, we only need to show $\Bun_G\gets \VinBun_{G,\lambda}^\gamma$ is QCAD. It is well-known that $\VinBun_G$ is locally QCAD. Hence we only need to show $\VinBun_{G,\lambda}^\gamma\to \Bun_G$ is quasi-compact. Then we are done because both the $\mBG_m$-locus and the $0$-fiber of $\VinBun_{G,\lambda}^\gamma$ is quasi-compact over $\BunG$.

\qed[Lemma \ref{lem-quasi-compactness-dri}]

We are going to construct a Drinfeld input from $\Theta_\lambda$ by taking quotients for the torus actions. We first introduce some notations.

\begin{notn}  \label{notn-group-action}
	Let $\on{ActSch}_\ft^{\on{aff}}$ be the category whose objects are $(H\act Y)$, where $H$ is an affine algberaic group and $Y\in \affSch_\ft$. We equip $\on{ActSch}_\ft^{\on{aff}}$ with the Cartesian symmetric monoidal structure. Note that the monoidal unit for it is $(\pt\act \pt)$. Also note that there is a symmetric monoidal forgetful functor $\oblv_{\on{Act}}: \on{ActSch}_\ft^{\on{aff}} \to \affSch_\ft.$

	As in Definition \ref{defn-PmbA1}, we define a category $\mbP_{\mBG_m\act \mBA^1}$ enriched in $\on{ActSch}_\ft^{\on{aff}}$ such that 
	$$\bHom_{\mbP_{\mBG_m\act \mBA^1}}(\mbb,\mbb) = (\mBG_m\act \mBA^1)$$
	and the other three mapping objects are $ (\mBG_m \act 0)$. We use the same symbols $\alpha^+$ and $\alpha^-$ to denote the canonical morphisms.

	Note that $\mbP_{\mBA^1}$ can be obtained from $\mbP_{\mBG_m\act \mBA^1}$ by the procudure of \emph{changing of enrichment} along $\oblv_{\on{Act}}$. In particular, there is a forgetful functor $\mbP_{\mBG_m\act \mBA^1} \to  \mbP_{\mBA^1}$ that intertwines the enrichment via $\oblv_{\on{Act}}$.

	Let $\on{ActPreStk}_\lft$ be the similarly defined category. A morphism $(H\act Y_1)\to (H_2\act Y_2)$ is said to be an open embedding if $H_1\simeq H_2$ and $Y_1\to Y_2$ is a schemtaic open embedding. It is said to be QCAD if $Y_1\to Y_2$ is QCAD.
\end{notn}

\begin{constr} (c.f. \cite[$\S$ C.13.4]{drinfeld2013algebraic})

	In the previous connstruction of $\Theta_\lambda$, we ignored the various $\mBG_m$-actions. If we keep tracking them, we can obtain a weakly $\on{ActSch}_\ft^{\on{aff}}$-enriched right-lax functor 
	$$ \Theta_\lambda^{\on{Act}}:  \mbP_{\mBG_m\act \mBA^1} \laxto \bCorr(\on{ActPreStk}_\lft)_{\on{QCAD},\all}^{\open,2\on{-op}}$$
	such that
	\begin{itemize}
	\item It sends $\alpha^+$ and $\alpha^-$ respectively to 
	\begin{eqnarray*} 
	(\mBG_m\act  \Bun_G^{P^-\hgen}) \gets (\mBG_m\act \Bun_{P,\lambda}^{M\hgen}) \to (\mBG_m\act \Bun_{M,\lambda}),\\
	(\mBG_m\act \Bun_{M,\lambda}) \gets (\mBG_m\act \Bun_{P^-,\lambda}) \to  (\mBG_m\act  \Bun_G^{P^-\hgen}) .
	\end{eqnarray*}

	\item The map $\Hom(\mbb,\mbb) \mt  \Theta_\lambda^{\on{Act}}(\mbb) \to  \Theta_\lambda^{\on{Act}}(\mbb)$ is provided by the diagram:
	$$
	\xyshort
	\xymatrix{
	 (\mBG_m\act \Bun_G^{P^-\hgen})  & (\mBG_m\mt \mBG_m \act \VinBun_{G,\lambda}^{P^-\hgen,\gamma})  \ar[r] \ar[l] \ar[d] & (\mBG_m \act \Bun_G^{P^-\hgen}).\\
	 &( \mBG_m\act \mBA^1  )
	},$$
	which is induced by the morphism
	$$  (\mBG_m \mt \mBG_m \act \wt{Z} )\to (\mBG_m\act Z)\mt (\mBG_m\act Z)\mt (\mBG_m\act \mBA^1) $$
	that exists for any Drinfeld-Gaitsgory interpolation $\wt Z$ (see \cite[$\S$ 2.2.3]{drinfeld2014theorem}).

	\item It is compatible with $\Theta_\lambda$ via the forgetful functors.
	\end{itemize}

	Then as in \cite[Footnote 41]{drinfeld2013algebraic}, we obtain the desired Drinfeld input by passing to quotients and changing enrichment.
\end{constr}

\qed[Proposition-Construction \ref{propconstr-drinfeld-input}]

\section{Step 3} \label{section-diagram-chasing}

We have two results to prove in this step: Proposition-Construction \ref{prop-constr-BunPm-olBunG-gen} and Lemma \ref{lem-diagram-chasing}. Each subsection corresponds to a result.

The main idea of the proof of Proposition-Construction \ref{prop-constr-BunPm-olBunG-gen} is to find a suitable open substack of $P^-\backslash \,_0\!\Vin_{G,\ge C_P}/G $ whose fibers over $T_{\ad,\ge C_P}^+$ are all isomorphic to $\mBB P^-$. This is achieved by studying the stablizer group $\wt{G}$ (reviewed in Appendix \ref{appendix-wtG}).

The proof of Lemma \ref{lem-diagram-chasing} is a tedious diagram-chasing.

\ssec{Proof of Proposition-Construction \ref{prop-constr-BunPm-olBunG-gen}}
\label{ssec-step-3-1}
\begin{goal} Construct a certain open embedding
	$$ (\BunPm\mt_{\BunG} \ol{\Bun}_{G,\ge P})^\gen \to \BunPm\mt_{\BunG} \ol{\Bun}_{G,\ge P}$$
	whose restriction to the $G$-stratum and $P$-stratum are canonically isomorphic to the maps
	\begin{eqnarray*}
		\BunPm \mt_{  \BunG } \ol{\Bun}_{G,G} &\to & \BunPm \mt_{  \BunG } \ol{\Bun}_{G,G},\\
		\BunPm^{M\hgen} \mt_{  \BunG^{P\hgen} } \ol{\Bun}_{G,P} &\to &  \BunPm \mt_{  \BunG } \ol{\Bun}_{G,P}.
	\end{eqnarray*}
\end{goal}

	By definition, we have
	$$ \BunPm\mt_{\BunG} \VinBun_{G,\ge C_P} \simeq \bMap_\gen(X,  P^-\backslash \Vin_{G,\gCP}/G \supset  P^-\backslash \,_0\!\Vin_{G,\gCP}/G ).$$
	Note that 
	$$P^-\backslash \,_0\!\Vin_{G,\gCP}/G \simeq \mBB P^-\mt_{\mBB G}\mBB \wtG_\gCP,$$ where $\wtG_\gCP:=\wtG\mt_{\Tadp}T_{\ad,\gCP}^+$. By Corollary \ref{cor-open-BPmwtG}, the map
	$$  \mBB(  P^-\mt_G \wtG_{\gCP} ) \to \mBB P^-\mt_{\mBB G}\mBB \wtG_\gCP \simeq  P^-\backslash \,_0\!\Vin_{G,\gCP}/G$$
	is a schematic open embedding. We define
	$$(\BunPm\mt_{\BunG} \VinBun_{G,\ge C_P})^\gen := \bMap_\gen(X,  P^-\backslash \Vin_{G,\gCP}/G \gets  \mBB(  P^-\mt_G \wtG_{\gCP} )  ).$$
	Then we have a schematic open embedding
	$$
	(\BunPm\mt_{\BunG} \VinBun_{G,\ge C_P})^\gen \to  \BunPm\mt_{\BunG} \VinBun_{G,\ge C_P}.
	$$

	As in the proof of Lemma \ref{lem-unique-lifting-group-action-VinleP-to-BPBPm}, a direct calculation shows that the $Z_M$-action on $_0\!\Vin_{G,\gCP}$ preserves the open substack
	$$   \mBB (  P^-\mt_G \wtG_{\gCP} )\mt_{ (P^-\backslash \,_0\!\Vin_{G,\gCP}/G)  } \,_0\!\Vin_{G,\gCP}.$$
	Hence it makes sense to define
	$$  (\BunPm\mt_{\BunG} \ol{\Bun}_{G,\ge P})^\gen := (\BunPm\mt_{\BunG} \VinBun_{G,\ge C_P})^\gen /Z_M.$$
	It is obvious that the restriction of
	$$ (\BunPm\mt_{\BunG} \ol{\Bun}_{G,\ge P})^\gen\to \BunPm\mt_{\BunG} \ol{\Bun}_{G,\ge P} $$
	to the $G$-stratum is an isomorphism. It remains to identify its restriction to the $P$-stratum with the map 
	$$\BunPm^{M\hgen} \mt_{  \BunG^{P\hgen} } \ol{\Bun}_{G,P} \to  \BunPm \mt_{  \BunG } \ol{\Bun}_{G,P}.$$
	Unwinding the definitions, we only need to identify the $C_P$-fiber of the open embedding 
	$$\mBB(  P^-\mt_G \wtG_{\gCP} )  \to \mBB P^-\mt_{\mBB G}\mBB \wtG_\gCP$$
	with the $C_P$-fiber of the map
	$$  \mBB M\mt_{\mBB P} \mBB \wt{G}_{\le P} \to  \mBB P^-\mt_{\mBB G} \mBB \wt{G}_{\le P}.$$
	However, this follows from $\wtG_{C_P}\simeq P\mt_M P^-$.

\qed[Proposition-Construction \ref{prop-constr-BunPm-olBunG-gen}]

\ssec{Proof of Lemma \ref{lem-diagram-chasing}}
\label{ssec-step-3-2}
We will introduce many temporary notations in this subsection. When we use an english letter, like $c$, to denote a correspondence, or when we use a letter of plain font, like $K$, to denote a D-module, it means such notations are only used in this subsection.

\begin{goal} 
	The morphism $\gamma_P$ and/or $\gamma_P'$ are equivalent to the morphism
	$$ \DMod_\hol^{!,*}( \beta ) \circ \mbK(P)\to \DMod_\hol^{!,*}( \beta ) \circ \mbK(G).$$
\end{goal}

\sssec{The arrow $\gamma_P$} \label{sssec-description-gammaP}
We first give the following tautological description of
$$\gamma_P: \CT^\gen_{P\mt G,*}(\mCF_P) \to \CT_{P\mt G,*}(\mCF_G).$$
Recall the morphism (\ref{eqn-Eis-mCF-to-Del}): 
$$\Eis^\enh_{ P\mt G\to G\mt G }(\mCF_P') \to \mCF_G'.$$
Its underlying morphism in $\Dmod_\hol(\BunG\mt \BunG)$ is a map 
$$ \vartheta_P:\mfp_{P\mt G,!}^{\enh} (\mCF_P)\to \mCF_G,$$
which by adjunction induces a morphism
$$ \theta_P: \mCF_P \to  \mfp_{P\mt G}^{\enh,!} (\mCF_G).$$
Then we have
$$\gamma_P\simeq  \CT^\gen_{P\mt G,*}(\theta_P).$$
Note that we indeed have $\CT_{P\mt G,*} \simeq \CT^\gen_{P\mt G,*}\circ \mfp_{P\mt G}^{\enh,!}$.

\sssec{Second adjointness for left functors} \label{sssec-second-adjointness-leftadjoint-version}
Next, we give a more convenient description for the second adjointness, when restricted to ind-holonomic objects.

Let $\CT_{P,*}$ be the restriction of $\CT_{P,*}$ to the full subcategory of ind-holonomic objects. By construction, the natural transformation $\CT_{P,*} \simeq \, \CT_{ P^-,! }$ is obtained as follows. We apply $\DMod^{\bt,!}$ to the $2$-morphism 
$$\alpha_{P,\lambda}^{+} \circ \alpha_{P,\lambda}^{-} \to \on{Id}_{\Bun_{M,\lambda}}$$
in $\bCorr(\PreStk_\lft)_{\on{QCAD},\all}^{\open,2\on{-op}}$
and obtain a natural transformation
$$  \CT_{P,*,\lambda}\circ (\CT_{ P^-,!,\lambda})^R  \to \Id_{\Dmod_\hol(\BunM)}.$$
Then we obtain the natural transfomation $\CT_{P,*,\lambda} \to \,\CT_{ P^-,!,\lambda}$ by adjunction. Equivalently, we have the left adjoint version of the above picture. Namely, we start from the $2$-morphism\footnote{The superscript ``$\rev$'' means exchanging the two arms of a correspondence.}
$$  \on{Id}_{\Bun_{M}}\to (\alpha_{P}^{-})^\rev \circ (\alpha_{P}^{+} )^\rev $$
in $\bCorr(\PreStk_{\lft})_{\all,\on{Stacky}}^{\open}$, and use $\Dmod_\hol^{!,*}$ to obtain a natural transformation
$$ \Id\to\, \CT_{ P^-,! }\circ  (\CT_{P,*})^L.$$
Then we can obtain the same natural transfomation $\CT_{P,*} \to \,\CT_{ P^-,!}$ by adjunction.

The advantage is: if we use left functors, we can work with all the connected components simultaneously.

Similarly, the natural transformation of $\CT_{P\mt G,*} \simeq \, \CT_{ P^-\mt G,!}$ can be obtained by the same procedure from the correspondences
\begin{eqnarray*}
 c^+ &:& (  \BunG \mt \BunG \gets \Bun_{P}\mt \BunG \to \Bun_{M} \mt \BunG  ), \\
  c^- &:& (  \Bun_{M}\mt  \BunG \gets  \Bun_{P^-} \mt \BunG \to \BunG \mt \BunG ),
\end{eqnarray*}
and the $2$-morphism
\begin{equation} \label{eqn-2-morphism-beta-+-}
  \on{Id}_{ (\Bun_{M}\mt \BunG) }\to  c^- \circ  c^+
 \end{equation}
in $\bCorr(\PreStk_{\lft})_{\all,\on{Stacky}}^{\open}$.

Similarly, the natural transformation of $\CT^\gen_{P\mt G,*} \simeq \, \CT^\gen_{ P^-\mt G,!}$ can be obtained by the same procedure from the correspondences
\begin{eqnarray*}
 c^{+,\gen} &:& (  \BunG^{P\hgen} \mt \BunG \gets \Bun_{P}\mt \BunG \to \Bun_{M} \mt \BunG  ), \\
  c^{-,\gen} &:& (  \Bun_{M}\mt  \BunG \gets  \Bun_{P^-}^{M\hgen} \mt \BunG \to \BunG^{P\hgen} \mt \BunG ).
\end{eqnarray*}
and the $2$-morphism
\begin{equation} \label{eqn-2-morphism-beta-+--gen}
 \on{Id}_{(\Bun_{M}\mt \BunG)}\to  c^{-,\gen} \circ  c^{+,\gen}
\end{equation}
in $\bCorr(\PreStk_{\lft})_{\all,\on{Stacky}}^{\open}$.

\begin{notn} To simplify the notations, for a correspondence $c$ (in english letter), we use the symbol $\bf{c}$ to denote the corresponding functor $\DMod_\hol^{!,*}(c )$. These shorthands are only used in this subsection.
\end{notn}

\sssec{Translation}
Using the above shorthands, the results in $\S$ \ref{sssec-second-adjointness-leftadjoint-version} are translated as below. The $2$-morphisms (\ref{eqn-2-morphism-beta-+-}) and (\ref{eqn-2-morphism-beta-+--gen}) induce natural transformations\footnote{The functor $\Id$ below is the identity functor for $\Dmod_\hol(\BunM\mt \BunG)$.} 
$$\bf \mu: \Id \to  c^-\circ c^+ \;\on{ and }\; \mu^\gen: \Id \to c^{-,\gen}\circ c^{+,\gen}$$
such that the following compositions are isomorphisms
\begin{eqnarray} \label{eqn-sadj-1}
    (\mathbf{c}^{+,\gen})^R \xto{\mu^\gen} & \mathbf{c}^{-,\gen}\circ \mathbf{c}^{+,\gen} \circ (\mathbf{c}^{+,\gen})^R  &\xto{\counit} \mathbf{c}^{-,\gen} \\\label{eqn-sadj-2}
	 (\mathbf{c}^{+})^R \xto{\mu} & \mathbf{c}^{-}\circ \mathbf{c}^{+} \circ (\mathbf{c}^{+})^R & \xto{\counit}  \mathbf{c}^{-} .
\end{eqnarray}

\sssec{}
Consider the map $\mfp_{P\mt G}^{\enh}:  \BunG^{P\hgen} \mt \BunG\to   \BunG \mt \BunG$. Let 
$$p: (   \BunG \mt \BunG \gets   \BunG^{P\hgen} \mt \BunG \xto{=} \BunG^{P\hgen} \mt \BunG )$$
be the corresponding correspondence. Note that we have $ \mbp \simeq \mfp_{P\mt G,!}^\enh$.

By definition, we have $c^+ \simeq p\circ c^{+,\gen}$, which provides
$$ \mbc^+ \simeq \mbp\circ \mbc^{+,\gen}. $$
We proved in $\S$ \ref{ssec-step-2-1} that the map $\Bun_{P^-}^{M\hgen} \to \BunPm \mt_{\BunG} \BunG^{P\hgen}$ is a schematic open embedding. Hence we also have a $2$-morphism $c^{-,\gen} \to  c^- \circ p$, which provides
$$ \nu: \mbc^{-,\gen} \to \mbc^-\circ\mbp. $$

By constuction, the $2$-morphism (\ref{eqn-2-morphism-beta-+-}) is equivalent to the composition
$$  \on{Id}_{(\Bun_{M}\mt \BunG)} \xto{ (\ref{eqn-2-morphism-beta-+--gen}) }   c^{-,\gen} \circ  c^{+,\gen} \to  c^- \circ p\circ  c^{+,\gen} \simeq  c^- \circ c^+.$$
Hence $\mu$ is isomorphic to
\begin{equation} \label{eqn-compare-unit-sadj}
 \Id\xto{\mu^\gen}  \mbc^{-,\gen}\circ \mbc^{+,\gen} \xto{\nu(\mbc^{+,\gen})}\mbc^-\circ \mbp \circ\mbc^{+,\gen} \simeq \mbc^-\circ \mbc^+.  
\end{equation}

\begin{lem} \label{lem-description-gamma_p'}
The arrow $\gamma_P$ is equivalent to the composition
\begin{eqnarray*}
 (\mbc^{-,\gen})(\mCF_P) \xto{\nu}   \mbc^-  \circ \mbp (\mCF_P)  \xto{  \mbc^- (\vartheta_P) }   \mbc^- (\mCF_G).
\end{eqnarray*}
\end{lem}

\proof 

By definition, the arrow $\theta_P:\mCF_P \to \mbp^R(\mCF_G)$ is isomorphic to
$$ \mCF_P\xto{\unit} \mbp^R\circ \mbp(\mCF_P) \xto{\mbp^R(\vartheta_P)}  \mbp^R(\mCF_G).$$
Hence by definition, $\gamma_P$ is isomorphic to
$$ (\mbc^{+,\gen})^R (\mCF_P)\xto{\unit} (\mbc^{+,\gen})^R \circ \mbp^R\circ \mbp(\mCF_P) \simeq  (\mbc^+)^R\circ \mbp(\mCF_P)   \xto{ (\mbc^+)^R(\vartheta_P)} (\mbc^+)^R(\mCF_G).$$
Hence we only need to show the following diagram of functors commute
\begin{equation} \label{eqn-diagram-chasing-functors}
\xyshort
\xymatrix{
	(\mbc^{+,\gen})^R  \ar[r]^-{\unit} \ar[d]^-{\simeq}_-{(\ref{eqn-sadj-1})} &
	(\mbc^{+,\gen})^R \circ \mbp^R\circ \mbp \ar[r]^-\simeq &
	(\mbc^+)^R\circ \mbp  \ar[d]^-{\simeq}_-{(\ref{eqn-sadj-2})} \\
	\mbc^{-,\gen} \ar[rr]^-\nu & & \mbc^-\circ \mbp.
}
\end{equation}
Note that we have
$$\Map( (\mbc^{+,\gen})^R, \mbc^-\circ \mbp ) \simeq \Map( \Id , \mbc^-\circ \mbp \circ \mbc^{+,\gen} ) \simeq \Map( \Id , \mbc^-\circ \mbc^{+} ).$$
Via this isomorphism, the top arc in (\ref{eqn-diagram-chasing-functors}), which is a point of the LHS, is given by the following point of the RHS:
\begin{eqnarray*}
 \Id \xto{\unit}  (\mbc^{+,\gen})^R\circ \mbc^{+,\gen} \xto{\unit}  (\mbc^{+,\gen})^R\circ  \mbp^R \circ  \mbp \circ \mbc^{+,\gen} \simeq (\mbc^{+})^R\circ \mbc^{+} \to \\
 \xto{\mu} \mbc^{-}\circ \mbc^{+} \circ (\mbc^{+})^R\circ \mbc^{+}  \xto{\counit}  \mbc^-\circ \mbc^{+}.
\end{eqnarray*}
The first row in the above composition is just $\unit:\Id\to (\mbc^{+})^R\circ \mbc^{+}$. Hence this composition is isomorphic to
$$ \Id \xto{\mu} \mbc^{-}\circ \mbc^{+} \xto{\unit} \mbc^{-}\circ \mbc^{+}\circ  (\mbc^{+})^R\circ \mbc^{+} \xto{\counit}   \mbc^-\circ \mbc^{+}, $$
which is just $\Id\to \mbc^{-}\circ \mbc^{+}$ by the axioms for $\unit$ and $\counit$.

Similarly, one shows that the bottom arc corresponds to natural transformation (\ref{eqn-compare-unit-sadj}). Then we are done by the discussion above the lemma.

\qed[Lemma \ref{lem-description-gamma_p'}]

\sssec{Finish of the proof}
We give temporary labels to the following correspondences
\begin{eqnarray*}
	i: (\ol{\Bun}_{G, \ge P}  \gets &  \ol{\Bun}_{G,P} & \xto{=}  \ol{\Bun}_{G, P}),\\
	d^\gen:( \Bun_G^{P\hgen}\mt \BunG  \gets &  \ol{\Bun}_{G,P} & \xto{=}  \ol{\Bun}_{G, P}), \\
	d:( \Bun_G\mt \BunG  \gets &  \ol{\Bun}_{G,\ge P} & \xto{=}  \ol{\Bun}_{G, \ge P}),\\
	j: (\ol{\Bun}_{G, \ge P} \xgets{=}  & \ol{\Bun}_{G, \ge P} & \to  \ol{\Bun}_{G}),\\
	b: ( \BunM\mt \BunG   \gets &  (\BunPm\mt_{\BunG} \ol{\Bun}_{G,\ge P})^\gen  & \to  \ol{\Bun}_{G, \ge P}).
\end{eqnarray*}
Note that we have an obvious isomorphism $\beta\simeq b\circ j$, hence
$$ \Dmod_\hol^{!,*}(\beta) \simeq \mbb \circ \mbj.$$
We have an isomorphism $p\circ d^\gen \simeq d\circ i$ because both sides are just 
$$\Bun_G\mt \BunG \gets \ol{\Bun}_{G,P} \xto{=}  \ol{\Bun}_{G, P}.$$
Hence $\mbp \circ \mbd^\gen \simeq \mbd\circ \mbi$. We have an isomorphism $b\circ i \simeq c^{-,\gen} \circ d^\gen $ because both sides are just
$$  \BunM\mt \BunG \gets \BunPm^{M\hgen}\mt_{\BunG^{P\hgen}} \ol{\Bun}_{G,P} \to \ol{\Bun}_{G,P}.$$
Hence $\mbb\circ \mbi\simeq  \mbc^{-,\gen} \circ \mbd^\gen$. We have a $2$-morphism $b\to c^-\circ d$ induced by the open embedding
$$  (\BunPm\mt_{\BunG} \ol{\Bun}_{G,\ge P})^\gen \subset \BunPm\mt_{\BunG} \ol{\Bun}_{G,\ge P}.  $$
Hence we have a natural transformation
$$ \xi:  \mbb \to \mbc^- \circ \mbd.$$
Moreover, the $2$-morphism
$$  b\circ i \simeq c^{-,\gen} \circ d^\gen \to c^-\circ p\circ d^\gen \simeq c^-\circ d\circ i $$
is isomorphic to the $2$-morphism induced from $b\to c^-\circ d$. Hence we have the following commutative diagram of functors
\begin{equation} \label{eqn-compare-2-morphism}
\xyshort
\xymatrix{
	\mbb \circ \mbi \ar[d]^-\simeq \ar[r]^-{ \xi(\mbi) } &
	  \mbc^-\circ \mbd\circ \mbi \ar[d]^-\simeq  \\
	\mbc^{-,\gen}\circ \mbd^\gen \ar[r]^-{\nu(\mbd^\gen )} &
	\mbc^-\circ \mbp \circ \mbd^\gen.
}
\end{equation}

After these preparations, we are ready to finish the proof. Recall that $\mbK(P)$ is a $!$-extension along $\ol{\Bun}_{G,P}\to \ol{\Bun}_G$. Let $K_1$ be the corresponding object in $\Dmod_\hol(\ol{\Bun}_{G,P})$. We also write $K_2:=j_{\ge P}^*(\mbK(G))$, where $j_{\ge P}^*: \ol{\Bun}_{G,\ge P}\to \ol{\Bun}_G $ is the open embedding. The morphism $\mbK(P) \to \mbK(G)$ is sent by $\mbj=j_{\ge P}^*$ to a morphism 
$$\eta:  \mbi( K_1) \to K_2.$$
It follows from definition that the arrow $\vartheta_P:\mbp(\mCF_P) \to \mCF_G$ is equivalent to
$$  \mbp \circ \mbd^\gen ( K_1) \simeq \mbd\circ \mbi (K_1) \xto{\mbd(\eta)}  \mbd(K_2) , $$
where $\mCF_P \simeq  \mbd^\gen ( K_1)$ and $\mCF_G \simeq  \mbd(K_2)$. Hence by Lemma \ref{lem-description-gamma_p'}, the arrow $\gamma_P$ is equivalent to
$$ \mbc^{-,\gen}\circ \mbd^\gen (K_1) \xto{  \nu(\mbd^\gen )  } \mbc^-\circ \mbp \circ \mbd^\gen(K_1) \simeq  \mbc^-\circ \mbd\circ \mbi(K_1) \xto{\mbc^-\circ \mbd(\eta)} \mbc^-\circ \mbd(K_2). $$
By (\ref{eqn-compare-2-morphism}), this arrow is equivalent to 
$$  \mbb\circ \mbi(K_1) \xto{\xi(\mbi(K_1))}  \mbc^-\circ \mbd\circ \mbi(K_1) \xto{\mbc^-\circ \mbd(\eta)} \mbc^-\circ \mbd(K_2),$$
or equivalently
$$  \mbb\circ \mbi(K_1) \xto{\mbb(\eta)} \mbb(K_2) \xto{\xi(K_2)}  \mbc^-\circ \mbd(K_2).$$
We claim $\xi(K_2)$ is invertible. Indeed, this is because $K_2$ is a $!$-extension from the $G$-stratum, and the open embedding
$$  (\BunPm\mt_{\BunG} \ol{\Bun}_{G,\ge P})^\gen \subset \BunPm\mt_{\BunG} \ol{\Bun}_{G,\ge P}$$
is an isomorphism when restricted to the $G$-stratum. Hence $\gamma_P$ is equivalent to $\mbb(\eta)$, which by definition is the image of $\mbK(P) \to \mbK(G)$ under $\mbb\circ \mbj\simeq  \Dmod_\hol^{!,*}(\beta) $.

\qed[Lemma \ref{lem-diagram-chasing}]

\section{Step 4} \label{section-restoring-symmetry}

We have two results to prove in this step: Proposition-Construction \ref{propconstr-step-4-1} and Proposition \ref{prop-step-4-2}. Each subsection corresponds to a result.

To avoid jumping between topics, we also prove Lemma \ref{lem-appear-local-model} (from Step 5) in $\S$ \ref{ssec-step-4-1}.

The proofs are similar to those in Step 1: Proposition-Construction \ref{propconstr-step-4-1} is similar to Proposition-Construction \ref{propconstr-factor-through-gen}, while Proposition \ref{prop-step-4-2} is similar to Proposition \ref{prop-FP-equivariant}.

\ssec{Proof of Proposition-Construction \ref{propconstr-step-4-1} and Lemma \ref{lem-appear-local-model}}
\label{ssec-step-4-1}
\begin{goal} Construct a canonical factorization of the map 
	$$(\BunPm\mt_{\BunG} \ol{\Bun}_{G,\ge P})^\gen \to \BunM\mt \BunG$$
	via $\BunM\mt \Bun_G^{P^-\hgen}$ such that we have an isomorphism
	\begin{equation} \label{eqn-fiber-product-local-model}
	(\BunPm\mt_{\BunG} \ol{\Bun}_{G,\ge P})^\gen  \mt_{ \Bun_G^{P^-\hgen} } \BunP^{M\hgen} \simeq Y^P_\rel/Z_M
	\end{equation}
	defined over $\BunM \mt \ol{\Bun}_G \mt \BunM$.
\end{goal}

The proof below is similar to that in $\S$ \ref{ssec-step-1-1}. Hence we omit some details.

Recall in $\S$ \ref{ssec-step-3-1}, we defined 
$$(\BunPm\mt_{\BunG} \VinBun_{G,\ge C_P})^\gen := \bMap_\gen(X,  P^-\backslash \Vin_{G,\gCP}/G \gets  \mBB(  P^-\mt_G \wtG_{\gCP} )  ).$$
By Lemma \ref{lem-Pm-wtG-gCP}, the \emph{right} projection map $P^-\mt_G \wtG_{\gCP}\to G$ factors through $P^-$. Hence we obtain the following commutative diagram of algebraic stacks
$$
\xyshort
\xymatrix{
	 P^-\backslash \Vin_{G,\gCP}/G \ar[d] &
	 \mBB(  P^-\mt_G \wtG_{\gCP} ) \ar[d]  \ar[l] \\
	 \mBB M\mt \mBB G  & \mBB M\mt \mBB P^- \ar[l].
}
$$
Taking $\bMap_\gen(X,-)$, we obtain a map
$$
 (\BunPm\mt_{\BunG} \VinBun_{G,\ge C_P})^\gen \to \BunM\mt \BunG^{P\hgen}.
$$
To obtain the map $(\BunPm\mt_{\BunG} \ol{\Bun}_{G,\ge P})^\gen \to \BunM\mt \BunG^{P\hgen}$, as before, we show that the map $ \mBB(  P^-\mt_G \wtG_{\gCP} ) \to \mBB M\mt \mBB P^- $ can be uniquely lifted to a morphism
$$ (Z_M\act   \mBB(  P^-\mt_G \wtG_{\gCP} ) ) \to ( \pt \act \mBB M\mt \mBB P^-  ) $$
fitting into the diagram
$$
\xyshort
\xymatrix{
	(Z_M\act P^-\backslash \Vin_{G,\gCP}/G ) \ar[d] &
	(Z_M\act \mBB(  P^-\mt_G \wtG_{\gCP} ) ) \ar[d]  \ar[l] \\
	(\pt \act \mBB M\mt \mBB G)  &  (\pt\act  \mBB M\mt \mBB P^-). \ar[l]
}
$$

It remains to compare both sides of (\ref{eqn-fiber-product-local-model}). By construction, 
$$(\BunPm\mt_{\BunG} \VinBun_{G,\ge C_P})^\gen \mt_{ \Bun_G^{P^-\hgen} } \BunP^{M\hgen}$$
is isomorphic to the image of 
$$   P^-\backslash \Vin_{G,\gCP}/P \gets   \mBB(  P^-\mt_G \wtG_{\gCP} )\mt_{\mBB P^-} \mBB M $$
under the functor $\bMap_\gen(X,-)$. Using Lemma \ref{lem-Pm-wtG-gCP}(1), the map
$$ \mBB (  P^-\mt_G \wtG_{\gCP} \mt_{P^-} M )\to \mBB(  P^-\mt_G \wtG_{\gCP} )\mt_{\mBB P^-} \mBB M $$ 
is an isomorphism. Also, the LHS is just
$$\mBB (  P^-\mt_G \wtG_{\gCP} \mt_{P^-} M )\simeq \mBB ( P^-\mt_G \wt{G}_\gCP \mt_G P )\simeq P^-\backslash \Vin_{G,\gCP}^\Bru /P^-.$$
Hence we obtain a $Z_M$-equivariant isomorphism
$$ (\BunPm\mt_{\BunG} \VinBun_{G,\ge C_P})^\gen \mt_{ \Bun_G^{P^-\hgen} } \BunP^{M\hgen} \simeq Y^P_\rel.$$
It follows from construction that it is defined over $\BunM \mt \VinBun_G \mt \BunM$. Then we obtain the isomorphism (\ref{eqn-fiber-product-local-model}) by taking quotients for the $Z_M$-actions.
 
\qed[Proposition-Construction \ref{propconstr-step-4-1} and Lemma \ref{lem-appear-local-model}]

\ssec{Proof of Proposition \ref{prop-step-4-2}}
\label{ssec-step-4-2}
\begin{goal} The objects $\DMod_\hol^{!,*}( \beta' ) \circ \mbK(P)$ and $\DMod_\hol^{!,*}( \beta' ) \circ \mbK(G)$ are both contained in the full subcategory
$$ \mbI(M\mt G, M\mt P^- ) \subset \Dmod( \BunM\mt \Bun_G^{P^-\hgen} ). $$
\end{goal}

	We first prove the claim for the second object. Using the base-change isomorphisms, it is easy to see $\DMod_\hol^{!,*}( \beta' ) \circ \mbK(G)$ is isomorphic to the image of $k_{\BunPm}$ under the $!$-pushforward functor along
	$$ \BunPm \to \BunM \mt \BunG^{P^-\hgen}. $$
	This map has a factorization
	$$ \BunPm \xto{f} \BunM \mt \BunPm \to \BunM \mt \BunG^{P^-\hgen}.$$
	It is clear that $f_!(k_{\BunPm})$ is $\mbU_P^-$-equivariant, i.e., is a $*$-pullback along
	$$  \BunM \mt \BunPm \to \BunM \mt \BunM.$$
	Then we are done by applying Proposition \ref{prop-well-defined-Eis_enh}(1) to the reductive group $M\mt G$ and the parabolic subgroup $M\mt P^-$.

	Now we prove the claim for the first object. Consider the restriction of $\beta$ on the $P$-stratum:
	$$ \beta_P': ( \BunM \mt \BunG^{P^-\hgen} \gets \BunPm^{M\hgen} \mt_{  \BunG^{P\hgen} } \ol{\Bun}_{G,P} \to \ol{\Bun}_{G,P}).$$
	Since $\mbK(P)$ is a $!$-extension along $i_P:\ol{\Bun}_{G,P}\to \ol{\Bun}_G$, 
	$$\DMod_\hol^{!,*}( \beta' ) \circ \mbK(P) \simeq \DMod_\hol^{!,*}(\beta_P')\circ i_P^*(\mbK(P)).$$
	It follows from construction that $\beta_P'$ is isormophic to the composition of
	$$  \Bun_G^{P\hgen}\mt \Bun_G^{P^-\hgen} \xgets{\ol{\Del}^\enh_P} \ol{\Bun}_{G,P} \xto{=}  \ol{\Bun}_{G,P} $$
	and
	$$ \delta^-: (\BunM \mt \Bun_G^{P^-\hgen} \gets  \BunPm^{M\hgen}\mt \Bun_G^{P^-\hgen} \to  \Bun_G^{P\hgen}\mt \Bun_G^{P^-\hgen}),$$
	where the map $\ol{\Del}^\enh_P$ is provided by Proposition-Construction \ref{propconstr-factor-through-gen}.
	Hence we only need to show 
	\begin{equation}\label{eqn-proof-prop-step-4-2-1}
	\Dmod_\hol^{!,*}(\delta^-)\circ \ol{\Del}^\enh_{P,!}\circ i_P^*(\mbK(P)) 
	\end{equation}
	is contained in $\mbI(M\mt G,M\mt P^-)$. In other words, we need to show its $!$-pullback along
	$$ \iota_{M\mt P^-}: \BunM \mt \BunPm \to \BunM \mt \Bun_G^{P^-\hgen} $$
	is $\mbU^-_P$-equivariant.

	Consider the correspondence
	$$ \delta^+: ( \BunM \mt \Bun_G^{P^-\hgen} \gets  \BunP\mt \Bun_G^{P^-\hgen} \to  \Bun_G^{P\hgen}\mt \Bun_G^{P^-\hgen}).$$
	As before, $\Dmod_\hol^{!,*}(\delta^-)$ is isomorphic to the restriction of $\Dmod^{\bt,!}(\delta^+)$. Hence we can rewrite (\ref{eqn-proof-prop-step-4-2-1}) as $\Dmod^{\bt,!}(\delta^+)\circ \ol{\Del}^\enh_{P,!}\circ i_P^*(\mbK(P))$.

	Consider the correspondence
	$$   e: (\BunM \mt \BunPm  \gets \BunP\mt \BunPm \to   \Bun_G^{P\hgen}\mt \Bun_G^{P^-\hgen}).$$
	By the base-change isomorphisms, the functor $ \iota_{M\mt P^-}^!\circ \Dmod^{\bt,!}(\delta^+)$ is just $\Dmod^{\bt,!}(e)$. Hence we only need to show $\ol{\Del}^\enh_{P,!}\circ i_P^*(\mbK(P))$ is contained in the full subcategory
	$$ \mbI(G\mt G,P\mt P^-) \subset \Dmod(  \Bun_G^{P\hgen}\mt \Bun_G^{P^-\hgen} ). $$
	Now this can be proved similarly to that in $\S$ \ref{ssec-step-1-2}. Namely, one only needs to replace Lemma \ref{lem-!-push-equivariant-to-IGP} by the following lemma, whose proof is similar. 

	\begin{lem}
	The functor 
	$$\ol{\Delta}_{P,!}^{\enh}:\Dmod_\hol(\ol{\Bun}_{G,P}) \to \Dmod_\hol( \Bun_G^{P\hgen} \mt \BunG^{P^-\hgen} )$$
	sends objects contained in $\Dmod_\hol(\ol{\Bun}_{G,P})^{\mbU_P\mt \mbU_{P}^-}$ into objects contained in $\mbI(G\mt G,P\mt P^-)$.
	\end{lem}

	\qed[Proposition \ref{prop-step-4-2}]

\section{Step 5} \label{section-interplay}

We have two results to prove in this step: Proposition-Construction \ref{propconstr-contration} and Theorem \ref{thm-interplay} (Lemma \ref{lem-appear-local-model} was proved in $\S$ \ref{ssec-step-4-1}). Each subsection corresponds to a result.

The main ingredient of the proof of Proposition-Construction \ref{propconstr-contration} is the monoid $\ol{M}$ associated to $P$ in \cite{wang2017reductive}. For parabolic subgroups $P\subset Q$, we will construct a map $\ol{L}\to \ol{M}$ from the monoid associated to $Q$ to that associated to $P$. This map will allow us to construct the desired correspondence from $H_{ L,G\on{-pos}}/Z_L $ to $H_{ M,G\on{-pos}}/Z_M $.

The proof of Theorem \ref{thm-interplay} is similar to that of Theorem \ref{thm-generic-second-adjointness}. Namely, we prove it using Drinfeld's framework.

\ssec{Proof of Proposition-Construction \ref{propconstr-contration}}
\label{ssec-contraction}

\begin{goal}
The correspondence
$$ \psi_{P,\ge Q}: (   H_{\MGPos}/Z_M \gets   (Y^P_{\rel}/Z_M)_{\ge Q}  \to   \ol{\Bun}_{G,\ge Q}  ) $$
is isomorphic to the composition of 
$$\psi_{Q,\ge Q}: (  H_{ L,G\on{-pos}}/Z_L    \gets   Y^Q_{\rel}/Z_L   \to   \ol{\Bun}_{G,\ge Q}  )$$
by a certain correspondence from $H_{ L,G\on{-pos}}/Z_L $ to $H_{ M,G\on{-pos}}/Z_M $.
\end{goal}

Let us first recall the definition of the map $Y^P_{\rel}\to H_{\MGPos}$.

\begin{constr} \label{constr-project-local-model-hecke}
By theorem \cite[Theorem 4.2.10]{wang2017reductive}, the closed subscheme $\ol{M}\inj \Vin_{G,C_P}$ is equal to $\mfs(C_P)\cdot \Vin_{G,C_P}\cdot \mfs(C_P)$, where $\mfs:T_\ad^+ \to \Vin_G$ is the canonical section. Hence the map
$$ \Vin_{G,\ge C_P} \to \Vin_{G,C_P},\;  x\mapsto   \mfs(C_P)\cdot x \cdot \mfs(C_P)$$
factors through $\ol{M}$. It is easy to see the obtained map $\Vin_{G,\ge C_{P}} \to \ol{M}$ intertwines the actions of $P^-\mt P\to  M\mt M $ and is $Z_M$-equivariant\footnote{This $Z_M$-action on $\Vin_{G,\ge C_P}$ is induced by the $T$-action on $\Vin_G$.}. Moreover, the map
$$ P^-\backslash \Vin_{G,\ge C_P}/P \to  M \backslash \ol{M}/M $$
sends $P^- \backslash  \Vin_{G,C_P}^\Bru  / P$ into $ M \backslash M / M$. Hence we obtain a map
$$ \mfq^-_{P,\Vin}:  Y^P_{\rel}/Z_M \to H_{\MGPos}/Z_M.$$
\end{constr}

\begin{notn} Recall $L$ is the Levi subgroup of $Q$. We write $Z_L$ for the center of $L$. Let $P_L=P\cap L$ and $P^-_L=P^-\cap L$ be the parabolic subgroups of $L$ corresponding to $P$ and $P^-$. Let $L^{P\on{-}\Bru}$ be the open Bruhat cell $P^-_LP_L$ in $L$. 
\end{notn}

\begin{notn} \label{notn-parabolic-stratification-on-local-model}
The projection map
$$Y^P_{\rel}/Z_M \to T_{\ad,\gCP}^+/Z_M \simeq T_{\ad,\ge P}^+/T$$
induces a stratification on $Y^P_{\rel}/Z_M$ labelled by the poset $\Par_{\ge P}$. As usual, for $Q\in \Par_{\ge P}$ we use the notation:
$$ (Y^P_{\rel}/Z_M)_{\ge Q} :=(Y^P_{\rel}/Z_M) \mt_{ T_{\ad,\ge P}^+/T  } (T_{\ad,\ge Q}^+/T).$$
The stack
$$  Y^P_{\rel,\ge C_Q}:= Y^P_{\rel}\mt_{T_{\ad,\gCP}^+} T_{\ad,\ge C_Q}^+ $$
inherits a $Z_L$-action from the $Z_M$-action on $Y^P_{\rel}$. Note that we have an isomorphism
$$  Y^P_{\rel,\ge C_Q}/Z_L \simeq (Y^P_{\rel}/Z_M)_{\ge Q}. $$
\end{notn}

\begin{constr}

By construction, we have
$$Y^P_{\rel,\ge C_Q} \simeq  \bMap_\gen (X ,  P^-\backslash \Vin_{G,\ge C_Q}/P  \supset  P^- \backslash  \Vin_{G,\ge C_Q}^{P\on{-}\Bru}  / P  ),$$
where
$$ \Vin_{G,\ge C_Q}^{P\on{-}\Bru} := \Vin_{G,\ge C_P}^{\Bru} \cap \Vin_{G,\ge C_Q}.$$
Note that the open locus $\Vin_{G,\ge C_Q}^{P\on{-}\Bru}$ is contained in $ \Vin_{G,\ge C_Q}^{\Bru}$. Indeed, the former is the $(P^-\mt P)$-orbit of the canonical section, while the latter is the $(Q^-\mt Q)$-orbit. Hence the map
$$  P^-\backslash \Vin_{G,\ge C_Q}/P \to Q^-\backslash \Vin_{G,\ge C_Q}/Q  $$
induces a $Z_L$-equivariant map $Y^P_{\rel,\ge C_Q}\to Y^Q_{\rel}$. Hence we obtain a map
$$ \pi_{P,Q}:  (Y^P_{\rel}/Z_M)_{\ge Q} \simeq   Y^P_{\rel,\ge C_Q}/Z_L \to Y^Q_{\rel}/Z_L  .$$
By construction, we have the following commutative diagram
\begin{equation} \label{eqn-local-model-functorial}
\xyshort
\xymatrix{
	  Y^P_{\rel}/Z_M  \ar[d]^-{ \mfp^-_{P,\Vin} } &  (Y^P_{\rel}/Z_M)_{\ge Q} \ar[d]  \ar[r]^-{\pi_{P,Q}} \ar[l]_-{\supset} &
	Y^Q_{\rel}/Z_L \ar[d]^-{\mfp^-_{Q,\Vin}}  \\
	\ol{\Bun}_{G,\ge P}  & \ol{\Bun}_{G,\ge Q}\ar[l]_-\supset & \ol{\Bun}_{G,\ge Q},\ar[l]_-=
}
\end{equation}
where the left square is Cartesian.
\end{constr}

\begin{propconstr} \label{propconstr-functorial-local-model-hecke}
Consider the lft algebraic stack\footnote{When $Q=G$, $W_{P,G}$ is just the open Zastava stack. When $Q=P$, $W_{P,P}$ is $H_\MGPos$.}
$$ W_{P,Q}:=  \bMap_\gen(X, P^-_L \backslash \ol{L} / P_L \supset P^-_L \backslash L^{P\on{-}\Bru} / P_L ).$$
Then there exists a canonical commutative diagram
	\begin{equation} \label{eqn-functorial-proj-local-model-hecke}
	\xyshort
	\xymatrix{
		Y^P_{\rel}/Z_M \ar[d]^-{ \mfq^-_{P,\Vin}    }  & (Y^P_{\rel}/Z_M)_{\ge Q}  \ar[l]_-{\supset} \ar[r]^-{\pi_{P,Q}} \ar[d] &
		Y^Q_{\rel}/Z_L \ar[d]^-{ \mfq^-_{Q,\Vin}} \\
		H_{\MGPos}/Z_M & W_{P,Q}/Z_L \ar[r] \ar[l] & H_{ L,G\on{-pos}}/Z_L 
	}
	\end{equation}
such that the right square in it is Cartesian.
\end{propconstr}

\proof 
Via the identification 
$$P^-_L \backslash \ol{L} / P_L  \simeq  \mBB P^-_L \mt_{\mBB L} (L\backslash \ol{L}/L)  \mt_{\mBB L}   \mBB P_L,$$
the open substack $ P^-_L \backslash L^\Bru / P_L$ of the LHS is contained in the open substack $\mBB P^-_L \mt_{\mBB L}  \mBB L  \mt_{\mBB L}   \mBB P_L $ of the RHS. Hence we obtain a $Z_L$-equivariant schematic open embedding
$$ W_{P,Q} \to   \Bun_{P_L^-}\mt_{\Bun_L}  H_{L,G\on{-pos}} \mt_{\Bun_L}  \Bun_{P_L}.$$
In particular, we obtain a map
\begin{equation} \label{propconstr-functorial-local-model-hecke-1}
W_{P,Q}/Z_L\to  H_{L,G\on{-pos}}/Z_L.
\end{equation}

As explained in Construction \ref{constr-project-local-model-hecke}, the map
$$ \Vin_{G,\ge C_Q} \to \Vin_{G,C_Q},\;  x\mapsto   \mfs(C_Q)\cdot x \cdot \mfs(C_Q)$$
factors through $\ol{L}$. It is easy to see the obtained map $\Vin_{G,\ge C_Q} \to \ol{L}$ intertwines the actions of $Q^-\mt Q\to  L\mt L $ and is $Z_L$-equivariant\footnote{This $Z_L$-action on $\Vin_{G,\ge C_Q}$ is induced by the $T$-action on $\Vin_G$.}. Moreover, the map
$$ P^-\backslash \Vin_{G,\ge C_Q}/P  \to P_L^- \backslash \ol{L}/P_L $$
sends the $P$-Bruhat cell to the $P$-Bruhat cell. Hence we obtain a $Z_L$-linear map $Y^P_{\rel,\ge C_Q}\to W_{P,Q}$. By taking quotient, we obtain a map
\begin{equation} \label{propconstr-functorial-local-model-hecke-2} 
(Y^P_{\rel}/Z_M)_{\ge Q}\simeq Y^P_{\rel,\ge C_Q} /Z_L \to W_{P,Q}/Z_L.
\end{equation}

Note that we have $\mfs(C_P)\cdot x \cdot \mfs(C_P) = \mfs(C_P)\cdot \mfs(C_Q)\cdot x \cdot \mfs(C_Q) \cdot \mfs(C_P)$ for $x\in \Vin_{G,\ge C_Q}$. Hence the composition
$$\Vin_{G,\ge C_Q} \to \Vin_{G,\gCP} \to \ol{M} $$
factors through $\ol{L}$. Since the above composition intertwines the action of $P^-\mt P \to M\mt M$ and is $Z_L$-equivariant, the obtained map $\ol{L}\to \ol{M}$ intertwines the actions of $P_L^-\mt P_L \to M\mt M$ and is $Z_L$-equivariant. Moreover, the map
$$ P^-_L \backslash \ol{L} / P_L  \to M\backslash \ol{M}/M  $$
sends the $P$-Bruhat cell into $M\backslash M/M$. Hence we obtain a map
\begin{equation} \label{propconstr-functorial-local-model-hecke-3} 
  W_{P,Q}/Z_L \to H_{\MGPos}/Z_L \to H_{\MGPos}/Z_M. 
 \end{equation}

It follows from constructions that the above maps (\ref{propconstr-functorial-local-model-hecke-1}), (\ref{propconstr-functorial-local-model-hecke-2}) and (\ref{propconstr-functorial-local-model-hecke-3})
fit into a commutative diagram (\ref{eqn-functorial-proj-local-model-hecke}). It remains to show its right square is Cartesian. We only need to show the maps
\begin{eqnarray*}
	P^-\backslash \Vin_{G,\ge C_Q}/P &\to& (Q^-\backslash \Vin_{G,\ge C_Q}/Q) \mt_{ ( L\backslash \ol{L}/L)  } (P_L^- \backslash \ol{L}/P_L) ,\\
	P^-\backslash \Vin_{G,\ge C_Q}^{P\on{-}\Bru}/P &\to& (Q^-\backslash \Vin_{G,\ge C_Q}^\Bru/Q) \mt_{ ( L\backslash L/L)  } (P_L^- \backslash L^{P\on{-}\Bru} /P_L)
\end{eqnarray*}
are isomorphisms. To prove the claim for the first map, we only need to show $\mBB P \simeq \mBB Q \mt_{\mBB L} \mBB P_L$, but this follows from the fact that $Q\to L$ is surjective. The claim for the second map follows from the fact that the maps
\begin{eqnarray*}
	P^-\backslash \Vin_{G,\ge C_Q}^{P\on{-}\Bru}/P &\to& M\backslash M/M \mt T_{\ad,\ge C_Q}, \\
   Q^-\backslash \Vin_{G,\ge C_Q}^\Bru/Q &\to& L\backslash L/L \mt T_{\ad,\ge C_Q}, \\
   P_L^- \backslash L^{P\on{-}\Bru} /P_L &\to& M\backslash M/M.
\end{eqnarray*}
are all isomorphisms.

\qed[Proposition-Construction \ref{propconstr-functorial-local-model-hecke}]

\sssec{Proof of Proposition-Construction \ref{propconstr-contration}}
The desired correspondence is 
$$ H_{ M,G\on{-pos}}/Z_M\gets  W_{P,Q}/Z_L \to H_{L,G\on{-pos}}/Z_L.$$
It satisfies the requirement because of (\ref{eqn-local-model-functorial}) and (\ref{eqn-functorial-proj-local-model-hecke}).

\qed[Proposition-Construction \ref{propconstr-contration}]

\ssec{Proof of Theorem \ref{thm-interplay}}
\label{ssec-interplay}
\begin{goal} Consider the diagram
\begin{eqnarray*}
     H_{\MGPos}/Z_M & \xgets{\mfq^+_{P,\Vin}}  \, _\df\ol{\Bun}_{G,P} \xto{\mfp^+_{P,\Vin}} &   \ol{\Bun}_{G,\ge P}, \\
     H_{\MGPos}/Z_M  &  \xgets{\mfq^-_{P,\Vin}}   \, Y_\rel^P/Z_M\xto{\mfp^-_{P,\Vin}}   &  \ol{\Bun}_{G,\ge P}  . \\
\end{eqnarray*}
Then
	$$ \mfq^\mp_{P,\Vin,!}\circ \mfp^{\mp,*}_{P,\Vin}\simeq \mfq_{P,\Vin,*}^\pm\circ \mfp_{P,\Vin}^{\pm,!}   $$
on ind-holonomic objects.
\end{goal}

The proof is similar to that of Theorem \ref{thm-generic-second-adjointness} hence we omit some details.

Let $\gamma$ and $\ol{\gamma}$ be as in $\S$ \ref{ssec-step-2-2}. Using the homomorphism 
$$\mBG_m \xto{\gamma} Z_M \to T_\ad \xto{t\mapsto (t^{-1},t)}  T_\ad \mt  T_\ad,$$
we obtain a $\mBG_m$-action on $G\mt G$, whose attractor, repellor and fixed loci are respectively given by $P^-\mt P$, $P\mt P^-$ and $M\mt M$.

On the other hand, consider the action
$$ \mBG_m \mt \Vin_{G,\ge C_P},\;  (s, x)\mapsto \mfs( \ol{\gamma}(s)) \cdot x \cdot \mfs( \ol{\gamma}(s)).$$
This action can actually be extended to an $\mBA^1$-action using the same formula. Hence its attractor, repellor and fixed loci are respectively given by $\Vin_{G,\ge C_P}$, $\ol{M}$ and $\ol{M}$. Also, the attractor, repellor and fixed loci for the restricted action on $\Vin_{G,\gCP}^\Bru$ are respectively given by $\Vin_{G,\ge C_P}^\Bru$, $M$ and $M$.

We claim the above $\mBG_m$-actions are compatible with the action $G\mt G\act \Vin_{G,\gCP}$. Indeed, one only need to prove this claim for the restricted actions on $\Vin_{G,\gCP}\mt_{T_\ad^+} T_\ad$, which can be checked directly (see Lemma \ref{lem-compare-torus-action-pre} below). As a corollary of this claim, we obtain an action (relative to $\mBA^1$) of the Drinfeld-Gaitsgory interpolation for $G\mt G$ on that for $ \Vin_{G,\gCP}$.

Let $(\on{ActSch}_\ft^{\on{aff}})_\rel$ be the category defined similarly as $\on{ActSch}_\ft^{\on{aff}}$ (see Notation \ref{notn-group-action}) but we replace ``algebraic groups'' by ``affine group schemes over an affine base scheme''. In other words, its objects are $(H\act Y)_{/S}$, where $S$ is an affine scheme, $H\to S$ is an affine group scheme and $Y\to S$ is an affine scheme equipped with an $H$-action. There is an obvious $\affSch_\ft$-action on $(\on{ActSch}_\ft^{\on{aff}})_\rel$. By the previous discussion,
$$ (G\mt G\act  \Vin_{G,\gCP})_{/\pt}.$$
is a $\mBG_m$-module object. Then Example \ref{exam-drinfeld-pre-input-Gm} provides a weakly $\affSch_\ft$-enriched functor 
$$\Theta_{(G\mt G\act   \Vin_{G,\gCP} )}:\mbP_{\mBA^1} \to \Corr( (\on{ActSch}_\ft^{\on{aff}})_\rel )_{\all,\all},$$
sending $\alpha^+$ and $\alpha^-$ respectively to
	\begin{eqnarray*}
	  (G\mt G \act \Vin_{G,\ge C_P}) &\gets (P^- \mt P \act   \Vin_{G,\ge C_P}   ) \to& (M \mt M \act \ol{M}),\\
	   (M \mt M \act \ol{M})  &\gets (P \mt P^- \act   \ol{M}) \to& (G\mt G \act \Vin_{G,\ge C_P}).
	 \end{eqnarray*}

Passing to quotients, we obtain a weakly $\affSch_\ft$-enriched right-lax functor 
$$\Theta_{(G\backslash   \Vin_{G,\gCP}/G )}:\mbP_{\mBA^1} \laxto \bCorr( \AlgStk_\lft )_{\all,\all}^{\all,2\on{-op}}.$$
It is easy to see it is strict at the composition $\alpha^+\circ \alpha^-$. Moreover, we claim it factors through $\bCorr( \AlgStk_\lft )_{\all,\all}^{\open,2\on{-op}}$. To prove the claim, one first proves Fact \ref{fact-cartesian-quotient} below, then uses it to deduce the desired claim from Lemma \ref{lem-opennes-drinfeld-input}.

In the previous construction, we ignored the open Bruhat cell. If we keep tracking it, we would obtain a certain weakly $\affSch_\ft$-enriched right-lax functor 
$$\mbP_{\mBA^1}  \laxto  \bCorr( \on{Arr}(\AlgStk_\lft) )_{\all,\all}^{\open,2\on{-op}}.$$
By taking $\bMap_\gen(X,-)$ for it, we obtain a weakly $\affSch_\ft$-enriched right-lax functor 
$$\Theta:  \mbP_{\mBA^1} \laxto \bCorr( \PreStk_\lft )_{\all,\all}^{\open,2\on{-op}}$$
sending $\alpha^+$ and $\alpha^-$ respectively to
	\begin{eqnarray*}
	 \VinBun_{G,\gCP}  &\gets  Y^P_\rel \to& H_{\MGPos},\\
	  H_{\MGPos}  &\gets \,_\df \VinBun_{G,C_P} \to&   \VinBun_{G,\gCP}  .
	 \end{eqnarray*}
Also, $\Theta$ is strict at the composition $\alpha^+\circ \alpha^-$.

As before, we can restrict to each connected component $H_\MGPos^{\lambda,\mu}$ of $H_{\MGPos}^{\lambda,\mu}$ and obtain a Drinfeld pre-input
$$\Theta_{\lambda,\mu}:  \mbP_{\mBA^1} \laxto \bCorr( \PreStk_\lft )_{\on{safe},\on{safe}}^{\open,2\on{-op}}.$$
In fact, the right arms of the relevant correspondences are schematic.

Also, by taking quotients for the $\mBG_m$-actions, we can obtain a Drinfeld input sending $\alpha^+$ and $\alpha^-$ respectively to
	\begin{eqnarray*}
	  \VinBun_{G,\gCP}/\mBG_m  &\gets  Y^{P,\lambda,\mu}_\rel/\mBG_m \to& H_{\MGPos}^{\lambda,\mu}/\mBG_m,\\
	  H_{\MGPos}^{\lambda,\mu}/\mBG_m  &\gets \,_\df \VinBun_{G,C_P}^{\lambda,\mu}/\mBG_m \to&   \VinBun_{G,\gCP}/\mBG_m.
	 \end{eqnarray*}
By Lemma \ref{lem-compare-torus-action} below, we see that the above $\mBG_m$-action on $\VinBun_{G,\gCP}$ can be obtained from the $Z_M$-actions by restriction along $2\gamma:\mBG_m\to Z_M$. Hence Theorem \ref{thm-drinfeld-framework} implies $\mfq_{P,\Vin,*}^+\circ \mfp_{P,\Vin}^{+,!}$ is left adjoint to 
$$\prod_{\lambda,\mu} \mfp_{P,\Vin,*}^{-,\lambda,\mu}\circ \mfq_{P,\Vin}^{-,\lambda,\mu,!}.$$
Note that the above functor is also the right adjoint of $\mfq^-_{P,\Vin,!}\circ \mfp^{-,*}_{P,\Vin}$. Hence we obtain
$$  \mfq_{P,\Vin,*}^+\circ \mfp_{P,\Vin}^{+,!}\simeq \mfq^-_{P,\Vin,!}\circ \mfp^{-,*}_{P,\Vin} .$$

The equivalence $\mfq_{P,\Vin,*}^-\circ \mfp_{P,\Vin}^{-,!}\simeq \mfq^+_{P,\Vin,!}\circ \mfp^{+,*}_{P,\Vin}$
can be obtained by exchanging the roles of $\alpha^+$ and $\alpha^-$. 

\qed[Theorem \ref{thm-interplay}]

\begin{facts} \label{fact-cartesian-quotient}
For a diagram
$$ (H_1\act Y_1)_{/S_1} \to (H_2\act Y_2)_{/S_2} \gets (H_3\act Y_3)_{/S_3} $$
in $(\on{ActSch}_\ft^{\on{aff}})_\rel $, if $H_1$, $H_2$, $H_3$ and $H_1\mt_{H_2} H_3$ are all flat over their base schemes, then the following square is Cartesian

$$
\xyshort
\xymatrix{
	(Y_1\mt_{Y_2} Y_3)/(H_1\mt_{H_2} H_3) \ar[r] \ar[d] &
	(Y_1/H_1) \mt_{(Y_2/H_2)} (Y_3/H_3) \ar[d] \\
	\mBB ( H_1\mt_{H_2} H_3 ) \ar[r] & 
	\mBB H_1\mt_{\mBB H_2} \mBB H_3.
}
$$
\end{facts}

\begin{lem} \label{lem-compare-torus-action-pre}
Consider the actions
\begin{eqnarray*}
	T_\ad\act  \Vin_G, & t\cdot x :=  \mfs(t)\cdot x\cdot \mfs(t),\\
  T_\ad\act (  G\mt \Vin_G \mt G ), & t\cdot (g_1,x,g_2) :=  ( \Ad_{t^{-1}}(g_1), \mfs(t) \cdot x\cdot \mfs(t), \Ad_t(g_2) ).
 \end{eqnarray*}
The map
	$$ G\mt \Vin_G\mt G \to \Vin_G ,  (g_1,x,g_2) \mapsto g_1\cdot x\cdot g_2^{-1}$$
is equivariant for these actions.
\end{lem}

\proof We only need to prove the lemma after restricting to the subgroup of invertible elements in $\Vin_G$, which is given by $G_\enh:=(G\mt T)/Z_G$. Then we are done by a direct calculation. (Recall that the canonical section $T/Z_G\to (G\mt T)/Z_G$ is given by $t\mapsto (t^{-1},t)$).

\qed[Lemma \ref{lem-compare-torus-action-pre}]

\begin{lem} \label{lem-compare-torus-action}
Consider the following two $T$-actions on $G\backslash \Vin_{G}/G$:

\begin{itemize}
	\item[(i)] The action provided by Lemma \ref{lem-compare-torus-action-pre} via the homomorphism $T\to T_\ad$.

	\item[(ii)] The one obtained from the $T$-action on $\Vin_{G}$, which commutes with the $(G\mt G)$-action.
\end{itemize}

The action in (i) is isomorphic to the square of the action in (ii).

\end{lem}

\proof Recall that the subgroup of invertible elements in $\Vin_G$ is isomorphic to $G_\enh:=(G\mt T)/Z_G$. We have a short exact sequence $1\to G\to G_\enh \to T_\ad\to 1$. The canonical section $\mfs:T_\ad^+\to \Vin_G$ provides a splitting to the above sequence. Explicitly, this splitting is given by $t\mapsto (t^{-1},t)$. Note that the corresponding $T_\ad$ on $G$ is the \emph{inverse} of the usual adjoint action.

Consider the sequence:
$$1\to G\mt G\to G_{\enh}\mt G_{\enh} \to T_{\ad}\mt T_{\ad}\to 1.$$
Recall that the $(G\mt G)$-action on $\Vin_G$ is defined to be the restriction of the $(G_\enh \mt G_\enh)$-action on $\Vin_G$. Hence the quotient stack $G\backslash \Vin_{G}/G$ inherits a $ (T_{\ad}\mt T_{\ad}) $-action. By the last paragraph, the action in (i) is obtained from this $(T_{\ad}\mt T_{\ad})$-action by restriction along the homomorphism
\begin{equation} \label{eqn-proof-lem-compare-torus-action}
 a: T\to T_{\ad}\mt T_{\ad},\; t\mapsto (t,t^{-1}).
 \end{equation}

On the other hand, consider the certer $Z(G_{\enh}) \mt Z(G_{\enh})$ of $G_{\enh}\mt G_\enh$. Then $G_{\enh}\mt G_\enh$-action on $\Vin_{G}$ induces a $Z(G_{\enh}) \mt Z(G_{\enh})$-action on $G\backslash \Vin_{G}/G$. By construction, this action factors through the homomorphism
$$q: Z(G_{\enh}) \mt Z(G_{\enh})\to  Z(G_{\enh}),\, (s_1,s_2)\mapsto s_1s_2^{-1}.$$
In summary, we obtain compatible actions on $G\backslash \Vin_{G}/G$ by
$$ Z(G_{\enh}) \xgets{q} Z(G_{\enh}) \mt  Z(G_{\enh}) \xto{p}  T_{\ad}\mt T_{\ad},$$
where $p$ is the composition $Z(G_{\enh}) \mt  Z(G_{\enh})\to G_\enh\mt G_\enh \to T_{\ad}\mt T_{\ad}$.

Recall that the homomorphism $T\to (G\mt T)/Z_G, t\mapsto (1,t)$ provides an isomorphism between $T\simeq Z(G_\enh)$ and the $T$-action on $\Vin_G$ is defined by using this identification. Hence the square of the action in (ii) can be obtained from the $Z(G_{\enh}) \mt Z(G_{\enh})$-action via the homomorphism
$$  T \simeq Z(G_{\enh}) \xto{s\mapsto (s,s^{-1})  } Z(G_{\enh})\mt  Z(G_{\enh})  $$
(because its composition with $q$ is the square map). Then we are done because the composition of this map by $p$ is equal to $a$.

\qed[Lemma \ref{lem-compare-torus-action}]

\appendix

\section{Theory of D-modules}	\label{Appendix-corr}

	We use the theory of $2$-categories of correspondences developed in \cite[Part III]{GR-DAG1} to encode the theory of D-modules. We will use two types of this theory:

	\begin{itemize}
		\item We study all the D-modules and mainly work with the right (or standard) functors, i.e., $!$-pullbacks and $*$-pushforwards. See $\S$ \ref{sssec-standard-functors}.

		\item We study ind-holonomic D-modules and mainly work with the left functors, i.e. $*$-pullbacks and $!$-pushforwards. See $\S$ \ref{sssec-holonomic-D-modules}.
	\end{itemize}

\setcounter{subsection}{1}

\sssec{Standard functors}
\label{sssec-standard-functors}

	Consider the $(3,2)$-category
	$$ \bCorr(\PreStk_{\lft})_{\on{QCAD},\on{all}}^{\on{open}}$$
	defined as follows:
\begin{itemize}
	\item Its objects are lft prestacks;
	
	\item The $(2,1)$-category $\bMap_{\bCorr}(Y_1,Y_2)$ is the $1$-full subcategory of $(\PreStk_\lft)_{/Y_1\mt Y_2}$ where:

	\begin{itemize}
		\item we restrict to those objects $Y_2\gets Z\to Y_1$ such that $Z\to Y_2$ is \emph{QCAD};

		\item we restrict to those morphisms $Z_1\to Z_2$ in $(\PreStk_\lft)_{/Y_1\mt Y_2}$ that are schematic open embeddings;
	\end{itemize}
	
	\item the composition functor 
	$$ \bMap_{\bCorr}(Y_1,Y_2) \mt \bMap_{\bCorr}(Y_2,Y_3) \to \bMap_{\bCorr}(Y_1,Y_3) $$
	sends $Y_2\gets U \to Y_1$ and $Y_3\gets V\to Y_2$ to
	$$ Y_3\gets V\mt_{Y_2} U \to Y_1.$$
\end{itemize}

	There exists\footnote{This claim was made in \cite[Remark 9.3.13]{drinfeld2013some} even for QCA maps. A detailed construction for QCAD maps is provided in the author's thesis, see \cite[Construction C.2.13]{chen2021thesis}.} a canonical functor
	\begin{eqnarray} \label{eqn-defn-dmod-bt-!}
	\DMod^{\bt,!}:\bCorr(\PreStk_{\lft})_{\on{QCAD},\on{all}}^{\on{open},2\on{-op}} \to  \bDGCat_\cont; \\ \nonumber
	Y \mapsto \Dmod(Y),\; (Y_2\xgets{f} Z\xto{g} Y_1 ) \mapsto (f_\bt\circ g^!:\Dmod(Y_1)\to \Dmod(Y_2)). 
	\end{eqnarray}
	The content of the claim is:
\begin{itemize}
	\item for any lft prestack $Y$, there is a DG-category $\Dmod(Y)$;
	
	\item for any morphism $f:Y_1\to Y_2$, there is a \emph{$!$-pullback functor} $f^!$;
	
	\item for any QCAD morphism $f:Y_1\to Y_2$, there is a \emph{renormalized pushforward functor} $f_\bt$ defined in \cite{drinfeld2013some};
	
	\item there are base-change isomorphisms for these $!$-pullback and $\bt$-pushforward functors;
	
	\item for any schematic open embedding $f:Y_1\to Y_2$, there is an adjoint pair $(f^!,f_\bt)$;
	
	\item there are certain higher compatibilities for the above data.
\end{itemize}

	As shown in \loccit, for a \emph{safe} map $f:Y_1\to Y_2$, the renormalized pushforward functor $f_\bt$ can be identified with the usual de-Rham pushforward functor $f_*$. Therefore we keep the notation $f_*$ and only use $f_\bt$ for non-safe map $f$.

\sssec{Holonomic D-modules} \label{sssec-holonomic-D-modules}
	For any finite type affine scheme $S\in \affSch_\ft$, we write $\Dmod_\hol(S)$ for the full subcategory of $\Dmod(S)$ generated by holonomic objects (under extensions and colimits).

	For any lft prestack $Y$, we write $\Dmod_\hol(Y)$ for the full subcategory of $\Dmod(Y)$ containing objects $\mCF$ such that $f^!(\mCF) \in \Dmod_\hol(S)$ for any map $f:S\to Y$ with $S\in \affSch_\ft$. Equivalently, we define
	$$ \Dmod_\hol(Y):= \lim_{S\in (\on{AffSch}_\ft)_{/Y} } \Dmod_\hol(S),$$
	with the connecting functors given by $!$-pullbacks. An object in $\Dmod_\hol(Y)$ is called an \emph{ind-holonomic} object in $\Dmod(Y)$.

	The following results are either well-known or formal\footnote{More details are provided in the author's thesis, see \cite[$\S$ C.5]{chen2021thesis}.}:

\begin{itemize}
	
	\item[(1)] For any map $f:Y_1\to Y_2$ between lft prestacks, the functor $f^!$ preserves ind-holonomic objects. Also, the partially defined left adjoint $f_!$ of $f^!$ is well-defined on $\Dmod_\hol(Y_1)$ and sends it into $\Dmod_\hol(Y_2)$. Hence we have a functor
	\begin{equation} \label{eqn-Dhol-!-push}
	 \Dmod_\hol:  \PreStk_{\lft} \to \DGCat 
	\end{equation}
	sending morphisms to $!$-pushforward functors.

	\item[(2)] For any lft prestack $Y$, we have an equivalence
	$$  \Dmod_\hol(Y):= \colim_{S\in (\on{AffSch}_\ft)_{/Y} } \Dmod_\hol(S),$$
	with the connecting functors given by $!$-pushforwards. In particular, $\Dmod_\hol(Y)$ is compactly generated by objects of the form $g_!(\mCF)$, where $g:S\to Y$ is an object in $(\on{AffSch}_\ft)_{/Y}$ and $\mCF$ is a compact object in $\Dmod_\hol(S)$.

	\item[(3)] For any quasi-compact schematic map $f:Y_1\to Y_2$ between lft prestacks, the functor $f_*$ preserves ind-holonomic objects. Also, the partially defined left adjoint $f^*$ of $f_*$ is well defined on $\Dmod_{\hol}(Y_2)$ and sends it into $\Dmod_\hol(Y_1)$.

	\item[(4)] For any lft \emph{algebraic stacks}, there is an equivalence
	$$ \Dmod_\hol(Y):= \lim_{S\in (\on{AffSch}_\ft)_{/Y} } \Dmod_\hol(S),$$
	with the connecting functors given by $*$-pullbacks. This is implicit in \cite[$\S$ 6.2.1-6.2.2]{drinfeld2013some}.

	\item[(5)] For any map $f:Y_1\to Y_2$ between lft algebraic stacks, there is a functor 
	$$f^*: \Dmod_\hol(Y_2)\to \Dmod_\hol(Y_1)$$
	uniquely characterized by its compatibility with (4). Moreover, there exists a functor
	\begin{eqnarray} \label{eqn-Dhol-1-cat}
	\DMod_\hol^{!,*}:\Corr(\AlgStk_{\lft})_{\all,\all} \to  \DGCat_\cont,\; Y \mapsto \Dmod_\hol(Y),\\\nonumber
	(Y_2\xgets{f} Z\xto{g} Y_1 ) \mapsto (f_!\circ g^*:\Dmod_\hol(Y_1)\to \Dmod_\hol(Y_2)).
	\end{eqnarray}
	We also have its $(\infty,2)$-categorical enrichment
	\begin{equation} \label{eqn-Dhol}
	\DMod_\hol^{!,*}:\bCorr(\AlgStk_{\lft})_{\all,\all}^\open \to  \bDGCat_\cont
	\end{equation}
	obtained by using the ``no cost'' extension in \cite[Chapter 7, $\S$ 4]{GR-DAG1}.

	\item[(6)] For any stacky map $f:Y_1\to Y_2$ between lft prestacks. (2) and (5) implies there is a functor $f^*:\Dmod_\hol(Y_2)\to \Dmod_\hol(Y_1)$ equipped with base-change isomorphisms against $!$-pushforwards. In fact, by left Kan extension along
	$$ \Corr(\AlgStk_{\lft})_{\all,\all} \to \Corr(\PreStk_{\lft})_{\all,\on{Stacky}},$$
	we obtain from (\ref{eqn-Dhol-1-cat}) a functor
	$$
	\DMod_\hol^{!,*}:\Corr(\PreStk_{\lft})_{\all,\on{Stacky}} \to  \DGCat_\cont.
	$$
	It follows from (2) that its restriction on $\PreStk_{\lft} \simeq \Corr(\PreStk_{\lft})_{\all,\on{iso}}$ can be identified with (\ref{eqn-Dhol-!-push}). 
	We also have its ``no cost'' extension 
	\begin{equation} \label{eqn-Dhol-ult}
	\DMod_\hol^{!,*}:\bCorr(\PreStk_{\lft})_{\all,\on{Stacky}}^{\open} \to  \bDGCat_\cont
	\end{equation}

\end{itemize}

\section{Well-definedness results in \texorpdfstring{\cite{gaitsgory2015outline}}{[Gai15]}}
	\label{appendix-well-defined}

\setcounter{subsection}{1}

\begin{lem} \label{lem-compact-generation-IGP}
The category $\IGP$ is compactly generated by the image of compact objects under the functor $\iota_{M,!}:\Dmod(\BunM)\to \IGP$. Also, the fully faithful functor $\IGP \to \Dmod(\Bun_G^{P\hgen})$ preserves compact objects.

\end{lem}

\proof Since $\iota_M^!$ is conservative, the image of its left adjoint $\iota_{M,!}$ generates $\IGP$. Hence the first claim follows from compact generation of $\Dmod(\BunM)$. The second claim follows from the fact that $\iota_{P,!} \circ \mfq_P^*$ preserves compact objects.

\qed[Lemma \ref{lem-compact-generation-IGP}]

\begin{rem} 
S Hence $\IGP$ is compactly generated because so is $\Dmod(\BunM)$. Note that $\IGP \to \Dmod(\Bun_G^{P\hgen})$ preserves compact objects because so is $\iota_{P,!} \circ \mfq_P^*$. 
\end{rem}

\sssec{Proof of Proposition \ref{prop-well-defined-iota_M_!}}
\label{appendix-proof-prop-well-defined-iota_M_!}
Let $\wt{\Bun}_P$ be the Drinfeld's compactification constructed in \cite{braverman2002geometric}. Recall it is defined as
$$ \wt{\Bun}_P := \bMap_\gen(G\backslash \ol{G/U}/M \supset G\backslash (G/U)/M ),$$
where $\ol{G/U}$ is the affine closure of $G/U$. By \cite[Remark 4.1.9]{barlev2014d}, the map $\BunP \to \Bun_G^{P\hgen}$ factors as
$$ \BunP \xto{j} \wt{\Bun}_P \xto{\wt{\iota}_P} \Bun_G^{P\hgen},$$
and the restriction of the map $\wt{\iota}_P$ on each connected component of $\wt{\Bun}_P$ is proper. Also, the map $\wt{\iota}_P$ is obtained by applying $\bMap_\gen(X,-)$ to the morphism
$$ (G\backslash \ol{G/U}/M \supset G\backslash (G/U)/M) \to (\mBB G\gets \mBB P).$$

The above properness implies $\wt{\iota}_{P,!}$ is well-defined. On the other hand, it was proved in \cite[$\S$ 1.1.6]{drinfeld2016geometric} that the composition
$$ \Dmod( \wt{\Bun}_P ) \xto{j^!} \Dmod(\BunP) \xto{ \mfq_{P,*}}  \Dmod(\BunM) $$
has a left adjoint isomorphic to
\begin{equation}\label{eqn-ULA-wtBunP}
 j_!\circ  \mfq_{P}^*(-) \simeq j_!(k_{\BunP}) \ot^! \mfq_P^!(-) [\on{shift}],
 \end{equation}
where $[\on{shift}]$ is a cohomological shift locally constant on $\BunM$. Combining the above two results, we obtain the well-definedness of $\iota_{P,!} \circ \mfq_P^*$.

To prove the second claim, we need to calculate $\iota_{P}^!\circ \iota_{P,!} \circ \mfq_P^*$. Consider the diagram
\begin{equation}\label{eqn-proof-prop-well-defined-iota_M_!-2}
\xyshort
\xymatrix{
	\BunP \mt_{\Bun_G^{P\hgen}} \wt{\Bun}_P \ar[r]^-{\onpr_1} \ar[d]^-{\onpr_2} & 
	\BunP \ar[d]^-{\iota_P} \\
	 \wt{\Bun}_P \ar[r]^-{\wt{\iota}_P} &  \Bun_G^{P\hgen}.
}
 \end{equation}
By the base-change isomorphism, we have 
$$\iota_{P}^! \circ \wt{\iota}_{P,!} \simeq \on{pr}_{1,!}\circ \on{pr}_{2}^!.$$
A direct calculation shows
$$ \BunP \mt_{\Bun_G^{P\hgen}} \wt{\Bun}_P \simeq \bMap_\gen(X, P\backslash \ol{G/U}/M \gets P\backslash (P/U)/M ).$$
Let $\ol{M}$ be the closure of $P/U$ in $\ol{G/U}$, then we have
$$\bMap_\gen(X, P\backslash \ol{G/U}/M \gets P\backslash (P/U)/M ) \simeq  \bMap_\gen(X, P\backslash \ol{M}/M \gets P\backslash (P/U)/M ).$$
Now the RHS is isomorphic to $\BunP\mt_{\BunM} H_\MGPos$, where
$$ H_\MGPos := \bMap_\gen(X, M\backslash \ol{M}/M \supset M\backslash M/M )$$
is the $G$-positive Hecke stack for $M$-tosors (see \cite[$\S$ 3.1.5]{schieder2016geometric}). Recall that the map
$$i: \BunP\mt_{\BunM} H_\MGPos \to \wt{\Bun}_P$$
is bijective on geometric points, and the connected components of the source provide a stratification on $\wt{\Bun}_P$ (known as the \emph{defect stratification}, see \cite{braverman2002geometric}). 

Now the diagram \ref{eqn-proof-prop-well-defined-iota_M_!-2} induces
\[
\xyshort
\xymatrix{
 &  \BunP\mt_{\BunM} H_\MGPos \ar[r]^-{\onpr_1} \ar[d]^-{i} & 
 \BunP \ar[d]^-{\iota_P} \\
   \BunP\ar[r]^-j & \wt{\Bun}_P \ar[r]^-{\wt{\iota}_P} &  \Bun_G^{P\hgen}
}
\]
We obtain
$$ \iota_{P}^!\circ \iota_{P,!} \circ \mfq_P^* \simeq \iota_{P}^!\circ \wt{\iota}_{P,!} \circ j_! \circ \mfq_P^* \simeq \on{pr}_{1,!} \circ  i^! \circ j_! \circ \mfq_P^*,$$
where the last isomorphism is due to the proper base-change isomorphism. Hence it remains to show the functor $  i^! \circ j_! \circ \mfq_P^*$ factors through
\begin{equation} \label{eqn-proof-prop-well-defined-iota_M_!-1}
 \Dmod(H_{\MGPos}) \xto{*\on{-pull}}  \Dmod(\BunP\mt_{\BunM} H_\MGPos).
\end{equation}
By (\ref{eqn-ULA-wtBunP}), we only need to show $ i^! \circ j_!(k_{\BunP})$ is contained in the image of (\ref{eqn-proof-prop-well-defined-iota_M_!-1}). However, this is well-known and can be proved using the Hecke actions in \cite[$\S$ 6.2]{braverman2002geometric}.

\qed[Proposition \ref{prop-well-defined-iota_M_!}]

\sssec{Proof of Proposition \ref{prop-well-defined-Eis_enh}}
\label{appendix-proof-prop-well-defined-Eis_enh}
Let $M$ (resp. $L$) be the Levi quotient group of $P$ (resp. $Q$). Let $P_L$ be the image of $P$ in $L$, which is a parabolic subgroup of $L$. Consider the correspondence
$$ \Bun_{L} \gets \Bun_{P_{L}} \to \Bun_{M}$$
and the corresponding geometric Eisenstein series functor
$$ \Eis_{P_{L},!}:\Dmod(\Bun_{M}) \to \Dmod(\Bun_{L})$$
defined in \cite{braverman2002geometric}. Recall that it is defined as the $*$-pull-$!$-push along the above correspondence.

It is easy to check the composition\footnote{See Appendix \ref{Appendix-corr} for the definition of compositions of correspondences.} of the correspondences
\begin{eqnarray*}
	 \Bun_{L} & \gets \Bun_{P_{L}} \to&  \Bun_{M},\\
	 \Bun_G^{Q\hgen} &\xgets{\iota_Q} \Bun_Q  \xto{\mfq_Q}& \Bun_{L}
\end{eqnarray*}
is isomorphic to the composition of the correpondences
\begin{eqnarray*}
	 \Bun_G^{P\hgen} & \xgets{\iota_P}  \Bun_{P}  \xto{\mfq_P}& \Bun_{M},\\
	\Bun_G^{Q\hgen} &\xgets{\mfp_{P\to Q}^\enh} \Bun_G^{P\hgen} \xto{=}&  \Bun_G^{P\hgen}
\end{eqnarray*}
Hence by the base-change isomorphisms, we have
\begin{equation} \label{eqn-Eisenh-vs-Eis!-changing-parabolic}
 \mfp_{P\to Q,!}^\enh  \circ \iota_{P,!} \circ \mfq_P^* \simeq  \iota_{Q,!}\circ \mfq_Q^* \circ \Eis_{P_{L},!}.
\end{equation}
In particular, the LHS is well-defined. Hence by Lemma \ref{lem-compact-generation-IGP}, $\mfp_{P\to Q,!}^\enh$ is well-defined. This proves (1).

To prove (2), since $\IGP$ is compactly generated (see Lemma \ref{lem-compact-generation-IGP}), we only need to prove $\Eis_{P\to Q}^\enh$ preserve compact objects. By Lemma \ref{lem-compact-generation-IGP} again, it suffices to prove $\Eis_{P\to Q}^\enh\circ \iota_{M,!}$ preserves compact objects. By (\ref{eqn-Eisenh-vs-Eis!-changing-parabolic}), we have
\[
 \Eis_{P\to Q}^\enh\circ \iota_{M,!}\simeq \iota_{L,!}\circ \Eis_{P_{L},!}.\]
Then we are done because both $\iota_{L,!}$ and $\Eis_{P_{L},!}$ preserves compact objects.

\qed[Proposition \ref{prop-well-defined-Eis_enh}]

\section{D-modules on stacks stratified by power sets}
\label{appendix-stratification}

\setcounter{subsection}{1}
We begin with the following definition.

\begin{defn} Let $Y$ be an algebraic stack stratified by a power poset $P(I)$. We define
$$ \Funct(P(I), \Dmod_\hol(Y))_{!} \subset \Funct(P(I), \Dmod_\hol(Y)) $$
to be the full subcategory consisting of those functors $F:P(I) \to \Dmod(Y)$ such that $F(J)$ is $!$-extended from the stratum $Y_J$.
\end{defn}

\begin{lem} \label{lem-gluing-functors}
Let $Y$ be an algebraic stack stratified by a power poset $P(I)$. The functor
	\begin{eqnarray*}
	\mbC_Y: \Funct(P(I), \Dmod_\hol(Y))_{!} \to \Dmod_\hol(Y),\\
	F \mapsto \on{coFib}( \colim_{J\subsetneq I} F(J) \to F(I) )  
	\end{eqnarray*}
	is an equivalence. Also, its inverse sends an object $\mCF\in \Dmod_\hol(Y)$ to a certain functor
	$$  P(I) \to \Dmod_\hol(Y),\; J\mapsto i_{J,!}\circ i^*_{J}(\mCF)[|J|-|I|].$$

\end{lem}

\proof We first recall the following formal fact. et $[1]$ be the poset $\{0,1\}$ and $I$ be an index category obtained by removing the final object from $[1]^r$ ($r\ge 1$). Let $C$ be any stable category. Suppose $F:I\to C$ is a functor such that $F(x)\simeq 0$ unless $x$ is the initial object $i_0$. Then $\colim F \simeq F(i_0)[r-1]$. This fact can be proven by induction on $r$. 

It follows that
$$i_K^*(\on{coFib}( \colim_{J\subsetneq I} F(J) \to F(I) )) \simeq i_K^*(F(K))[|I|-|K|]. $$
Hence the second claim follows from the first one.

It remains to show $\mbC_Y$ is an equivalence. 

The case $I=\{*\}$ is well-known. Namely, for any open substack $j:U\to Y$ and its complement $i: Z\to Y$, we have $j_! \circ j^! \simeq \on{Fib}( \Id\to i_!\circ i^*)$ and therefore the two functors in the statement of lemma are inverse to each other.

The general case can be proved by induction as follows. Suppose $I = I^\flat \bigsqcup \{a\}$ and $I$ is nonempty. Note that apart from the embedding $P(I^\flat) \subset P(I)$, we also have a map
$$ P(I^\flat) \to P(I),\; J\mapsto J^\sharp:=J\bigsqcup \{a\} .$$
The open substacks $\{U_i\}_{i\in I^\flat}$ provide a stratification of $Y$ labelled by $P(I^\flat)$. We use the notation $Z$ to denote the same stack $Y$ equipped with this new stratification. For any $K\in P(I^\flat)$, the stratum $Z_K$ inherits a stratification by $P(\{a\})$. The stratum $(Z_K)_{\{a\}}$ labelled by $\{a\}$ is isomorphic to $Y_{K^\sharp}$ and the stratum $(Z_K)_{\{\emptyset\}}$ labelled by $\emptyset$ is isomorphic to $Y_K$. Consider the functor
$$A: \Funct(P(I), \Dmod_\hol(Y))_! \to \Funct(P(I^\flat), \Dmod_\hol(Z))_! ,F\mapsto A(F), $$
where $A(F)(K):= \on{coFib}( F(K) \to F(K^\sharp) )$. Note that this is well-defined, i.e. $A(F)(K)$ is indeed a $!$-extension from $Z_K$. Moreover, $A$ is an equivalence by the $I=\{*\}$ case of the lemma (applying to each $Z_K$).

Hence by induction hypothesis, $\mbC_Z\circ A$ is also an equivalence. It remains to show $\mbC_Y \simeq \mbC_Z\circ A$. Note that we have the following pushout diagram
\[
\xyshort
\xymatrix{
	\colim_{K\subsetneq I^\flat} F(K)  \ar[r] \ar[d] &
	\colim_{J\subsetneq I, J\ne I^\flat} F(J) \ar[d]  \\
	\colim_{K\subset I^\flat} F(K) \ar[r] & \colim_{J\subsetneq I} F(J),
}
 \]
which is obtained by writing the simplicial nerve of $P(I)\setminus\{I\}$ as a pushout. By cofinality, the above diagram is equivalent to
\[
\xyshort
\xymatrix{
	\colim_{K\subsetneq I^\flat} F(K)  \ar[r] \ar[d] &
	\colim_{K\subsetneq I^\flat} F(K^\sharp) \ar[d]  \\
	F(I^\flat) \ar[r] & \colim_{J\subsetneq I} F(J).
}
 \]
Then we have
\begin{eqnarray*}
 & & \on{coFib}( \colim_{J\subsetneq I} F(J) \to F(I) ) \\
 & \simeq & \on{coFib}(  \on{coFib}(\colim_{K\subsetneq I^\flat} F(K) \to \colim_{K\subsetneq I^\flat} F(K^\sharp))   \to \on{coFib} (  F(I^\flat) \to F(I)  )  )   \\
 & \simeq & \on{coFib}(  \colim_{K\subsetneq I^\flat} ( \on{coFib}(F(K)\to F(K^\sharp)  )  )   \to \on{coFib} (  F(I^\flat) \to F(I)  )  )  \\
 & \simeq & \on{coFib}(  \colim_{K\subsetneq I^\flat} A(F)(K)   \to A(F)(I^\flat)  )
\end{eqnarray*} 
as desired. This proves the claim.

\qed[Lemma \ref{lem-gluing-functors}]

The above lemma implies

\begin{cor} \label{Cor-gluing-functors}
Let $Y$ be an algebraic stack stratified by a power poset $P(I)$. The functor
	$$
	\mbJ_Y: \Funct(P(I), \Dmod_\hol(Y))_{!} \to \Dmod_\hol(Y_I),\; F \mapsto j_I^*\circ F(I) $$
	has a right adjoint sending an object $\mCF\in \Dmod_\hol(Y_I)$ to a certain functor
	$$ \mbG^*_{\mCF,Y}: P(I) \to \Dmod_\hol(Y),\; J\mapsto i_{J,!}\circ i^*_{J}\circ j_{I,*} (\mCF)[|J|-|I|].$$

\end{cor}

\proof Follows from the fact that $\mbJ_Y \simeq j_I^*\circ \mbC_Y.$

\qed[Corollary \ref{Cor-gluing-functors}]

\begin{rem}  \label{rem-compare-gluing-functor}
	Note that the functor $\mbG^*_{\mCF,Y}$ sends the arrow $J\subset I$ to a morphism
$$ i_{J,!}\circ i^*_{J}\circ j_{I,*} (\mCF)[|J|-|I|] \to i_{I,!}(\mCF).$$
Applying $i_J^!$ to this map, we obtain a map
$$  i^*_{J}\circ j_{I,*} (\mCF)[|J|-|I|] \to i_{J}^!\circ i_{I,!}(\mCF).$$
Note that this map is invertible if $|J|=|I|-1$, but not for general $J$.
\end{rem}

\begin{lem}   \label{lem-gluing-functors-colimite}
Let $Y$ be an algebraic stack stratified by a power poset $P(I)$ and $J\in P(I)$. Consider the maps
$$ Y_I \xto{j_{I,\ge J} } \to Y_{\ge J} \xto{j_{\ge J}}\to Y.$$
For any $\mCF\in \Dmod_\hol(Y_I)$, we have
$$  \on{coFib}(\colim_{ J\subset K\subsetneq I } \mbG^*_{\mCF,Y}(K) \to \mbG^*_{\mCF,Y}(I)    )    \simeq    (j_{\ge J})_!  \circ  (j_{I,\ge J})_* (\mCF). $$
\end{lem}

\proof The case $J=I$ follows from definition. Indeed, the LHS is given by 
$$\mbC_Y\circ (\mbJ_Y)^R \simeq \mbC_Y\circ ( j_I^*\circ \mbC_Y)^R \simeq  (j_I^*)^R \simeq j_{I,*}.$$

In the general case, note that both sides are contained in the image of the functor $(j_{\ge J})_!$. Hence we only need to show
$$  \on{coFib}(\colim_{J\subset K\subsetneq I } j^*_{\ge J}\circ \mbG^*_{\mCF,Y}(K) \to  j^*_{\ge J}\circ \mbG^*_{\mCF,Y}(I)    )    \simeq   (j_{I,\ge J})_* (\mCF). $$
Consider the open substack $Y_{\ge J}$. It inherits a stratification by the poset $P(I\setminus J)$ with $(Y_{\ge J})_K \simeq Y_{J\bigsqcup K}$. Hence we also have a functor
$$  \mbG^*_{\mCF,Y_{\ge J}}:  P(I\setminus J) \to \Dmod_\hol(Y_{\ge J}).  $$
It follows from construction that this functor is isomorphic to
$$  P(I\setminus J)  \xto{-\bigsqcup J}\to P(I) \xto{ \mbG^*_{\mCF,Y} }\to   \Dmod_\hol(Y) \xto{  j^*_{\ge J} }  \to \Dmod_\hol(Y_{\ge J}).$$
Hence we only need to show
$$  \on{coFib}(\colim_{ K\subsetneq I\setminus J }\mbG^*_{\mCF,Y_{\ge J}}(K) \to  \mbG^*_{\mCF,Y_{\ge J}}( I\setminus J  )    )    \simeq   (j_{I,\ge J})_* (\mCF). $$
In other words, we have reduced the lemma to the case $J=I$. 

\qed[Lemma \ref{lem-gluing-functors-colimite}]

\section{The group scheme \texorpdfstring{$\wt{G}$}{tilta G}}\label{appendix-wtG}
	Consider the $(G\mt G)$-action on $\Vin_G$. Note that it preserves the fibers of $\Vin_G\to T_\ad^+$. We write $\wt{G}$ for the corresponding stabilizer of the canonical section $\mfs:T_\ad^+ \to \Vin_G$. In this appendix, we review some results about $\wt{G}$.

	\setcounter{subsection}{1}

	We begin by reviewing some facts:
	\begin{facts} \label{fact-wtg}
	We have the following facts\footnote{(1) and (2) are well-known, (4) and (5) follow from the identification $\Vin_G\mt_{T_\ad^+} T_\ad \simeq (G\mt T)/Z_G$.}: 
		
	\begin{itemize}
		\item[(1)] $\wtG$ is a closed subgroup of $G\mt G \mt T_\ad^+$ (relative to $T_\ad^+$), whose fiber at $C_P$ is
			$$\wt{G}_{C_P}\simeq P\mt_{M} P^-.$$

		\item[(2)] By \cite[Corollary D.5.4]{drinfeld2016geometric}, $\wtG$ is smooth over $T_\ad^+$, and we have 
			\begin{equation} \label{eqn-defect-free-vin-as-quotient}
 				_0\!\Vin_G \simeq (G\mt G \mt T_\ad^+)/\wt{G},\;  G\backslash \,_0\!\Vin_G/G \simeq \mBB \wt{G},
 			\end{equation}
			where $\mBB$ means taking relative classifying stack.

		\item[(3)] By (2), the $T$-action on $\Vin_G$ (which commutes with the $G\mt G$-action) induces a diagram between group actions:
		$$( \pt \act \mBB G\mt \mBB G )   \gets (T \act \mBB \wtG) \to (T\act T_\ad^+).$$

			\item[(4)] $\wt{G}$ contains the locally closed subscheme
		$$\Gamma: G\mt T_\ad \to G\mt G\mt T_\ad^+,\; (g,t)\mapsto ( g, \on{Ad}_t(g),t) $$

		\item[(5)] $\wt{G}$ is preserved by the action
			$$ (T_\ad\mt T_\ad) \act (G\mt G\mt T_\ad^+),\; (t_1,t_2)\cdot(g_1,g_2,s)\mapsto (\on{Ad}_{t_1^{-1}}(g_1), \on{Ad}_{t_2^{-1}}(g_2),  t_1\cdot s \cdot t_2^{-1}).$$
	\end{itemize}
\end{facts}

\begin{warn} \label{warn-torus-action}
	The $T$-action on $\Vin_G$ does not induce a $T$-action on $\wt{G}$ because this action does not preserve the canonical section $\mfs:T_\ad^+ \to \wt{G}$.
\end{warn}

\begin{rem}
In the rest of this appendix, we prove some new results about $\wt{G}$. If we replace $\wt{G}$ by the fibers of the map $\wt{G}\to T_\ad^+$, then all these results become obvious because of the explicit description of the fibers $\wt{G}_{C_P}$ in Fact \ref{fact-wtg}(1). Hence we focus on how to make the arguments work in family (over $T_\ad^+$).
\end{rem}

The following result generalizes \cite[Proposition D.6.4]{drinfeld2016geometric}:

\begin{lem} \label{lem-wtG-closure}
	$\wt{G}$ is isomorphic to the closure of the locally closed embedding
	$$\Gamma:G\mt T_\ad \to G\mt G\mt T_\ad^+,\; (g,t)\mapsto ( g, \on{Ad}_t(g),t). $$
\end{lem}

\proof 
	Let $\ol{\Gamma}$ be the desired closure. Hence we obtain a closed embedding $\ol{\Gamma} \to \wt{G}$. Since $\wt{G}$ is reduced, it remains to show $\ol{\Gamma} \to \wt{G}$ is surjective. Note that $\ol{\Gamma}$ is also preserved by the action in Fact \ref{fact-wtg}(5). Hence we only need to check the fiber of $\ol{\Gamma} \to \wt{G}$ at $C_P\in T_\ad^+$ is surjective. Then we are done by \cite[Proposition D.6.4]{drinfeld2016geometric}.

	\qed[Lemma \ref{lem-wtG-closure}]

\begin{lem} \label{lem-wtG-leP-containing}
The closed subscheme
 $$\wt G_{\le P}:= \wtG \mt_{\Tadp} T_{\ad,\le P}^+ \inj G\mt G \mt  T_{\ad,\le P}^+   $$
 is contained in $P\mt P^- \mt T_{\ad,\le P}^+$.
\end{lem}

\proof Using the action in Fact \ref{fact-wtg}(5), we only need to show $\wtG_{C_Q}$
is contained in $P\mt P^-$ for any $Q\subset P$. But this is obvious.

\qed[Lemma \ref{lem-wtG-leP-containing}]

\begin{lem} \label{lem-Pm-wtG-gCP} We write $\wtG_\gCP:= \wt{G} \mt_{T_\ad^+} T_{\ad,\gCP}^+$. We have:

	\begin{itemize}
	\item[(1)] The closed subscheme
	$$  P^-\mt_G \wt{G}_\gCP \inj P^-\mt G\mt T_{\ad,\gCP}^+ $$
	is contained in $P^-\mt P^-\mt T_{\ad,\gCP}^+$.

	\item[(2)] The composition
	\begin{equation} \label{eqn-lem-Pm-wtG-gCP}
	 P^-\mt_G \wt{G} \mt_{T_\ad^+} T_{\ad,\gCP}^+ \to P^-\mt P^-\mt T_{\ad,\gCP}^+ \xto{\on{pr}_{23}}  P^-\mt T_{\ad,\gCP}^+ 
	 \end{equation}
	is an isomorphism, where the first map is obtained by (1).
	\end{itemize}
\end{lem}

\begin{warn} The similar statement for $\on{pr}_{13}$ is false.
\end{warn}

\proof 
	We first prove (1). Using the action in Fact \ref{fact-wtg}(5), we only need to check the similar claim at any $C_{P'}\in T_{\ad,\gCP}^+$. But this is obvious.

	Similarly, it is easy to see (\ref{eqn-lem-Pm-wtG-gCP}) induces isomorphisms between fibers at any closed point of $T_{\ad,\gCP}^+$. To prove (2), we only need to show $P^-\mt_G \wt{G} \mt_{T_\ad^+} T_{\ad,\gCP}^+ $ is smooth over $T_{\ad,\gCP}^+$.

	We claim $P^-\mt G\mt T_{\ad,\gCP}^+$ and $\wtG_\gCP$ are transversal in $G\mt G\mt T_{\ad,\gCP}^+$. Indeed, by the last paragraph, the dimension of any irreducible component of their intersection is at most $\dim(P^-)+\dim( T_{\ad,\gCP}^+)$. But this number is equal to 
	$$ \dim( P^-\mt G\mt T_{\ad,\gCP}^+ ) + \dim (\wtG_\gCP)- \dim(G\mt G\mt T_{\ad,\gCP}^+).$$
	This proves the transversity. In particular, we obtain that $P^-\mt_G \wt{G} \mt_{T_\ad^+} T_{\ad,\gCP}^+$ is smooth.

	It remains to show $f:P^-\mt_G \wt{G} \mt_{T_\ad^+} T_{\ad,\gCP}^+\to T_{\ad,\gCP}^+$ induces surjections between tangent spaces. Note that the fibers of this map is smooth and of dimension $\dim(P^-)$. Hence at any closed point $x$ of the source, we have
	$$\dim(\on{ker}(df_x))= \dim(P^-) =\dim(  P^-\mt_G \wt{G} \mt_{T_\ad^+} T_{\ad,\gCP}^+ ) - \dim (  T_{\ad,\gCP}^+  ).$$
	This implies $df_x$ is surjective.

\qed[Lemma \ref{lem-Pm-wtG-gCP}]

\begin{lem} \label{lem-Pm-s-G-open-orbit}
Consider the $(P^-\mt G)$-action on $\Vin_{G,\ge C_P}$. Its stablizer for the canonical section is $P^-\mt_G \wt{G}_\gCP$. Then the map
\begin{equation} \label{proof-lem-Pm-s-G-open-orbit-2}
(P^-\mt G\mt T_{\ad,\gCP}^+)/( P^-\mt_G \wtG_{\gCP} ) \to \,_0\!\Vin_{G,\gCP} 
\end{equation}
induced by this action is an open embedding.
\end{lem}

\proof 
	We claim the LHS is a smooth scheme. By Lemma \ref{lem-Pm-wtG-gCP}, there is an isomorphism 
	\begin{equation} \label{proof-lem-Pm-s-G-open-orbit-1}
	P^-\mt_G \wtG_{\gCP} \simeq P^-\mt T_{\ad,\gCP}^+
	\end{equation}
	between group schemes over $T_{\ad,\gCP}^+$. Moreover, the projection map $P^-\mt G\mt T_{\ad,\gCP}^+ \to G\mt T_{\ad,\gCP}^+$ intertwines the actions of (\ref{proof-lem-Pm-s-G-open-orbit-1}). Hence we obtain a map
	$$ (P^-\mt G\mt T_{\ad,\gCP}^+)/( P^-\mt_G \wtG_{\gCP} ) \to  ( G\mt T_{\ad,\gCP}^+)/(P^-\mt T_{\ad,\gCP}^+) \simeq G/P^- \mt T_{\ad,\gCP}^+ .$$
	Since $P^-\mt G\mt T_{\ad,\gCP}^+\to  G\mt T_{\ad,\gCP}^+$ is affine and smooth, the above map is also affine and smooth. This proves the claim on smoothness. Then the lemma follows from the fact that both sides of (\ref{proof-lem-Pm-s-G-open-orbit-2}) have the same dimenison and that this map is injective on the level of closed points.

\qed[Lemma \ref{lem-Pm-s-G-open-orbit}]

\begin{cor} \label{cor-open-BPmwtG}
The map
$$  \mBB(  P^-\mt_G \wtG_{\gCP} ) \to \mBB P^-\mt_{\mBB G}\mBB \wtG_\gCP  $$ is a schematic open embedding
\end{cor}

\proof Follows from Lemma \ref{lem-Pm-s-G-open-orbit} by taking quotients for the $(P^-\mt G)$-actions.

\qed[Corollary \ref{cor-open-BPmwtG}]

\bibliography{mybiblio.bib}{}
\bibliographystyle{alpha}

\end{document}